\def\version{}

\documentclass[reqno,twoside,11pt]{amsart}
\usepackage{cite}
\usepackage{amsmath}
\usepackage{amsfonts}
\usepackage{amssymb}
\usepackage{epsfig}
\usepackage{verbatim}
\usepackage{color}

\IfFileExists{epsf.def}{\input epsf.def}{\usepackage{epsf}}
\IfFileExists{mathptmx.sty}{\usepackage{mathptmx}}{\usepackage{mathpazo}}
\usepackage{mathrsfs}

\DeclareFontFamily{OT1}{eusb}{} \DeclareFontShape{OT1}{eusb}{m}{n} {<5> <6> <7> <8> <9> <10> <11> <12> <14.4> eusb10}{}
\DeclareMathAlphabet{\eusb}{OT1}{eusb}{m}{n}

\DeclareFontFamily{OT1}{eusm}{} \DeclareFontShape{OT1}{eusm}{m}{n} {<5> <6> <7> <8> <9> <10> <11> <12> <14.4> eusm10}{}
\DeclareMathAlphabet{\eusm}{OT1}{eusm}{m}{n}

\DeclareFontFamily{OT1}{eufm}{} \DeclareFontShape{OT1}{eufm}{m}{n} {<5> <6> <7> <8> <9> <10> <11> <12> <14.4> eufm10}{}
\DeclareMathAlphabet{\mathfrak}{OT1}{eufm}{m}{n}

\DeclareFontFamily{OT1}{fraktura}{}
\DeclareFontShape{OT1}{fraktura}{m}{n} {<5> <6> <7> <8> <9> <10> <11> <12> <13> <14.4> [1.1] eufm10}{}
\DeclareMathAlphabet{\fraktura}{OT1}{fraktura}{m}{n}

\DeclareFontFamily{OT1}{cmfi}{} \DeclareFontShape{OT1}{cmfi}{m}{n} {<5> <6> <7> <8> <9> <10> <11> <12> <13> <14.4> [0.9] cmfi10}{}
\DeclareMathAlphabet{\cmfi}{OT1}{cmfi}{b}{n}

\DeclareFontFamily{OT1}{cmss}{} \DeclareFontShape{OT1}{cmss}{m}{n} {<5> <6> <7> <8> <9> <10> <11> <12> <13> <14.4> cmss10}{}
\DeclareMathAlphabet{\cmss}{OT1}{cmss}{m}{n}

\newtheoremstyle{thm}{1.5ex}{1.5ex}{\itshape\rmfamily}{} {\bfseries\rmfamily}{}{2ex}{}

\newtheoremstyle{def}{1.5ex}{1.5ex}{\rmfamily\sl}{} {\bfseries\rmfamily}{}{2ex}{}

\newtheoremstyle{rem}{1.3ex}{1.3ex}{\rmfamily}{} {\bfseries\rmfamily}{}{2ex}{}

\newtheoremstyle{ass}{1.5ex}{1.5ex}{\rmfamily\sl}{} {\bfseries\rmfamily}{}{2ex}{}

\newenvironment{proofsect}[1] {\vskip0.1cm\noindent{\rmfamily\itshape#1.}}{\qed\vspace{0.15cm}}

\theoremstyle{thm}
\newtheorem{theorem}{Theorem}[section]
\newtheorem{lemma}[theorem]{Lemma}
\newtheorem{proposition}[theorem]{Proposition}

\newtheorem*{Main Theorem}{Main Theorem.}
\newtheorem{corollary}[theorem]{Corollary}

\newtheorem{problem}[theorem]{Problem}

\newtheoremstyle{named}{}{}{\itshape}{}{\bfseries}{}{.5em}{\thmnote{#3}}
\theoremstyle{named}
\newtheorem*{namedtheorem}{Theorem}

\theoremstyle{def}

\theoremstyle{rem}
\newtheorem{remark}[theorem]{{Remark}}

\numberwithin{equation}{section}
\renewcommand{\theequation}{\arabic{section}.\arabic{equation}}

\renewcommand{\section}{\secdef\sct\sect}
\newcommand{\sct}[2][default]{\refstepcounter{section}
\addcontentsline{toc}{section}
{{\tocsection {}{\thesection}{\!\!\!\!#1\dotfill}}{}}
\vspace{0.7cm}
\centerline{ 
\scshape\arabic{section}.\ #1} \nopagebreak \vspace{0.2cm}}
\newcommand{\sect}[1]{
\vspace{0.4cm} \centerline{\large\scshape\rmfamily #1}
\vspace{0.2cm}}

\renewcommand{\subsection}{\secdef\subsct\sbsect}
\newcommand{\subsct}[2][default]{\refstepcounter{subsection}
\addcontentsline{toc}{subsection}
{{\tocsection{\!\!}{\hspace{1.2em}\thesubsection}{\!\!\!\!#1\dotfill}}{}}
\nopagebreak\vspace{0.45\baselineskip} {\flushleft\bf
\thesection.\arabic{subsection}~\bf #1.~}
\\*[3mm]\noindent
\nopagebreak}
\newcommand{\sbsect}[1]{
\vspace{0.1cm}\noindent
\textbf{#1.~}\vspace{0.1cm}}

\renewcommand{\subsubsection}{%
\secdef \subsubsect\sbsbsect}
\newcommand{\subsubsect}[2][default]{%
\refstepcounter{subsubsection} 
\addcontentsline{toc}{subsubsection}{{\tocsection{\!\!}
{\hspace{3.05em}\thesubsubsection}{\!\!\!\!#1\dotfill}}{}}
\nopagebreak
\vspace{0.15\baselineskip} \nopagebreak {\flushleft\rmfamily
\itshape\arabic{section}.\arabic{subsection}.\arabic{subsubsection}
\ \rmfamily #1\/.}\ }
\newcommand{\sbsbsect}[1]{\vspace{0.1cm}\noindent
\rmfamily \itshape
\arabic{section}.\arabic{subsection}.\arabic{subsubsection} \
\sffamily #1\/.\ }

\renewcommand{\caption}[1]{%
\vglue0.5cm
\refstepcounter{figure}
\begin{center}
\begin{minipage}[c]{0.8\textwidth}\small {\sc Fig.~\thefigure\ }#1\end{minipage}
\end{center}
}

\newcommand{\dist}{\operatorname{dist}}

\newcommand{\diam}{\operatorname{diam}}

\newcommand{\textd}{\text{\rm d}\mkern0.5mu}

\newcommand{\texte}{\text{\rm  e}\mkern0.7mu}

\newcommand{\Var}{\text{\rm Var}}
\newcommand{\1}{{1\mkern-4.5mu\textrm{l}}}
\renewcommand{\1}{\text{\sf 1}}

\newcommand{\FF}{\mathcal F}

\newcommand{\NN}{\mathcal N}

\newcommand{\C}{\mathbb C}

\newcommand{\E}{\mathbb E}

\newcommand{\N}{\mathbb N}

\newcommand{\Q}{\mathbb Q}
\newcommand{\R}{\mathbb R}

\newcommand{\Z}{\mathbb Z}

\newcommand{\scrX}{\mathscr{X}}

\newcommand{\twoeqref}[2]{(\ref{#1}--\ref{#2})}

\newcommand{\cc}{{\text{\rm c}}}

\newcommand{\frakg}{\fraktura g}
\newcommand{\fraka}{\fraktura a}

\def\myffrac#1#2 in #3{\raise 2.6pt\hbox{$#3 #1$}\mkern-1.5mu\raise 0.8pt\hbox{$#3/$}\mkern-1.1mu\lower 1.5pt\hbox{$#3 #2$}}
\newcommand{\ffrac}[2]{\mathchoice%
	{\myffrac{#1}{#2} in \scriptstyle}
	{\myffrac{#1}{#2} in \scriptstyle}
	{\myffrac{#1}{#2} in \scriptscriptstyle}
	{\myffrac{#1}{#2} in \scriptscriptstyle}
}

\newcommand{\wh}{\widehat}
\newcommand{\wt}{\widetilde}
\newcommand{\ol}{\overline}

\newcommand{\laweq}{\,\overset{\text{\rm law}}=\,}

\newcommand{\leb}{{\rm Leb}}

\newcommand{\Lawarrow}{{\,\overset{\text{\rm law}}\longrightarrow\,}}

\newcommand{\DGFF}{{\rm DGFF}}

\newcommand{\MB}{}
\newcommand{\eMB}{\normalcolor}

\newcommand{\OL}{}
\newcommand{\eOL}{\normalcolor}

\begin{document}

\title[Intermediate level sets of 2D \DGFF{} \hfill \version\hfill]
{\large On intermediate level sets of two-dimensional\\discrete Gaussian Free Field}

\author[\hfill  \version \hfill Biskup and Louidor]
{Marek~Biskup and Oren~Louidor}
\thanks{\hglue-4.5mm\fontsize{9.6}{9.6}\selectfont\copyright\,\textrm{2016}\ \ \textrm{M.~Biskup, O.~Louidor.
Reproduction, by any means, of the entire
article for non-commercial purposes is permitted without charge.\vspace{2mm}}}
\maketitle

\vspace{-5mm}
\centerline{\textit{
Department of Mathematics, UCLA, Los Angeles, California, USA}}
\centerline{\textit{
Faculty of Industrial Engineering and Management, Technion, Haifa, Israel}}

\vskip0.5cm
\begin{quote}
\footnotesize \textbf{Abstract:}
We consider the discrete Gaussian Free Field (\DGFF) in scaled-up (square-lattice) versions of suitably regular continuum domains~$D\subset\C$ and describe the scaling limit, including local structure, of the level sets at heights growing as a $\lambda$-multiple of the height of the absolute maximum, for any $\lambda\in(0,1)$. We prove that, in the scaling limit, the scaled spatial position of a typical point~$x$ sampled from this level set is distributed according to a Liouville Quantum Gravity (LQG) measure in~$D$  at parameter equal~$\lambda$-times its critical value,  the field value at~$x$ has an exponential intensity measure and the configuration near~$x$ reduced by the value at~$x$ has the law of a pinned \DGFF{} reduced by a suitable multiple of the  potential  kernel.  In particular, the law of the total size of the level set, properly-normalized, converges that that of the total mass of the LQG measure. This sharpens considerably an earlier conclusion by Daviaud~\cite{Daviaud}.
\end{quote}

\section{Introduction}
\noindent
It has long been recognized that the two-dimensional \emph{continuum} Gaussian Free Field (CGFF) offers a variety of constructions of random fractals with an underlying conformally-invariant structure. This has been used fruitfully in the work of Schramm and Shef\-field~\cite{Schramm-Sheffield} on the convergence to~SLE$_4$ of the level lines at specific heights of order unity and the ensuing coupling of the whole field to the Conformal Loop Ensemble by Sheffield and Werner~\cite{Sheffield-CLE,Sheffield-Werner}. Other examples include
the construction of the Liouville Quantum Gravity measures by Duplantier and Sheffield~\cite{Duplantier-Sheffield} as well as the recent research programs of Miller and Sheffield on imaginary geometry~\cite{MS1,MS2,MS3,MS4} and the connection between the Liouville Quantum Gravity and the Brownian Map~\cite{MS5,MS6}.

A parallel, and largely independent, line of recent research has focused on various quantitative aspects of the extremal values associated with the \emph{discrete} Gaussian Free Field (\DGFF). This is a Gaussian process $\{h_x\colon x\in\Z^2\}$ marked by a proper (typically finite) subset~$V$ of the square lattice (other infinite graphs can be considered as well) with the law determined by 
\begin{equation}
\label{E:1.1}
E(h_x)=0\quad\text{and}\quad E(h_x\,h_y)=G^V(x,y),
\end{equation}
where~$G^V$ denotes the Green function of the simple symmetric random walk in~$V$ killed upon exit from~$V$. (In particular, $h$ vanishes outside of~$V$ almost surely.) Here an early paper  of  Bolthausen, Deuschel and Giacomin~\cite{BDeuG} showed that the maximum of the \DGFF{} in square boxes~$V_N:=(0,N)^2\cap\Z^2$ grows as
\begin{equation}
\label{E:1.2}
\max_{x\in V_N}h_x\,\sim\, 2\sqrt g\,\log N,\qquad N\to\infty,
\end{equation}
where ``$\sim$'' designates that the ratio of the two quantities tends to one in the stated limit and $g:=2/\pi$ is a constant such that the Green function obeys $G^{V_N}(x,x)=g\log N+O(1)$ as~$N\to\infty$ for~$x$ ``deep'' inside~$V_N$. Daviaud~\cite{Daviaud} subsequently  finessed the approach of~\cite{BDeuG}  to capture some geometric aspects of the \emph{intermediate} level sets
\begin{equation}
\label{E:1.3}
\bigl\{x\in V_N\colon h_x\ge 2\sqrt g\,\lambda\log N\bigr\}
\quad\text{for}\quad\lambda\in(0,1).
\end{equation}
He showed that this set contains $N^{2(1-\lambda^2)+o(1)}$  points, where $o(1)\to0$ in probability as $N\to\infty$, and thus demonstrated a fractal nature of this set. (The structure of the exponent is quite universal; see Chatterjee, Dembo and Ding~\cite{Chatterjee-Dembo-Ding}. Continuum versions of this result exist, dealing with thick points of CGFF; cf Hu, Peres and Miller~\cite{HMP10}.) 
Other fractal properties were also proved; e.g.\ the growth-rate exponents of its intersection with balls of increasing radii.

\begin{figure}[t]
\centerline{\includegraphics[width=0.3\textwidth]{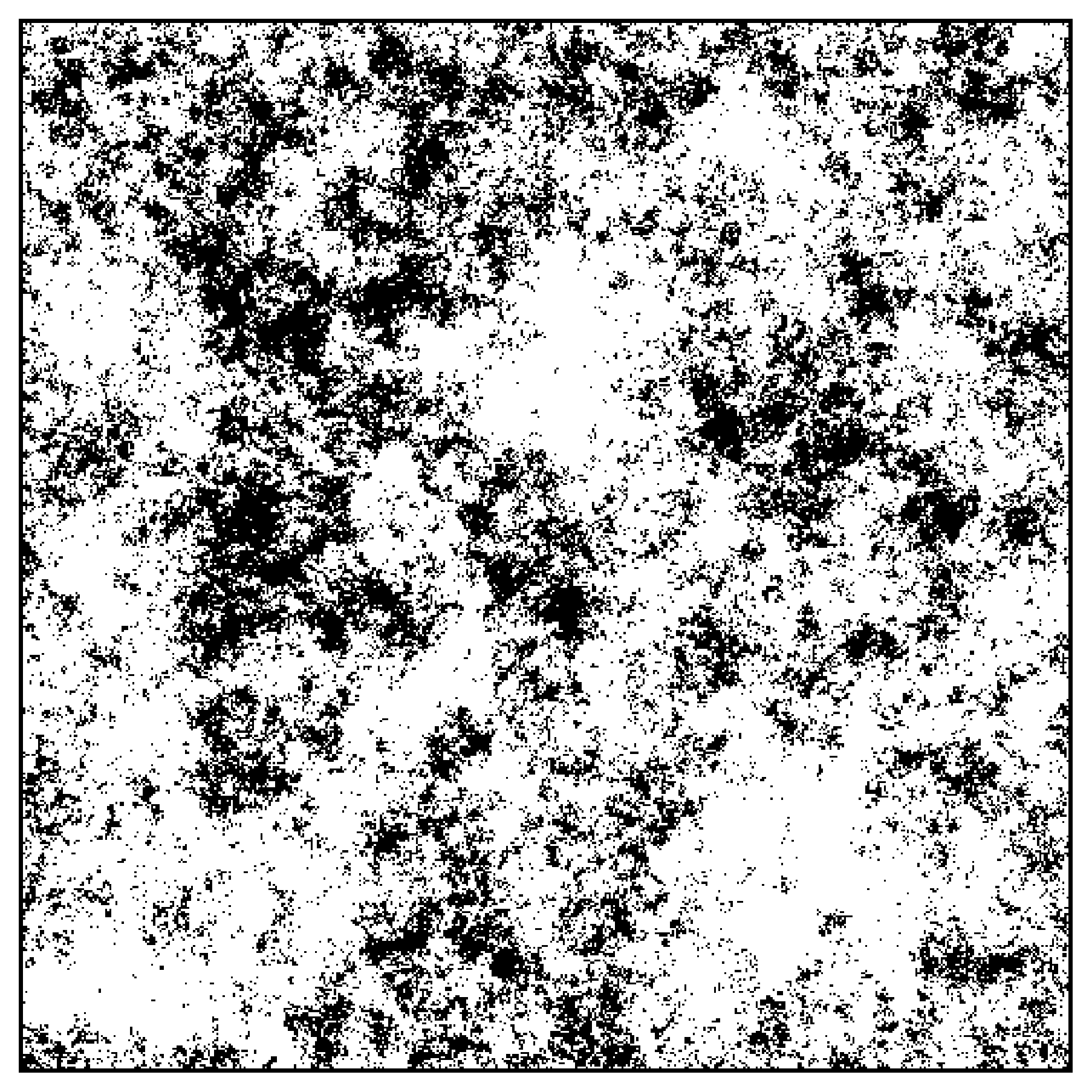}
\includegraphics[width=0.3\textwidth]{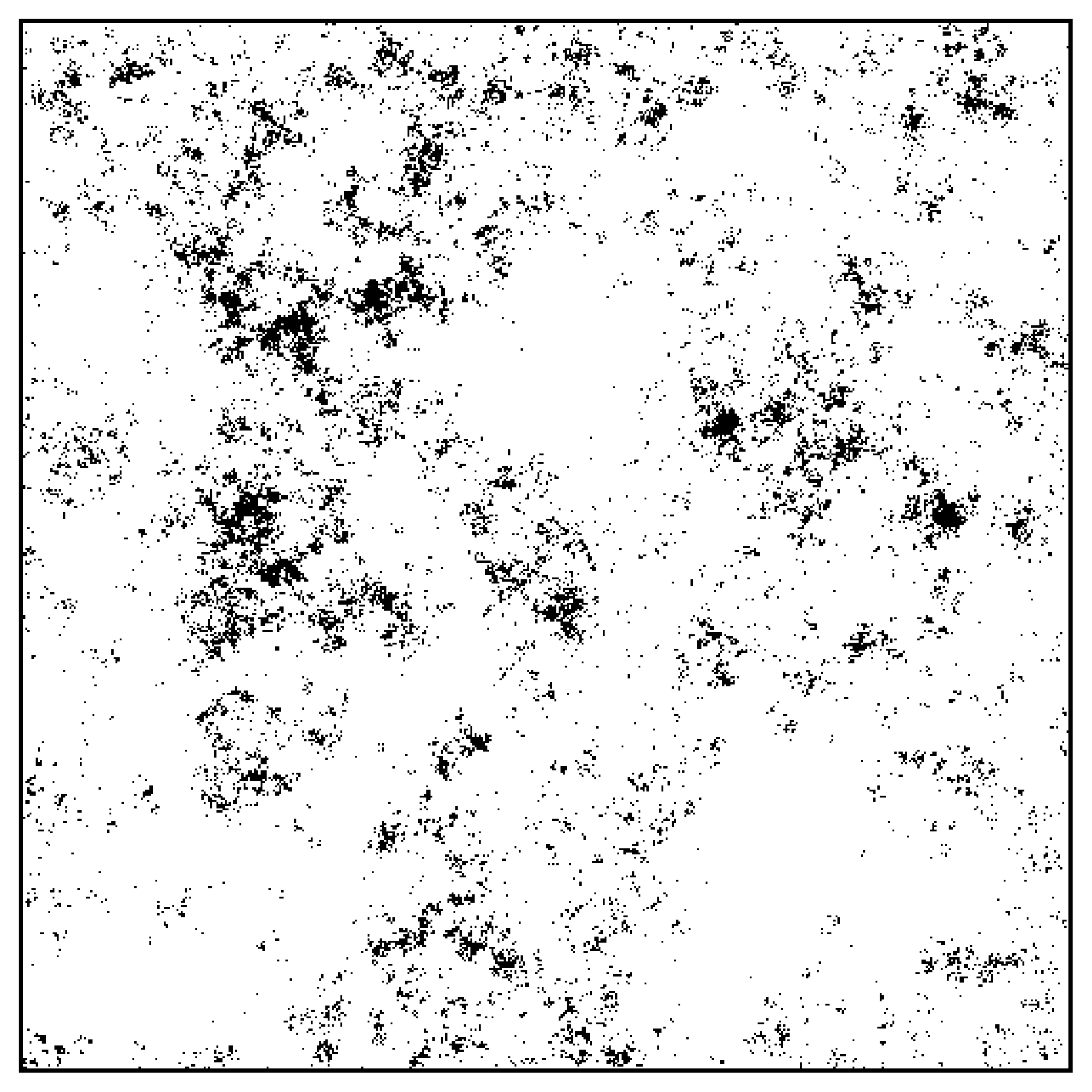}
\includegraphics[width=0.3\textwidth]{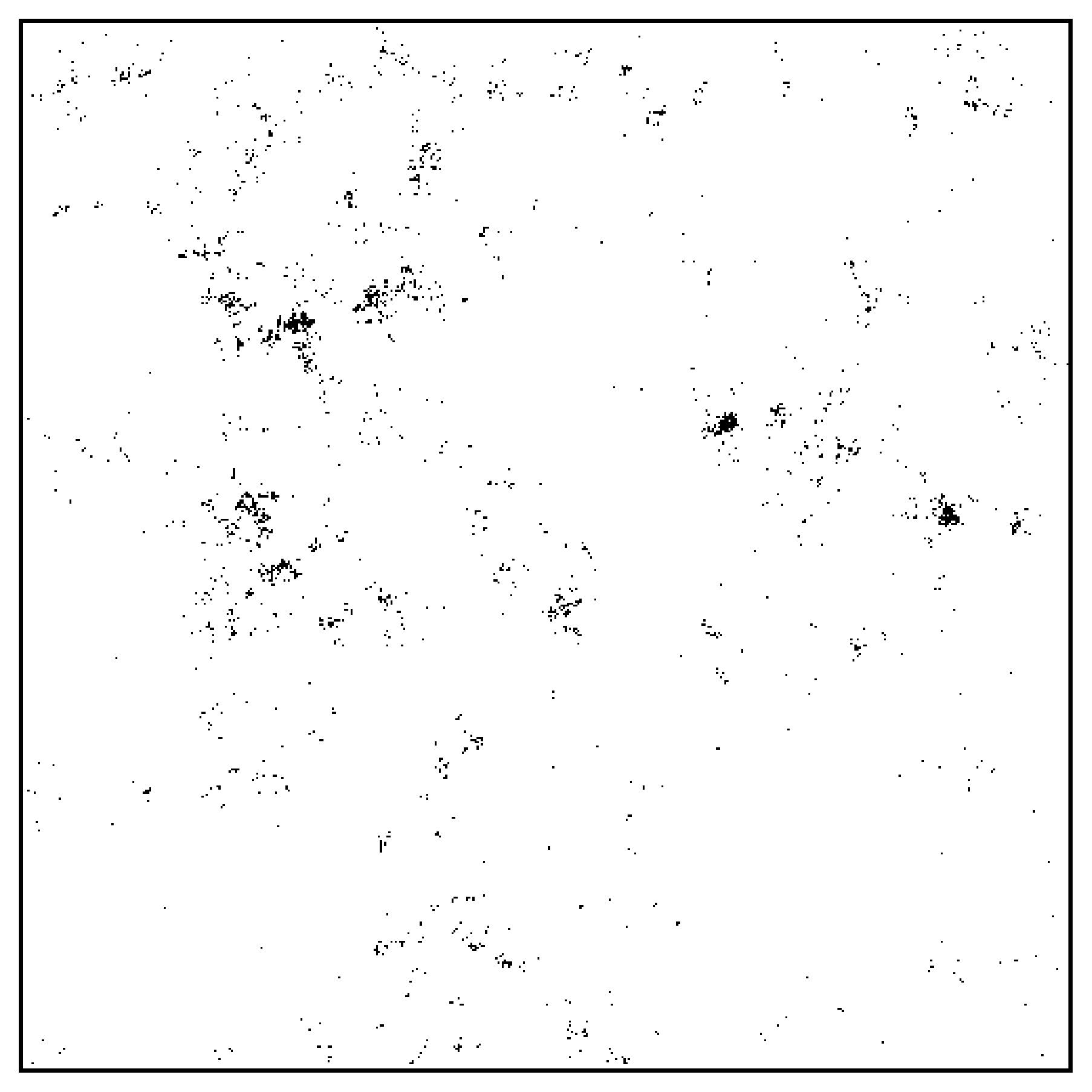}
}
\begin{quote}
\small 
\vglue-0.2cm
\caption{
\label{fig1}
Plots of the level sets \eqref{E:1.3} for a sample of the \DGFF{} on a square domain with $N:=500$ and $\lambda$ taking values (as labeled left to right) $0.1$, $0.3$ and $0.5$, respectively. The clustering (and fractal) nature of these sets is quite apparent.}
\normalsize
\end{quote}
\end{figure}

The objective of the present paper is to show that the intermediate level set~\eqref{E:1.3} admits a non-trivial \emph{scaling limit} which can be quite explicitly characterized. 
 A number of issues need to be addressed when setting the problem up mathematically. The first one is a proper formulation of the limit.   Indeed, after scaling the space by~$N$, the set \eqref{E:1.3} is increasingly dense everywhere in  the unit square $[0,1]^2$  and so taking its limit directly (e.g., in the topology of Hausdorff convergence) does not seem useful. We resolve this by encoding the level set into the point measure
\begin{equation}
\label{E:1.9}
\sum_{x\in  V_N}\delta_{x/N}\otimes\delta_{\,h_x-a_N},
\end{equation}
where $a_N$ is a scale sequence such that, in light of \eqref{E:1.2},
\begin{equation}
\label{E:1.10}
\frac{a_N}{\log N}\,\,\underset{N\to\infty}\longrightarrow\,\,2\sqrt g\,\lambda
\end{equation}
for some~$\lambda\in(0,1)$. The next issue is unbounded mass,  and that even under integration against compactly-supported continuous functions.  We resolve this by showing that \eqref{E:1.9} can be normalized by a \emph{deterministic} quantity so that a non-degenerate distributional limit becomes possible.  Another issue is dependence on the underlying domain; we resolve this by working in a class of lattice approximations~$D_N$ of a ``nice'' continuum set~$D\subset\C$. This will be useful technically and, moreover, will permit discussion of behavior under conformal maps.


\section{Main results}
\noindent 
Throughout the rest of the paper we will write~$h^V(x)$ to denote the \DGFF{} in~$V$ evaluated at~$x$. The presentation of our results opens up with the existence of the scaling limit.

\subsection{Existence of scaling limit}
We start by fixing the class of admissible domains. Let~$\mathfrak D$ be the class of all bounded open sets $D\subset\C$ with a finite number of connected components and with boundary~$\partial D$ that has only a finite number of connected components with each having a positive (Euclidean) diameter. Given~$D\in\mathfrak D$, we will restrict to sequences~$\{D_N\}$ of lattice domains such that
\begin{equation}
\label{E:2.1}
	 D_N \subseteq \bigl\{x \in \Z^2 \colon \textd_\infty(x/N,D^\cc)>\ffrac1N\bigr\}
\end{equation}
and, for each $\delta>0$ and all $N$ sufficiently large, also
\begin{equation}
\label{E:2.2}
D_N\supseteq\bigl\{x \in \Z^2 \colon \textd_\infty(x/N,D^\cc)>\delta\bigr\},
\end{equation}
where $\textd_\infty$ denotes the $\ell^\infty$-distance on~$\Z^2$.
Note that $x\in D_N$ implies $x/N\in D$. 

Next let us consider the \DGFF{} $\phi$ on~$\Z^2$ pinned to zero at the origin or, equivalently, the \DGFF{} on~$\Z^2\smallsetminus\{0\}$. This is a Gaussian process $\{\phi_x\colon x\in\Z^2\}$ with law to be denoted by~$\nu^0$ which is of mean zero and covariance given by
\begin{equation}
E_{\nu^0}\bigl(\phi_x\phi_y)=\fraka(x)+\fraka(y)-\fraka(x-y),
\end{equation}
where $\fraka\colon\Z^2\to[0,\infty)$ is the  potential kernel, i.e.,  the unique function that obeys $\fraka(0)=0$, is discrete harmonic on~$\Z^2\smallsetminus\{0\}$ and has the asymptotic form $\fraka(x)=g\log|x|+O(1)$ as~$|x|\to\infty$,
with~$|x|$ denoting the Euclidean norm of~$x$. 
Our main result is then:

\begin{theorem}
\label{thm-2.1}
For each $\lambda\in(0,1)$ and each $D\in\mathfrak D$, there is a random Borel measure~$Z^D_\lambda$ on~$\overline D$ with $E [Z^D_\lambda(\overline D)]\in(0,\infty)$ such that the following holds for each sequence $a_N$ satisfying \eqref{E:1.10} and each sequence~$D_N$ of scaled-up versions of~$D$ obeying \twoeqref{E:2.1}{E:2.2}: Set
\begin{equation}
\label{E:2.5}
K_N:=\frac{N^2}{\sqrt{\log N}}\,\,\texte^{-\frac{a_N^2}{2g\log N}}
\end{equation}
and, for each sample $h^{D_N}$ of the \DGFF{} in~$D_N$, define the point measure
\begin{equation}
\label{E:2.6}
\eta_N^D:=\frac1{K_N}\sum_{x\in D_N}\delta_{x/N}\otimes\delta_{\,h^{D_N}(x)-a_N}\otimes\delta_{\{h^{D_N}(x)-h^{D_N}(x+z)\colon z\in\Z^2\}}.
\end{equation}
Then, relative to the topology of vague convergence of measures on $\overline D\times\R\times\R^{\Z^2}$,
\begin{equation}
\label{E:2.7}
\eta^D_{N}\,\,\,\underset{N\to\infty}\Lawarrow\,\,\,Z_\lambda^D(\textd x)\,\otimes\,\texte^{-\alpha\lambda h}\textd h\,\otimes\,\nu_\lambda(\textd\phi),
\end{equation}
where $\alpha:=2/\sqrt g$ and~$\nu_\lambda$ is the probability measure on~$\R^{\Z^2}$ defined by
\begin{equation}
\label{E:2.8}
\nu_\lambda(\,\cdot\,)= \nu^0\Bigl(\,\phi+\frac2{\sqrt g}\lambda\,\fraka\in\cdot\Bigr),
\end{equation}
with~$\nu^0$ and $\fraka$ as specified above.
\end{theorem}

As it turns out, the convergence in \eqref{E:2.7} actually holds in   (a somewhat larger) space of Radon measures on  $\overline D\times(\R\cup\{+\infty\})\times\R^{\Z^2}$. As an immediate consequence we thus get:

\begin{corollary}
\label{cor2.2}
Let~$D\in\mathfrak D$. For any $D_N$ related to~$D$ as in \twoeqref{E:2.1}{E:2.2}, any sequence~$a_N$ such that~\eqref{E:1.10} holds with some~$\lambda\in(0,1)$, and~$K_N$ as in \eqref{E:2.5},
\begin{equation}
\label{E:2.9}
\frac1{K_N}\#\bigl\{x\in D_N\colon h^{D_N}(x)\ge a_N\bigr\}
\,\,\,\underset{N\to\infty}\Lawarrow\,\,\, (\alpha\lambda)^{-1}\, Z_\lambda^D(\overline D).
\end{equation}
Moreover, the random variables on the left are uniformly integrable and the convergence thus holds also under expectation.
\end{corollary}

 As is readily checked,  \eqref{E:1.10} yields $K_N = N^{2(1-\lambda^2)+o(1)}$ as~$N\to\infty$. 
Since $Z_\lambda^D(\overline D)\in(0,\infty)$ a.s.\ (see Theorem~\ref{thm-2.3}(7)), \eqref{E:2.9} determines the asymptotic size of the level set \eqref{E:1.3}. This strengthens considerably the aforementioned conclusion of Daviaud~\cite{Daviaud}.  The asymptotic positivity of the size of the level set normalized by its expectation appears already in the recent work by the first author, Ding and Goswami~\cite{BDG}. 

\subsection{Properties of $Z_\lambda$-measures}
In order to make our description of the limit law in \eqref{E:2.7} complete, we have to characterize the law of the random measure $Z^D_\lambda$. For this, we will first note a number of properties of these measures that will in fact be proved jointly with the above convergence theorem. 

We first need some additional notations. For each~$D\in\mathfrak D$ and each~$x\in D$, let~$\Pi^D(x,\cdot)$ denote the harmonic measure on~$\partial D$ relative to~$x$. As is well known (see, e.g., \cite[Lemma~2.3]{BL2}), if $\wt D\subseteq D$ are two admissible domains, then
\begin{equation}
\label{E:2.10}
C^{D,\wt D}(x,y) 
= g\int_{\partial D}\Pi^D(x,\textd z)\log|y-z|-g\int_{\partial \wt D}\Pi^{\wt D}(x,\textd z)\log|y-z|
\end{equation}
defines a symmetric, positive semi-definite function $C^{D,\wt D}\colon\wt D\times\wt D\to\R$ which is analytic in both variables. We may thus define $\{\Phi^{D,\wt D}(x)\colon x\in\wt D\}$ to be a Gaussian field with mean zero and covariance $C^{D,\wt D}$; this field has smooth sample paths a.s.
For~$\lambda\ge0$, we define $\psi^D_\lambda\colon D\to[0,\infty)$ by
\begin{equation}
\psi^D_\lambda(x):=\exp\Bigl\{\,2\lambda^2\int_{\partial D}\Pi^D(x,\textd z)\log|x-z|\Bigr\}.
\end{equation}
For~$D$ simply connected, $\psi_\lambda^D(x)$ is the $2\lambda^2$-th power of the conformal radius of~$D$ from~$x$. Setting $\psi^D_\lambda(x):=0$ for~$x\in\partial D$, the resulting function is continuous on~$\overline D$.  Writing $\leb(A)$ to denote the Lebesgue measure of a (measurable) set~$A\subset\R^2$, we then have:

\begin{theorem}
\label{thm-2.3}
Let $\lambda\in(0,1)$ and recall that~$\alpha:=2/\sqrt g$. Then the family of (laws of) random measures $\{Z_\lambda^D\colon D\in\mathfrak D\}$ obeys the following properties:
\begin{enumerate}
\item[(1)] For each~$D\in\mathfrak D$, the measure $Z^D_\lambda$ is supported on~$D$; i.e., $Z^D_\lambda(\R^2\smallsetminus D)=0$ a.s.
\item[(2)] If $A\subset D\in\mathfrak D$ is measurable with $\leb(A)=0$, then $Z_\lambda^D(A)=0$ a.s.
\item[(3)] There is $c=c(\lambda)\in(0,\infty)$ such that for each~$D\in\mathfrak D$ and each measurable $A\subset D$,
\begin{equation}
\label{E:2.12}
E Z^D_\lambda(A) = c\int_A\psi^D_\lambda(x)\textd x.
\end{equation}
\item[(4)] If $D,\wt D\in\mathfrak D$ obey $D\cap \wt D=\emptyset$, then 
\begin{equation}
\label{E:2.13}
Z_\lambda^{D\cup\wt D}(\textd x)\,\,\laweq\,\, Z_\lambda^{D}(\textd x)+Z_\lambda^{\wt D}(\textd x),
\end{equation}
with the measures $Z_\lambda^{D}$ and $Z_\lambda^{\wt D}$ on the right regarded as independent.
\item[(5)] (Gibbs-Markov) If $D,\wt D\in\mathfrak D$ obey $\wt D\subseteq D$ and $\leb(D\smallsetminus \wt D)=0$, then
\begin{equation}
\label{E:2.14}
Z_\lambda^D(\textd x)\,\laweq\,\texte^{\alpha\lambda\Phi^{D,\wt D}(x)}\,Z_\lambda^{\wt D}(\textd x),
\end{equation}
where $\{\Phi^{D,\wt D}(x)\colon x\in\wt D\}$ is as above and is regarded as independent of $Z_\lambda^{\wt D}$.
\item[(6)] The law of $Z^D_\lambda$ is translation invariant; $Z^{a+D}_\lambda(a+\textd x)\laweq Z^D_\lambda(\textd x)$ for each $a\in\R^2$.
\item[(7)] For each $A\subset D$ non-empty and open, $Z_\lambda^D(A)>0$ a.s.
\end{enumerate}
The properties (1-6),  for a given~$c>0$ in \eqref{E:2.12},  determine the laws of $\{Z_\lambda^D\colon D\in\mathfrak D\}$ uniquely.
\end{theorem}

Obviously, (2) is a special case of (3) although we prefer to state these separately. The constant~$c$ in \eqref{E:2.12} can be computed explicitly; just compare \eqref{E:2.9} with \eqref{E:3.6}. Perhaps the most important property of all is~(5). Here we note that the measure on the right of~\eqref{E:2.13} is well defined due to the fact that $Z_\lambda^{\wt D}(D\smallsetminus \wt D)=0$ a.s.\ thanks to property~(1), and this carries no loss on the left-hand side because also $Z_\lambda^{ D}(D\smallsetminus \wt D)=0$ a.s. thanks to property~(2). We will refer to~(5) --- sometimes also in conjunction with~(4) --- as the Gibbs-Markov property. This is because properties (4-5) arise directly from the Gibbs-Markov decomposition of the \DGFF{}; cf~\eqref{E:A.6}. 

By property~(6), the law of~$Z^D_\lambda$ transforms canonically under the spatial shifts. The behavior of~$Z^D_\lambda$ under scaling of~$D$ is more subtle as it is intimately tied to the existence of the limit \eqref{E:2.7} and its independence of the sequence of discrete domains~$D_N$ \emph{and} of how the centering sequence~$a_N$ achieves the overall asymptotic~\eqref{E:1.10}. Once a suitable scaling relation is established, the Gibbs-Markov property yields also rotation invariance and, in fact, leads to:

\begin{theorem}
\label{thm-2.4}
Let~$\lambda\in(0,1)$. Under any conformal bijection $f\colon D\to f(D)$ between the admissible domains $D,f(D)\in\mathfrak D$, the laws of the above measures transform as
\begin{equation}
\label{E:2.16}
Z_\lambda^{f(D)}\circ f(\textd x)\,\,\laweq\,\, |f'(x)|^{2+2\lambda^2}\, Z_\lambda^D(\textd x).
\end{equation}
\end{theorem}

\subsection{Connection to Liouville Quantum Gravity}
Although the above properties already determine the law of $\{Z_\lambda^D\colon D\in\mathfrak D\}$ uniquely, we are able to make even a more explicit connection with the so called Liouville Quantum Gravity measures that have been introduced and studied by Duplantier and Sheffield~\cite{Duplantier-Sheffield}. 

Again we start with some definitions. Let $\cmss H_0^1(D)$ denote the closure of the set of smooth, functions with compact support in~$D$ with respect to the norm induced by the Dirichlet inner product $\langle f,g\rangle_\nabla:= \frac14\int_D\nabla f(x)\cdot\nabla g(x)\textd x$. Given a sequence $\{X_n\colon n\ge1\}$ of i.i.d.\ standard normal random variables and an orthonormal basis $\{f_n\colon n\ge1\}$ in $\cmss H_0^1(D)$, define
\begin{equation}
\label{E:2.17}
\varphi_n(x):=\sum_{k=1}^n X_k f_k(x).
\end{equation}
For each~$\beta\in[0,\infty)$, define the random measure
\begin{equation}
\label{E:2.17a}
\mu_n^{D,\beta}(\textd x):=1_D(x)\texte^{\beta\varphi_n(x)-\frac{\beta^2}2E[\varphi_n(x)^2]}\,\textd x.
\end{equation}
As goes back to Kahane~\cite{Kahane}, there exists a random, a.s.\ finite (albeit possibly trivial) Borel measure $\mu_\infty^{D,\beta}$ --- called the \emph{Gaussian multiplicative chaos} associated with the continuum Gaussian Free Field --- which is concentrated on~$D$ and such that, for each measurable set~$A$,
\begin{equation}
\label{E:2.18}
\mu_n^{D,\beta}(A)\,\underset{n\to\infty}\longrightarrow\,\mu_\infty^{D,\beta}(A)\quad\text{a.s.}
\end{equation}
It is also known (cf a remark after Rhodes and Vargas~\cite[Theorem~5.5]{RV-review}) that for each $\beta\in(0,\beta_\cc)$, where (in our normalization) $\beta_\cc:=\alpha=2/\sqrt g$, we have $\mu_\infty^{D,\beta}(D)>0$ a.s. Moreover, as was shown in \cite[Theorem~5.5]{RV-review}, the law of the limit measure does not depend on the choice of the above orthonormal basis. (In fact, thanks to Shamov~\cite[Corollary~5]{Shamov}, the law of~$\mu_\infty^{D,\beta}$ is determined solely by its expectation and the way the measure transforms under the Cameron-Martin shifts of the underlying CGFF.) With this stated, we now claim:

\begin{theorem}
\label{thm-LQG}
Let $\lambda\in(0,1)$, $\alpha:=2/\sqrt g$ and consider the family of measures $\{Z_\lambda^D\colon D\in\mathfrak D\}$ as above. Then, for $c\in(0,\infty)$ as in \eqref{E:2.12} and for each~$D\in\mathfrak D$,
\begin{equation}
\label{E:2.20}
Z_\lambda^D(\textd x) \,\laweq\,c\psi^D_\lambda(x)\,\mu_\infty^{D,\,\lambda\alpha}(\textd x).
\end{equation}
In particular, $Z_\lambda^D$ has the law of the Liouville Quantum Gravity measure in~$D$ corresponding to (subcritical) parameter $\beta:=\lambda\alpha$. 
\end{theorem}

\begin{figure}[t]
\vglue0.2cm
\centerline{\includegraphics[width=0.7\textwidth]{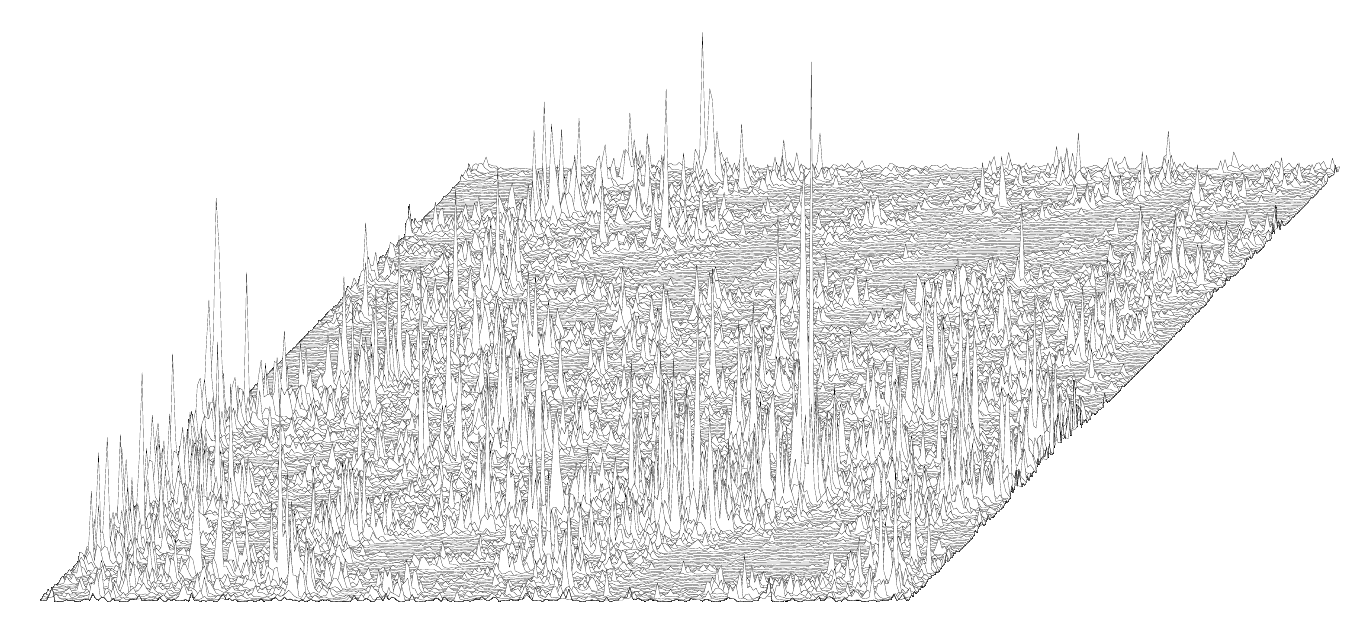}
}
\begin{quote}
\small 
\caption{
\label{fig3}
A sample of the LQG measure $\psi_\lambda^D(x)\mu^{D,\beta}_\infty(\textd x)$ for~$D$ a unit square and parameters $\beta=\lambda\alpha$ and~$\lambda=0.3$. The high points indicate places of high local intensity.}
\normalsize
\end{quote}
\end{figure}

\subsection{Remarks and open problems}
We proceed by a series of remarks and questions left to be studied.

\begin{figure}[t]
\centerline{\includegraphics[width=0.58\textwidth]{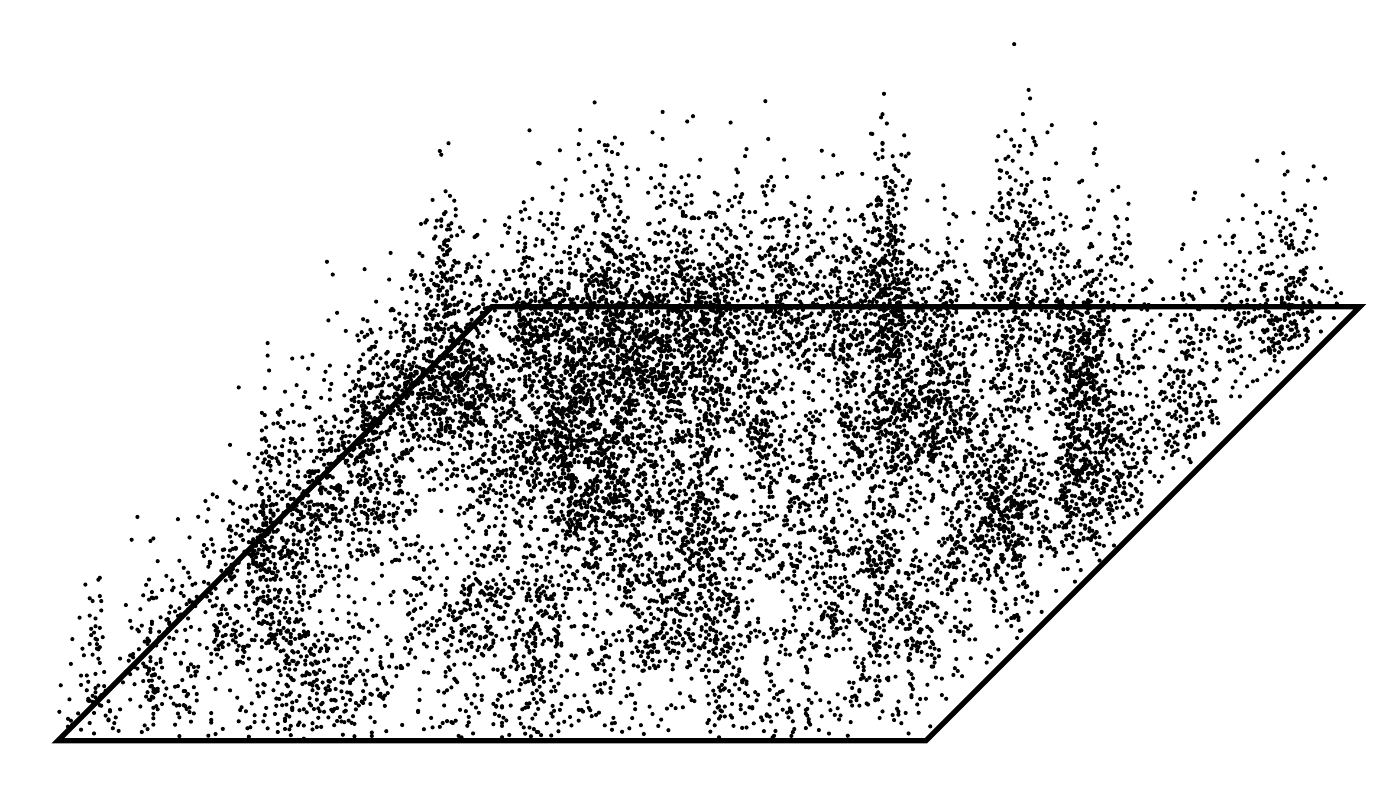}
\hglue-2.2cm
\includegraphics[width=0.58\textwidth]{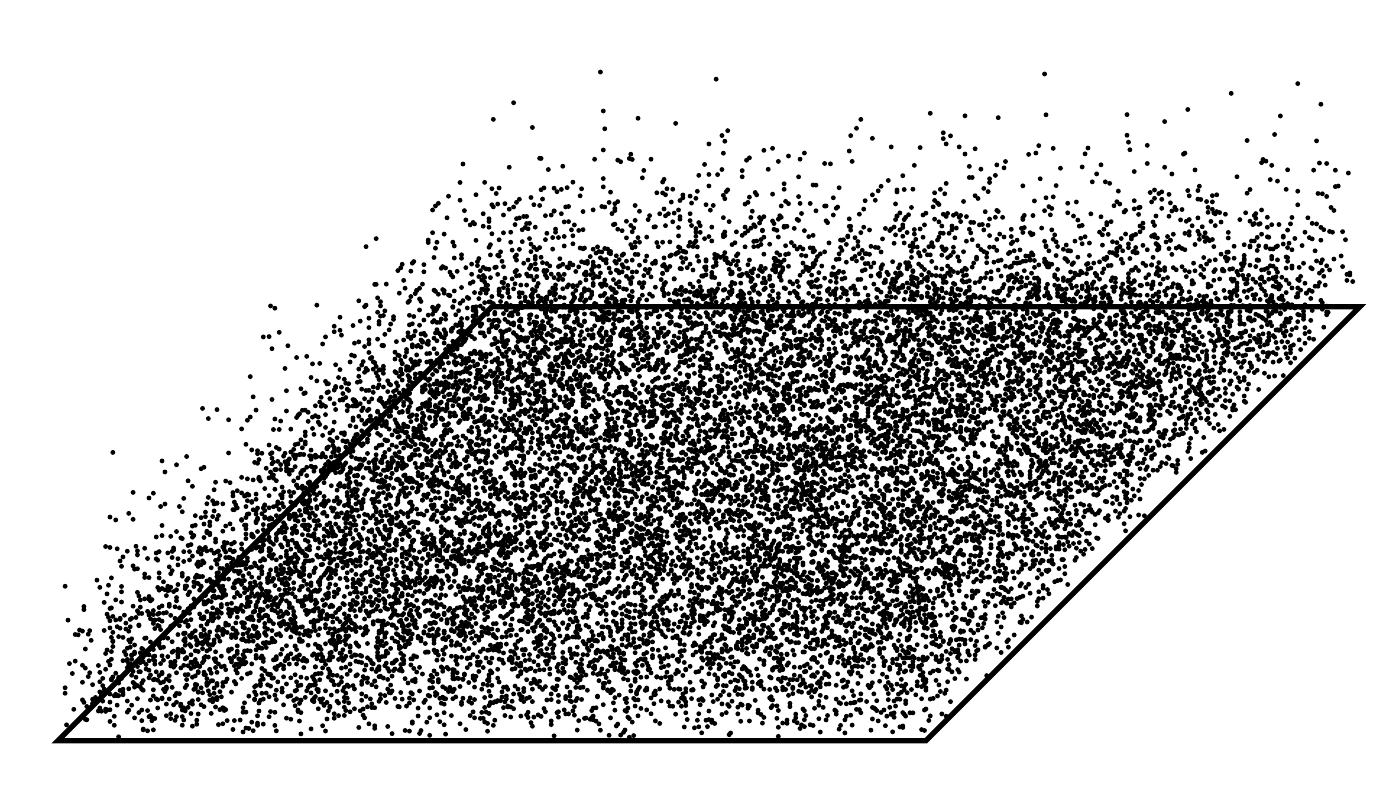}
}
\begin{quote}
\small 
\vglue-0.2cm
\caption{
\label{fig2}
Left: A sample of the measure in \eqref{E:1.9} on a square of side-length~$N:=300$ with~$\lambda:=0.2$. Right: A corresponding sample for i.i.d.\ normals with mean zero and variance~$g\log N$. Only the points with positive vertical coordinate are depicted.}
\normalsize
\end{quote}
\end{figure}

\medskip\noindent
(1) \textit{General Gaussian processes: }
We believe that the \emph{form} of the limit measure in~\eqref{E:2.7} is actually quite universal. For instance, for i.i.d.\ Gaussians indexed by the vertices in~$D_N$ (see Fig.~\ref{fig2}) with variance $g\log N$ with same~$K_N$ we get the same limit statement with~$Z^D_\lambda$ replaced by (a multiple of) the Lebesgue measure on~$D$ and~$\nu_\lambda$ by the point mass concentrated on $\phi$ defined by $\phi_0:=0$ and $\phi_x:=-\infty$ for~$x\ne0$. That~$Z^D_\lambda$ is itself random in the case of the \DGFF{} is a reflection of long-range correlations.

\medskip\noindent
(2) \textit{Simultaneous limit for all~$\lambda$:\ }
Our proofs are technically based on the computation of the first two moments of the measure~$\eta^D_N$ integrated against compactly-supported, continuous functions. (This is literally true when $\lambda<1/\sqrt2$ with a  truncation needed for complementary~$\lambda$.) One could use similar techniques to study the level sets for several values of~$\lambda$ simultaneously but the number of required moments seems to increase with the number of the levels to be controlled. In particular, we presently do not see a way to solve:

\begin{problem}
Find a way to extract a joint distributional limit of the level sets \eqref{E:1.3}, or their associated point measures \eqref{E:2.6}, simultaneously for all $\lambda\in(0,1)$. 
\end{problem} 

\noindent
Our belief that a joint limit should exist is supported by the fact that \twoeqref{E:2.17}{E:2.18}, and a suitable continuity argument, define the LQG measure for all~$\beta\in(0,\beta_\cc)$ at the same time. This is because the LQG measure is a measurable function of the underlying CGFF.

\medskip\noindent
\textit{(3) Relation to extremal process:\ }
Our point process-based approach is strongly motivated by recent advances in the understanding of the extremal values of the \DGFF{}; i.e., roughly speaking, the set \eqref{E:1.3} with~$\lambda:=1$. Here, first, building on the work of Bolthausen, Deuschel and Zeitouni~\cite{BDZ} and Bramson and Zei\-touni~\cite{BZ}, Bramson, Ding and Zeitouni~\cite{BDingZ}  showed that for 
\begin{equation}
\label{E:1.4}
m_N:=2\sqrt g\,\log N-\frac34\sqrt g\,\log\log N,
\end{equation}
the centered maximum, $\max_{x\in V_N}h_x-m_N$, converges to a non-degenerate distributional limit as $N\to\infty$. Then, in~\cite{BL1,BL2,BL3}, the present authors described the limit of the full extremal process for the \DGFF{} expressed in terms of the ``structured'' point process,
\begin{equation}
\label{E:1.5}
\eta^D_{N,r}:=\sum_{x\in D_N}
\1_{\{h_x=\max_{z\in \Lambda_r(x)}h_z\}}
\delta_{\,x/N}\otimes\delta_{\,h_x-m_N}\otimes\delta_{\,\{h_{x}-h_{x+z}\colon z\in\Z^2\}},
\end{equation}
where~$\Lambda_r(x):=\{z\in\Z^2\colon |z-x|\le r\}$ and the indicator thus effectively restricts the sum to the points where the field has an $r$-local maximum. 
The main result of~\cite{BL1,BL2,BL3} is that, for any $r_N$ with $r_N\to\infty$ and $N/r_N\to\infty$, relative to the topology of vague convergence,
\begin{equation}
\label{E:1.6}
\eta^D_{N,r_N}\,\,\,\underset{N\to\infty}\Lawarrow\,\,\,\text{\rm PPP}\bigl(\,Z^D(\textd x)\,\otimes\,\texte^{-\alpha h}\textd h\,\otimes\,\nu_1(\textd\phi)\bigr),
\end{equation}
where $\text{PPP}$ stands for Poisson point process, $Z^D$ is a random a.s.-finite Borel measure on~$D$,
and~$\nu$ is a (deterministic) probability measure on~$[0,\infty)^{\Z^2}$; cf~\eqref{E:2.23}.
In~\cite{BL2}, the~$Z^D$ measure was shown to coincide with a version of the critical Liouville Quantum Gravity;  unfortunately, the identification with the versions constructed in Duplantier, Rhodes, Sheffield and Vargas~\cite{DRSV2,DRSV1} has not yet been fully completed. Notwithstanding, in~\cite{BL2} we show that
\begin{equation}
\label{E:2.16b}
Z_1^{f(D)}\circ f(\textd x)\,\,\laweq\,\, |f'(x)|^{4}\, Z_1^D(\textd x).
\end{equation}
This corresponds, at least formally, to $\lambda\uparrow1$ limit of \eqref{E:2.16}. The $\lambda\downarrow0$ limit reduces~$Z^D_\lambda$ to (a multiple of) the Lebesgue measure on~$D$.

\medskip\noindent
(4) \textit{Conformal invariance and uniqueness of LQG measure:\ }
The previous remark brings us to the formulation of the conformal transformation rule in Theorem~\ref{thm-2.4}. In Duplantier and Sheffield~\cite[Proposition~1.2]{Duplantier-Sheffield}, a version of this rule is stated for the LQG measure in the following form: {Supposing} that the CGFF~$h^D$ in a domain~$D$ transforms under any conformal map $f\colon D\to f(D)$ (in our parametrization) as follows
\begin{equation}
\label{E:2.21}
h^{f(D)}\laweq h^D\circ f+\Bigl(\frac2\beta+\frac{2\beta}{\alpha^2}\Bigr)\log|f'|
\end{equation}
(note that this requires working with CGFF of non-zero mean and/or Cameron-Martin shifts),
the~LQG measure in~$D$ for parameter~$\beta$ transforms into the corresponding LQG measure in~$f(D)$. Unfortunately, this harbors a technical caveat: The measure must be realized as a \emph{unique} function of the CGFF function of the CGFF or, at least, one that is independent of the approximation scheme used to define it. This is in fact a subtle issue that has been fully settled only quite recently (cf the aforementioned references to \cite[Theorem~5.5]{RV-review} or~\cite[Corollary~5]{Shamov}). 

Our approach to Theorem~\ref{thm-2.4} has the advantage that it works solely with the family of random measures $\{Z^D_\lambda\colon D\in\mathfrak D\}$ and, in particular, avoids dealing with the uniqueness of LQG and/or its dependence on the underlying CGFF. In fact, Theorem~\ref{thm-2.4} could concisely be stated as: 

\begin{namedtheorem}[Theorem~\ref{thm-2.4}']
Every family of measures $\{Z_\lambda^D\colon D\in\mathfrak D\}$ satisfying properties (1-6) in Theorem~\ref{thm-2.3} obeys also \eqref{E:2.16}, for each conformal bijection $f\colon D\to f(D)$ with $D,f(D)\in\mathfrak D$.
\end{namedtheorem} 

\noindent
Our proof of \eqref{E:2.16} reduces, after some minor amount of preparation, to the same argument as the proof of \eqref{E:2.16b} in Theorem~7.2 in~\cite{BL2}.  The only time when we need to invoke uniqueness of the LQG measure is, quite naturally, when we identify $Z^D_\lambda$ with the (unique) LQG measure in the proof of Theorem~\ref{thm-LQG}.

\medskip\noindent
\textit{(5) Fluctuations around random limit:\ }
Our next remark concerns going beyond the limit statement~\eqref{E:2.7}. The point is that the limit measure, albeit random, captures only the leading-order growth of the level set. We thus pose:

\begin{problem}
\label{pb2.7}
Characterize the limit law of the (suitably scaled) fluctuations in the limit \eqref{E:2.7}.
\end{problem}

\noindent
To make the formulation easier, one may choose to work in the setting when all the \DGFF{}s are defined on the same probability space as the limit LQG measure. Perhaps the easiest underlying graph for this is the triangular lattice.

\medskip\noindent
\textit{(6) Crossover to critical regime:\ }
Of quite some interest is how the behavior for $\lambda<1$ blends with that at $\lambda=1$. Our proofs only apply for~$a_N$ such that $a_N/\log N$ is, in the limit, strictly less than~$2\sqrt g$. This is for a good reason: When $a_N/\log N\to2\sqrt g$, the growth rate of the requisite normalizing sequence should be slower than \eqref{E:2.5}. This stems from a subtle entropic-repulsion effect that lies at the heart of the paper~\cite{BL3} and can be seen by noting that plugging $a_N:=m_N$ for $m_N$ as in \eqref{E:1.4} results in~$K_N\sim\log N$. We thus pose:

\begin{problem}
Suppose that $a_N/\log N\to 2\sqrt g$ yet $m_N-a_N\to\infty$. Prove that, for a suitably re-defined~$K_N$, we still have~\eqref{E:2.6} with~$Z^D_\lambda$ replaced by the critical LQG measure~$Z^D$ from \eqref{E:1.6}.
\end{problem}

\noindent
An affirmative resolution of this problem may in fact require that $m_N-a_N$ tends to infinity at some minimal rate. A reasonable guess is that $m_N-a_N$ of order~$\sqrt{\log N}$ should already be enough. This scale appears naturally as it marks the level where the discrete approximation to the critical LQG measure is typically supported.

\medskip\noindent
\textit{(7) Beyond 2D DGFF:\ }
A natural question is of course whether the above results are in any sense \emph{universal} for other models that are, at least at large spatial scales, well captured by CGFF. These include general logarithmically-correlated Gaussian fields, gradient models and local time of a two-dimensional simple random walk. Some progress on these has already been made (e.g., Belius and Wu~\cite{Belius-Wu}, Abe~\cite{Abe}).

\subsection{Proof strategy}
The overall strategy of our proofs is rather simple. Through moment calculations for the size of the level set, we establish tightness and asymptotic non-triviality of the measures $\{\eta^D_N\colon N\ge1\}$ relative to the vague topology. This permits extraction of subsequential weak limits. We then proceed to derive various relations that such limits have to satisfy which ultimately characterize them uniquely. This proves existence of the limit as well as its desired properties.

The specific ``characterization'' steps are as follows. First we focus only on the measures \eqref{E:2.6} restricted to the first two coordinates. One more (subtle) second moment calculation shows that every subsequential limit $\eta^D$ of such two-coordinate measures admits the decomposition 
\begin{equation}
\eta^D(\textd x\,\textd h) = Z^D_\lambda(\textd x)\otimes\texte^{-\alpha\lambda h}\textd h,
\end{equation}
with $Z^D_\lambda$ a non-degenerate, a.s.\ finite measure whose law possibly depends on the subsequence, the sequence of approximation domains~$D_N$ as well as the way~$a_N$ approaches the limit~\eqref{E:1.10}. 

Next we demonstrate that the measures $Z^D_\lambda$, with~$D$ restricted to a suitable countable collection of domains (this is the best one can hope to have when extracting limits by subsequences), obey properties (1-7) in Theorem~\ref{thm-2.3}. Property~(5) is then particularly important, as it yields a representation of~$Z^D_\lambda$, for~$D$ a dyadic square, in terms of a multiplicative chaos measure. From here we get uniqueness of the law of~$Z^D_\lambda$ on dyadic squares; one more use of property~(5) then extends this to all $D\in\mathfrak D$. The existence of the limit of~$\eta^D_N$, and its independence of the approximation domains~$D_N$ or the specific way~$a_N$ achieves the limit~$a_N/\log N\to 2\sqrt g\,\lambda$, follow.  This pretty much completes the proof of Theorems~\ref{thm-2.1},~\ref{thm-2.3} and Theorem~\ref{thm-LQG} for the two-coordinate measures.

The sheer existence of the limit (and translation invariance of the \DGFF{}) now implies the transformation rule for shifts and scaling of the underlying domain:
\begin{equation}
\label{E:2.22}
Z^{a+r D}_\lambda(a+r\textd x) \laweq r^{2+2\lambda^2}Z^D_\lambda(\textd x),\qquad a\in\R^2,\,r>0.
\end{equation}
The representation using multiplicative chaos adds rotation invariance to this as well. One more use of property~(5) for a decomposition of a given~$D$ using a myriad of tiny dyadic squares then permits us to apply these symmetries ``infinitesimally'' thus proving, with the help of conformal invariance of the ``binding'' fields~$\Phi^{D,\wt D}$, Theorem~\ref{thm-2.4}. 

As a final step, we extend control to the full three-coordinate process \eqref{E:2.6}. This boils down to yet another moment calculation, which yields factorization of the limit into the product measure on the right of \eqref{E:2.7}. It is easy to see why the limit law of $h^{D_N}(x)-h^{D_N}(x+\cdot)$ should be described by \eqref{E:2.8}: Conditioning on~$h^{D_N}(x)$ to be roughly $2\sqrt g\,\lambda\log N$ changes the mean of the field at~$y$ by, roughly, $\frac 2{\sqrt g}\lambda\fraka(y-x)$ while the variance tends to that of the \DGFF{} on~$\Z^2\smallsetminus\{0\}$. We note that a similar reasoning applies also to the extremal process limit \eqref{E:1.6} except there the ``cluster law''~$\nu_1$ requires an extra conditioning (to ensure a local maximum at~$x$) and taking a limit,
\begin{equation}
\label{E:2.23}
\nu_1(\cdot) = \lim_{r\to\infty}\,\nu^0\biggl(\phi+\frac2{\sqrt g}\fraka\in\cdot\,\bigg|\, \phi(x)+\frac2{\sqrt g}\fraka(x)\ge0\colon|x|\le r\biggr).
\end{equation}
This limit is singular, which is a source of much headache in the proofs of \cite{BL3}.

\begin{remark}
The above strategy --- extract a subsequential limit and then prove its uniqueness --- also lies at the core of our earlier work~\cite{BL1,BL2,BL3} on the extrema of the \DGFF{}. However, the two approaches are technically quite different, both in the proof of the factorization (which, in~\cite{BL1,BL2,BL3}, relies on a connection with particle systems) and in the proof of uniqueness (which, for the extremal values, relies on the existence of the limit of the centered absolute maximum).
\end{remark}

Since we are dealing with scaling limits of the \DGFF{}, it is no surprise that we will need to invoke bounds on, as well as limits of, the Green function in various lattice domains of interest. The limit statements in particular require weak convergence of the harmonic measure on~$\frac1N\partial D_N$ to that on~$\partial D$. This is where the containment~$D\in\mathfrak D$ and relations \twoeqref{E:2.1}{E:2.2} are required. To make referencing easier, we collect the needed statements in the Appendix.

\section{Proofs in the second moment regime}
\noindent
We are now ready to commence the exposition of our proofs. As noted above, the starting point are calculations of the first two moments of the size of the level set. These  are straightforward for $\lambda\in(0,1/{\sqrt2})$ but harder in the complementary regime of~$\lambda$, where additional (albeit standard) truncations are required to keep the second moment comparable to the square of the first. For ease of exposition as well as pedagogical appeal, we will first deal with the former regime leaving the latter to a subsequent section. 

Throughout this section, we thus assume that~$a_N$ is a sequence such that~\eqref{E:1.10} holds for some~$\lambda\in(0,1/\sqrt2)$. We suppose that, for each~$D\in\mathfrak D$, a sequence $\{D_N\}$ of approximating lattice domains is given satisfying \twoeqref{E:2.1}{E:2.2}. Unless stated otherwise, all estimates will depend on the choice of~$D$ and the sequences $a_N$ and~$D_N$.

\subsection{Level-set size moments}
For each~$b\in\R$, define
\begin{equation}
\Gamma_N^D(b):=\bigl\{x\in D_N\colon h^{D_N}(x)\ge a_N+b\bigr\}.
\end{equation}
We begin by a bound on the overall size of~$\Gamma_N^D(b)$:

\begin{lemma}
\label{lemma-3.1}
For each~$D\in\mathfrak D$ there is~$c = c(D)\in(0,\infty)$ such that for all $b\in\R$ with $|b|\le\log N$, all non-negative sequences~$a_N$ satisfying $a_N\le 4\sqrt g\log N$, all sequences $D_N$ satisfying \twoeqref{E:2.1}{E:2.2}, all~$A\subset D_N$ and all~$N\ge1$, we have
\begin{equation}
E\bigl|\Gamma_N^D(b)\cap A\bigr|\le cK_N\,\frac{|A|}{N^2}\,\texte^{-\frac{a_N}{g\log N} b}\,.
\end{equation}
\end{lemma}

\begin{proofsect}{Proof}
The claim will follow by summing over~$x\in D_N$ once we prove that, for some constant~$c$ depending only on the diameter of~$D$, we have
\begin{equation}
P\bigl(h^{D_N}(x)\ge a_N+b\bigr)\le c\frac1{\sqrt{\log N}}\texte^{-\frac{a_N^2}{2g\log N}}\,\texte^{-\frac{a_N}{g\log N} b}
\end{equation}
uniformly in~$x\in D_N$ and in $b\in[-\log N,\log N]$. To this end we first invoke the Gibbs-Markov property of the \DGFF{} (cf~\eqref{E:A.6}) to note that, if $U\subset V$ are finite lattice domains, then by writing $h^V(x)\laweq h^U(x)+\varphi^{V,U}(x)$ and requiring $\varphi^{V,U}(x)\ge0$,
\begin{equation}
\label{E:3.4}
P\bigl(h^U(x)\ge a\bigr)\le 2 P\bigl(h^V(x)\ge a\bigr)\,.
\end{equation}
By enlarging~$D_N$ to, say, a square domain $\wt D_N$ of comparable diameter, we may thus assume that all~$x\in D_N$ lie deep inside~$\wt D_N$. The variance of $h^{\wt D_N}(x)$ is then within a constant of~$g\log N$ uniformly in~$x\in D_N$ and so we get, for some~$c>0$ independent of~$N$,
\begin{equation}
P\bigl(h^{\wt D_N}(x)\ge a_N+b\bigr)
\le\frac1{\sqrt{2\pi}}\frac1{\sqrt{g\log N-c}}\int_b^\infty\texte^{-\frac{(a_N+s)^2}{2g\log N+c}}\textd s.
\end{equation}
The $s\ge 2\log N$ portion of the integral is bounded directly; in the rest we use that~$c$ can be dropped at the cost of a suitable multiplicative term popping in the front. Bounding $(a_N+ s)^2\ge a_N^2+2a_N s$, the integral over~$s\in[b,2\log N]$ is then performed explicitly.
\end{proofsect}

With the overall scale under control, we can now calculate the leading-order asymptotic of the above expectation for nice-enough sets~$A$.

\begin{lemma}
\label{lemma-3.2}
There is a constant~$c_0\in\R$ such that for each~$b\in\R$ and each open set $A\subseteq D$,
\begin{equation}
\label{E:3.6}
E\bigl|\{x\in\Gamma_N^D(b)\colon x/N\in A\}\bigr| = 
\frac{\texte^{c_0\lambda^2+\lambda\alpha}}{\sqrt{2\pi g}}\texte^{-\alpha\lambda b}\Bigl[\,\int_A\psi^D_\lambda(x)\textd x+o(1)\Bigr]\,K_N.
\end{equation}
where $o(1)\to0$ as~$N\to\infty$ uniformly on compact sets of~$b$.
\end{lemma}

\begin{proofsect}{Proof}
Thanks to the uniform control from Lemma~\ref{lemma-3.1}, we may assume that the closure of~$A$ lies in~$D$.
We will need the asymptotic
\begin{equation}
\label{G-asymp}
G^{D_N}\bigl(\lfloor xN\rfloor,\lfloor xN\rfloor\bigr)
= g\log N + g\int_{\partial D}\Pi^D(x,\textd z)\log|x-z|+c_0+o(1)\,,
\end{equation}
with $c_0$ a constant and $o(1)\to0$ as $N\to\infty$ uniformly on compact sets in~$D$, and thus in~$x\in A$. (This is where the conditions on~$D_N$ are relevant, see \twoeqref{E:A.4}{G-asymp2}.) Now we repeat the calculation from the proof of Lemma~\ref{lemma-3.1} while keeping careful track of all non-vanishing terms. The boundedness and continuity of~$\psi^D_\lambda$ finally permit us to replace a Riemann sum by the integral in \eqref{E:3.6}.
\end{proofsect}

Our next lemma concerns the second moment estimate for the size of~$\Gamma_N^D(b)$. It is here where we need to limit the range of possible~$\lambda$:

\begin{lemma}
\label{lemma-3.3}
Suppose $0<\lambda<1/\sqrt2$. For each~$b_0\in\R$ and each~$D\in\mathfrak D$ there is $c_1=c_1(\lambda,b_0,D)\in(0,\infty)$ such that for each $b\in[-b_0,b_0]$ and each~$N\ge1$,
\begin{equation}
\label{e:2.3}
E \bigl(|\Gamma_N^D (b)|^2 \bigr)
\le c_1 K_N^2
\end{equation}
 Moreover, there is an absolute constant~$c_2\in(0,\infty)$ such that for all~$D\in\mathfrak D$,
\begin{equation}
\label{e:2.3b}
\limsup_{N\to\infty}\frac1{K_N^2}\,E \bigl(|\Gamma_N^D (0)|^2\bigr)\le
c_2\int_{D\times D}\biggl(\frac{[\diam D]^2}{|x-y|}\biggr)^{4\lambda^2}\,\textd x\textd y
\end{equation}
where $\diam D$ is the diameter of~$D$ in the Euclidean norm.
\end{lemma}

\begin{proofsect}{Proof of \eqref{e:2.3}}
Thinking, without much loss of generality, of~$b$ as absorbed into~$a_N$, we can assume $b:=0$ in the following.
Writing
\begin{equation}
\label{E:3.9a}
E \bigl(|\Gamma_N^D (0)|^2 \bigr) = \sum_{x,y\in D_N}P\bigl(h^{ D_N}(x)\ge a_N,\,h^{ D_N}(y)\ge a_N\bigr).
\end{equation}
we will need a good estimate on the probability on the right-hand side.  First we again take~$\wt D_N$ to be a neighborhood of~$D_N$ of diameter twice the diameter of~$D_N$ and note that, by the argument leading to \eqref{E:3.4} and the FKG inequality for $\varphi^{V,U}$ (implied by $\text{\rm Cov}(\varphi^{V,U}(x),\varphi^{V,U}(y))\ge0$)
\begin{equation}
\label{E:3.9ua}
P\bigl(h^{ D_N}(x)\ge a_N,\,h^{ D_N}(y)\ge a_N\bigr)
\le 4P\bigl(h^{\wt D_N}(x)\ge a_N,\,h^{\wt D_N}(y)\ge a_N\bigr)\,.
\end{equation}
Next we invoke the Gibbs-Markov decomposition (see~\eqref{E:A.6})
\begin{equation}
\label{E:3.10}
h^{\wt D_N}(y) = \frakg_x(y)h^{\wt D_N}(x)+\hat h^{\wt D_N\smallsetminus\{x\}}(y),
\end{equation} 
where $h^{\wt D_N}(x)$ and $\hat h^{\wt D_N\smallsetminus\{x\}}$ on the right-hand side are independent with $\hat h^{\wt D_N\smallsetminus\{x\}}$ having the law of the \DGFF{} in~$\wt D_N\smallsetminus\{x\}$ and where $\frakg_x$ is a function that is harmonic in~$\wt D_N\smallsetminus\{x\}$, vanishing outside~$\wt D_N$ and normalized such that $\frakg_x(x)=1$. Using this decomposition, the above probability is recast as
\begin{multline}
\label{E:3.11}
\quad
P\bigl(h^{\wt D_N}(x)\ge a_N,\,h^{\wt D_N}(y)\ge a_N\bigr)
\\
=\int_0^\infty P\Bigl(\,\hat h^{\wt D_N\smallsetminus\{x\}}(y)\ge a_N(1-\frakg_x(y))-s\frakg_x(y)\Bigr)P\bigl(h^{\wt D_N}(x)-a_N\in\textd s\bigr).
\end{multline}
 We will pick~$\delta>0$ and bound the right-hand side by $P(h^{\wt D_N}\ge a_N)$ when~$|x-y|\le\delta\sqrt{K_N}$ so let us suppose that $|x-y|>\delta\sqrt{K_N}$ from now on.

Observe that since~$x,y$ lie ``deep'' inside~$\wt D_N$ and $|x-y|>\delta\sqrt{K_N}=N^{1-\lambda^2+o(1)}$, we have
\begin{equation}
\label{E:3.13ua}
\frakg_x(y)=\frac{G^{\wt D_N}(x,y)}{G^{\wt D_N}(x,x)}
\le \frac{\log\frac{N}{|x-y|}+c}{\log N-c}
\le1-(1-\lambda^2)+o(1)=\lambda^2+o(1),
\end{equation}
where $o(1)\to0$ uniformly in~$x,y\in D_N$. Assuming~$s\in[0,a_N]$, from $\lambda<1/\sqrt 2$ we then have
\begin{equation}
a_N(1-\frakg_x(y))-s\frakg_x(y)>\epsilon a_N
\end{equation}
 for some~$\epsilon>0$ as soon as~$N$ is large enough, uniformly in~$x,y\in D_N$. The argument in Lemma~\ref{lemma-3.1} in conjunction with $\frakg_x(y)\in[0,1]$ and the asymptotic \eqref{E:1.10} then show
\begin{multline}
\qquad
P\Bigl(\,\hat h^{\wt D_N\smallsetminus\{x\}}(y)\ge a_N(1-\frakg_x(y))-s\frakg_x(y)\Bigr)
\\
\le \frac{c}{\sqrt{\log N}}\,\texte^{-\frac{[a_N(1-\frakg_x(y))-s\frakg_x(y)]^2}{2G(y,y)}}
\le c\frac{K_N}{N^2}\,\texte^{\frakg_x(y)\frac{a_N^2}{g\log N}+\frac{a_N}{G(y,y)}\frakg_x(y)s}\,,
\qquad
\end{multline}
where we wrote $G(y,y)$ for~$G^{\wt D_N\smallsetminus\{x\}}(y,y)$ to reduce clutter of indices and then used that $|G(y,y) - g\log N| \leq c$ uniformly in~$y\in D_N$.
The explicit form of the law of~$h^{D_N}$ with respect to the Lebesgue measure readily shows
\begin{equation}
\label{E:3.17e}
P\bigl(h^{\wt D_N}(x)-a_N\in\textd s\bigr)
\le c\frac{K_N}{N^2}\,\texte^{-\frac{a_N}{G(x,x)} s}\textd s.
\end{equation}
Since $G(x,x)/G(y,y)=1+o(1)$ and $\frakg_x(y)\le\lambda^2+o(1)<1$, the integral in~\eqref{E:3.11} over~$s\in[0, a_N]$ yields a harmless multiplicative factor. 
Also, the middle inequality in \eqref{E:3.13ua} implies
\begin{equation}
\label{E:3.18e}
\texte^{\frakg_x(y)\frac{a_N^2}{g\log N}}
\le c\biggl(\frac N{|x-y|}\biggr)^{4\lambda^2+o(1)}
\end{equation}
with~$o(1)\to0$ uniformly in~$x,y\in D_N$ with~$|x-y|>\delta\sqrt{K_N}$. From \eqref{E:3.9ua} we thus get
\begin{equation}
\label{E:3.16}
P\bigl(h^{ D_N}(x)\ge a_N,\,h^{ D_N}(y)\ge a_N\bigr)
\le P\bigl(h^{D_N}(x)\ge 2a_N\bigr)+c\Bigl(\frac {K_N}{N^2}\Bigr)^2\biggl(\frac N{|x-y|}\biggr)^{4\lambda^2+o(1)}
\end{equation}
uniformly in~$x,y\in D_N$ withÊ$|x-y|>\delta\sqrt{K_N}$.

In order to finish the proof, we now write
\begin{equation}
\label{E:3.17}
\begin{aligned}
E \bigl(|\Gamma_N^D (0)|^2 \bigr)
&\le\sum_{\begin{subarray}{c}
x,y\in D_N\\|x-y|\le\delta\sqrt{K_N}
\end{subarray}}
P\bigl(h^{ D_N}(x)\ge a_N\bigr)
\\
&\qquad\qquad\qquad+\sum_{\begin{subarray}{c}
x,y\in D_N\\|x-y|>\delta\sqrt{K_N}
\end{subarray}}
P\bigl(h^{ D_N}(x)\ge a_N,\,h^{ D_N}(y)\ge a_N\bigr).
\end{aligned}
\end{equation}
Summing over~$y$ and invoking Lemma~\ref{lemma-3.1} bounds the first term by a factor of order $(\delta K_N)^2$. The contribution of the first term on the right of \eqref{E:3.16} to the second sum is bounded via Lemma~\ref{lemma-3.1} as well:
\begin{equation}
P\bigl(h^{D_N}(x)\ge 2a_N\bigr)
\le \frac{c}{\sqrt{\log N}}\texte^{-2\frac{a_N^2}{g\log N}}
= c\Bigl(\frac{K_N}{N^2}\Bigr)^2 \texte^{-\frac{a_N^2}{g\log N}}\sqrt{\log N}\le c \delta\Bigl(\frac{K_N}{N^2}\Bigr)^2.
\end{equation}
Plugging in also the second term on the right of \eqref{E:3.16}, we thus get
\begin{equation}
\label{E:3.22e}
E \bigl(|\Gamma_N^D (0)|^2 \bigr)
\le 2c\delta(K_N)^2+c\Bigl(\frac {K_N}{N^2}\Bigr)^2\sum_{\begin{subarray}{c}
x,y\in D_N\\|x-y|>\delta\sqrt{K_N}
\end{subarray}}
\biggl(\frac N{|x-y|}\biggr)^{4\lambda^2+o(1)}\,.
\end{equation}
The standard domination by integrals bounds the sum by $c(N^2)^2\int_{D\times D}|x-y|^{-4\lambda^2+o(1)}\textd x\textd y$ regardless of~$\delta$, with the integral convergent since $\lambda<1/\sqrt2$ implies $4\lambda^2<2$. Hence, also the second term on the right is of order~$(K_N)^2$, thus proving \eqref{e:2.3}. 
\end{proofsect}

\begin{proofsect}{Proof of \eqref{e:2.3b}}
For the limit statement \eqref{e:2.3b}, we will have to reveal the~$D$-dependence of certain constants in the above derivation. The bound \eqref{E:3.22e} shows that we need to consider only pairs~$x,y\in D_N$ with~$|x-y|\ge\epsilon N$ as the contribution of the complementary pairs is negligible in the limit~$N\to\infty$ and~$\epsilon\downarrow0$. We only need to refine the bounds \eqref{E:3.17e} and \eqref{E:3.18e}.
Concerning \eqref{E:3.17e}, the asymptotic \eqref{G-asymp} along with the fact that $\diam(\wt D_N)\le 2\diam(D_N)$ gives
\begin{equation}
G^{\wt D_N}(x,x)\le g\log N+g\log(2\diam D)+c_0+o(1)
\end{equation}
and so, in light of $a_N=O(\log N)$, the constant~$c$ in \eqref{E:3.17e} is bounded by a numerical constant (which comes from bounding such constants in the probability density of~$h^{\wt D_N}(x)$) times
\begin{equation}
\label{E:3.17f}
\texte^{\frac{a_N^2}{g(\log N)^2}\log(\diam D)+o(1)} = (\diam D)^{4\lambda^2+o(1)}
\end{equation}
with~$o(1)\to0$ as~$N\to\infty$. Concerning  \eqref{E:3.18e}, the asymptotic of the Green function \eqref{G-asymp2} at points of distance order~$N$ in turns gives
\begin{equation}
\frakg_x(y)\le \frac1{\log N}\Bigl[-\log\frac{|x-y|}N+\log(2\diam D)+o(1)\Bigr]
\end{equation}
which then implies
\begin{equation}
\label{E:3.25f}
\texte^{\frakg_x(y)\frac{a_N^2}{g\log N}}
\le c\biggl(\frac {N\diam D}{|x-y|}\biggr)^{4\lambda^2+o(1)}\,.
\end{equation}
Using \eqref{E:3.17f} and \eqref{E:3.25f} in the derivation of \eqref{E:3.16} and taking~$N\to\infty$ followed by~$\epsilon\downarrow0$, we then readily get  \eqref{e:2.3b} as well.
\end{proofsect}

\subsection{Subsequential limits and factorization}
We will now start deriving consequences of the above lemmas for the random measures~$\eta_N^D$ from~\eqref{E:2.6}. Since our strategy is to first deal only with events/functions that are trivial in the third ``coordinate,''  we will temporarily abuse notation and set
\begin{equation}
\label{E:3.19}
\eta_N^D:=\frac1{K_N}\sum_{x\in D_N}\delta_{x/N}\otimes\delta_{\,h^{D_N}(x)-a_N}
\end{equation}
instead of the full definition in \eqref{E:2.6}.
As a direct consequence of the above lemmas and the fact that $\overline D\times (\R\cup\{\infty\})$ is a separable metric space, we then get:

\begin{corollary}
\label{cor-3.4}
Suppose $\lambda\in(0,1/\sqrt2)$. Then $\{\eta_N^D\colon N\ge1\}$ is tight with respect to the vague topology on the space of Radon measures on $\overline D\times (\R\cup\{\infty\})$. Moreover, every subsequential weak limit~$\eta^D$ of these measures satisfies, for each~$b\in\R$,
\begin{equation}
\label{E:3.20}
P\Bigl(\eta^D\bigl(D\times[b,\infty)\bigr)<\infty\Bigr)=1
\end{equation}
and, for each non-empty open~$A\subset D$ and each~$b\in\R$,
\begin{equation}
\label{E:3.21}
P\Bigl(\eta^D\bigl(A\times[b,\infty)\bigr)>0\Bigr)>0.
\end{equation}
Furthermore, we have $\eta^D(A\times\R)=0$ a.s.\ for each measurable~$A$ with~$\leb(A)=0$ and, in particular, $\eta^D(\partial D\times\R)=0$ a.s..
\end{corollary}

\begin{proofsect}{Proof}
The first part of the statement requires showing that, for any continuous compactly-supported function $f\colon \overline D\times (\R\cup\{\infty\})\to\R$, the family of random variables $\{\langle \eta_N^D,f\rangle\colon N\ge1\}$ is tight. For this it suffices to show that the family $\{\eta_N^D(\overline D\times[b,\infty))\colon N\ge1\}$ is tight for each~$b\in\R$. This is a consequence of Lemma~\ref{lemma-3.1} and the fact that $\eta_N^D(\overline D\times[b,\infty)) = |\Gamma_N^D(b)|$.

Let now~$\eta^D$ be a subsequential weak limit of the measures $\{\eta_N^D\colon N\ge1\}$. Fatou's lemma, a straightforward approximation argument and Lemma~\ref{lemma-3.1} then show $E\eta^D(\overline D\times[b,\infty))<\infty$ for each~$b\in\R$. Lemma~\ref{lemma-3.1} also gives $\eta^D(A\times\R)=0$ a.s.\ whenever $\leb(A)=0$.  It remains to show that~$\eta^D$ is non-trivial in the sense stated in \eqref{E:3.21}. Let $A\subset D$ be non-empty and open and pick~$b\in\R$. 
 Denote $X_N:=\eta_N^D(A\times[b,\infty))$. Lemma~\ref{lemma-3.3} shows that $\sup_{N\ge1}E(X_N^2)<\infty$ and so the family~$\{X_N\colon N\ge1\}$ is uniformly integrable. Since $\inf_{N\ge1} E(X_N)>0$ by Lemma~\ref{lemma-3.2} and the fact that~$\psi_\lambda^D>0$, any distributional limit of~$X_N$ has positive expectation as well. 
\end{proofsect}

Given a function $f\colon \overline D\times (\R\cup\{\infty\})\to\R$ and~$b\in\R$, define
\begin{equation}
f_b(x,h):=f(x,h+b)\texte^{-\alpha\lambda b}.
\end{equation}
A key step is now the proof of:

\begin{proposition}
\label{prop-factor}
Suppose $\lambda\in(0,1/\sqrt2)$.
Any subsequential limit~$\eta^D$ of $\{\eta_N^D\colon N\ge1\}$ obeys  the following: For each~$b\in\R$ and each $f\colon \overline D\times (\R\cup\{\infty\})\to\R$ of the form $f(x,h)=\1_A(x)\1_{[0,\infty)}(h)$ with $A\subset D$ open,
\begin{equation}
\label{E:3.24}
\langle\eta^D,f_b\rangle = \langle \eta^D,f\rangle \in\R
\end{equation}
holds with probability one. 
\end{proposition}

The proof of this proposition relies on a calculation that is formalized as:

\begin{lemma}
\label{lemma-aux}
For any $\lambda\in(0,1/\sqrt2)$, any open $A\subset D$, any $b\in\R$, and $A_N:=\{x\in \Z^2\colon x/N\in A\}$,
\begin{equation}
\lim_{N\to\infty}\,
\frac1{K_N}E\Bigl|\bigl|\Gamma_N^D(0)\cap A_N\bigr|-\texte^{\alpha\lambda b}\bigl|\Gamma_N^D(b)\cap A_N\bigr|\Bigr|=0.
\end{equation}
\end{lemma}

\begin{proofsect}{Proof}
Since any open set $A \subset D$ can be written as the union of an increasing sequence of open sets whose closure lies inside $D$, we can assume without loss of generality that $A$ has positive Euclidean distance to $D^{\text c}$. Then, using Cauchy-Schwarz, we may as well show
\begin{equation}
\lim_{N\to\infty}\,
\frac1{K_N^2}E\biggl(\Bigl(\bigl|\Gamma_N^D(0)\cap A_N\bigr|-\texte^{\alpha\lambda b}\bigl|\Gamma_N^D(b)\cap A_N\bigr|\Bigr)^2\biggr)=0.
\end{equation}
 Invoking \eqref{E:3.9a}, the calculation in the proof of Lemma~\ref{lemma-3.3} shows that the second moment of $|\Gamma_N^D(b)|$ is dominated by the sum of pairs $x,y\in D_N$ that are within distance of order~$N$.  It thus suffices to prove that, for any $\epsilon>0$,
\begin{equation}
\label{E:3.27}
\begin{aligned}
\max_{\begin{subarray}{c}
x,y\in A_N\\|x-y|\ge\epsilon N
\end{subarray}}
\biggl(P\bigl(h^{D_N}&(x)\ge a_N,\,h^{D_N}(y)\ge a_N\bigr)
\\*[-4mm]
&-\texte^{\alpha\lambda b}P\bigl(h^{D_N}(x)\ge a_N+b,\,h^{D_N}(y)\ge a_N\bigr)
\\*[1mm]
&\qquad-\texte^{\alpha\lambda b}P\bigl(h^{D_N}(x)\ge a_N,\,h^{D_N}(y)\ge a_N+b\bigr)
\\
&\qquad\qquad
+
\texte^{2\alpha\lambda b}P\bigl(h^{D_N}(x)\ge a_N+b,\,h^{D_N}(y)\ge a_N+b\bigr)\biggr)
= o\Bigl(\frac{K_N^2}{N^4}\Bigr)
\end{aligned}
\end{equation}
as~$N\to\infty$. For this we need to compute a sharp leading order asymptotic of the probability $P(h^{D_N}(x)\ge a_N+b_1,\,h^{D_N}(y)\ge a_N+b_2)$ for the four possible choices $b_1,b_2\in\{0,b\}$.

We will invoke the decomposition \eqref{E:3.10} and the representation \eqref{E:3.11}. Writing again~$G(y,y)$ for~$G^{D_N\smallsetminus\{x\}}(y,y)$, for any $s\in[0,a_N+b_2]$ we then get
\begin{equation}
\label{E:3.31u}
\begin{aligned}
P\Bigl(\,\hat h^{D_N\smallsetminus\{x\}}&(y)\ge a_N(1-\frakg_x(y))-s\frakg_x(y)+b_1\Bigr)
\\
&= \frac{c+o(1)}{\sqrt{\log N}}\texte^{-\frac{[a_N(1-\frakg_x(y))-s\frakg_x(y)+b_1]^2}{2G(y,y)}}
= \frac{c+o(1)}{\sqrt{\log N}}\texte^{-\frac{[a_N(1-\frakg_x(y))-s\frakg_x(y)]^2}{2G(y,y)}}\texte^{-\alpha\lambda b_1}
\\
&=\bigl(\texte^{-\alpha\lambda b_1}+o(1)\bigr)P\Bigl(\,\hat h^{D_N\smallsetminus\{x\}}(y)\ge a_N(1-\frakg_x(y))-s\frakg_x(y)\Bigr)\,,
\end{aligned}
\end{equation}
where $c>0$ is a numerical constant and where we used that $\frakg_x(y)=O(1/\log N)$ when $|x-y|\ge\epsilon N$ and then applied the asymptotic of~$a_N$ and~$G^{D_N\smallsetminus\{x\}}(y,y)=g\log N+O(1)$ enabled by the fact that now~$y$ is ``deep'' inside~$D_N\smallsetminus\{x\}$ as implied by the assumptions on $A$ and $x,y$. In addition, writing $G(x,x)$ for~$G^{D_N}(x,x)$, we also get
\begin{multline}
\label{E:3.32u}
\quad
P\bigl(h^{D_N}(x)-a_N-b_2\in\textd s\bigr)
= \frac{c+o(1)}{\sqrt{\log N}}\texte^{-\frac{[a_N+b_2+s]^2}{2G(x,x)}}\textd s
\\
=\texte^{-\alpha\lambda b_2}\,\frac{c+o(1)}{\sqrt{\log N}}\texte^{-\frac{[a_N+s]^2}{2G(x,x)}}\textd s
=\bigl(\texte^{-\alpha\lambda b_2}+o(1)\bigr)P\bigl(h^{D_N}(x)-a_N\in\textd s\bigr)\,.
\quad
\end{multline}
where~$c$ is again a positive constant and~$o(1)\to0$ as~$N\to\infty$ uniformly in~$s\in[0,a_N]$.
Putting \twoeqref{E:3.31u}{E:3.32u} together and integrating over~$s\in[0,a_N]$ we get 
\begin{multline}
\label{E:3.32ua}
\qquad
P\bigl(h^{D_N}(x)\ge a_N+b_1,\,h^{D_N}(y)\ge a_N+b_2\bigr)
\\=\bigl(\texte^{-\alpha\lambda (b_1+b_2)}+o(1)\bigr)P\bigl(h^{D_N}(x)\ge a_N,\,h^{D_N}(y)\ge a_N\bigr)\,,
\qquad
\end{multline}
where we used that, by the FKG inequality and Lemma~\ref{lemma-3.2},
\begin{equation}
P\bigl(h^{D_N}(x)\ge 2a_N\bigr)=o(1) P\bigl(h^{D_N}(x)\ge a_N,\,h^{D_N}(y)\ge a_N\bigr)
\end{equation}
with $o(1)\to0$ as~$N\to\infty$ uniformly in~$x,y\in A_N$ with~$|x-y|\ge\epsilon N$.
Plugging \eqref{E:3.32ua} into \eqref{E:3.27}, the desired claim follows. 
\end{proofsect}

We are now ready to give:

\begin{proofsect}{Proof of Proposition~\ref{prop-factor}}
Let $f(x,h):=\1_A(x)\1_{[0,\infty)}(h)$ with $A\subset D$ closed. Lemma~\ref{lemma-aux} can be rephrased as
\begin{equation}
\lim_{N\to\infty}E\bigl|\langle \eta^D_N,f\rangle-\langle \eta^D_N,f_{b}\rangle\bigr|=0\,,\qquad b\in\R.
\end{equation}
Taking the distributional limit (choosing further subsequence if necessary) of $\langle \eta^D_N,f-f_{b}\rangle$ then shows,  by Fatou's lemma,  $\langle \eta^D,f-f_{b}\rangle=0$ a.s.\ which is the desired claim. (The finiteness follows from the tightness proved in Lemma~\ref{lemma-3.1}.)
\end{proofsect}

The identity \eqref{E:3.24} now readily yields the desired factorization property:

\begin{lemma}[Factorization]
\label{lemma-3.7}
 Suppose $\eta^D$ is a Borel measure on $\overline D\times\R$ with $\eta^D(\partial D \times \R) = 0$ a.s.\  such that \eqref{E:3.24} holds for each function $f\colon \overline D\times \R\to\R$ of the form $f(x,h):=\1_A(x)\1_{[0,\infty)}(h)$ with~$A\subset D$ open and each~$b\in\R$. Then, with probability one,~$\eta^D$ takes the form
\begin{equation}
\label{E:3.31}
\eta^D(\textd x\textd h)=Z^D_\lambda(\textd x)\otimes \texte^{-\alpha\lambda h}\textd h
\end{equation}
for some finite random Borel measure $Z^D_\lambda$ on~$\overline D$.
\end{lemma}

\begin{proofsect}{Proof}
For $A\subset \ol{D}$ Borel, define $Z^D_\lambda(A):=\alpha\lambda\eta^D(A\times[0,\infty))$. This is automatically a finite Borel measure on~$D$. The condition \eqref{E:3.24} now shows that, almost surely if $A$ is open then
\begin{equation}
\begin{aligned}
\eta^D\bigl(\,A\times[b,\infty)\bigr)
& =\texte^{-\alpha\lambda b}\langle\eta^D,f_{b}\rangle= \texte^{-\alpha\lambda b}\langle\eta^D,f\rangle 
\\
&=(\alpha\lambda)^{-1}\texte^{-\alpha\lambda b}Z^D_\lambda(A)
=\int_{A\times[b,\infty)}Z^D_\lambda(\textd x) \texte^{-\alpha\lambda h}\textd h.
\end{aligned}
\end{equation}
The null set in this statement may depend on~$A$ and~$b$ but we can choose a common null set for all sets in the class $\{A\times[b,\infty)\colon A\subset D \text{ open dyadic square},\, b\in\Q\}$ as it is countable. The equality of the measures \eqref{E:3.31} on $D$ then follows from the fact that this class is a $\pi$-system (in the sense of Dynkin) which generates the product Borel $\sigma$-algebra on $D\times\R$. As $\eta^D(\partial D \times R) = Z^D_\lambda(\partial D) = 0$, the equality in~\eqref{E:3.31} extends to all of $\ol{D}$.
\end{proofsect}


\subsection{Uniqueness of subsequential limit}
At this point we have shown that, assuming $\lambda\in(0,1/\sqrt2)$, every subsequential limit~$\eta^D$ of the measures $\{\eta^D_N\colon N\ge1\}$, with $\eta^D_N$ as in \eqref{E:3.19}, factors into the form \eqref{E:3.31}. The goal of this subsection is to show that the measure~$Z^D_\lambda$, and thus also the subsequential limit~$\eta^D$, is in fact unique. This will in particular show that~$\eta^D_N$ converges in distribution to the same limit, regardless of the approximating sequence~$D_N$ or the way~$a_N$ achieves the asymptotic \eqref{E:1.10}.

As our first lemma we will check that $Z^D_\lambda$ obeys the properties listed in Theorem~\ref{thm-2.3}. However, these require extracting subsequential limits for multiple domains at the same time.  Cantor's diagonal argument makes this possible provided we restrict ourselves to  a countable class $\mathfrak D_0$ of domains in~$\mathfrak D$. We will assume that $\mathfrak D_0$ contains all open squares of the form
\begin{equation}
\label{E:3.33}
\bigl(k2^{-n},(k+1)2^{-n}\bigr)\times\bigl(\ell 2^{-n},(\ell+1)2^{-n}\bigr),\qquad k,\ell,n\in\Z.
\end{equation}
For each domain~$D\in\mathfrak D_0$ we then fix a sequence~$\{D_N\}$ of lattice approximations satisfying \twoeqref{E:2.1}{E:2.2}. All (simultaneous) subsequential limits will naturally pertain to the specific choice of~$\mathfrak D_0$ as well as the lattice approximations~$\{D_N\}$.

\begin{proposition}
\label{lemma-3.8}
Assume $\lambda\in(0,1/\sqrt2)$ and let~$\{\eta^D\colon D\in\mathfrak D_0\}$ be subsequential limits (along the same subsequence) of $\{\eta^D_N\colon N\ge1\}$ for~$D\in\mathfrak D_0$. For each~$D\in\mathfrak D_0$, let~$Z^D_\lambda$ be the measure associated with~$\eta^D$ as in \eqref{E:3.31}. Then $\{Z^D_\lambda\colon D\in\mathfrak D_0\}$ obeys properties (1-7) in~Theorem~\ref{thm-2.3}, with all domains restricted to be contained in~$\mathfrak D_0$ or translates thereof.
\end{proposition}

\begin{proofsect}{Proof of properties (1-6)}
Properties~(1,2) are direct consequences of Lemma~\ref{lemma-3.1}. Property~(3) holds for all open $A\subset D$  thanks to Lemma~\ref{lemma-3.2}; the equality for general measurable $A\subset D$ is then obtained by realizing that \eqref{E:2.12} represents equality of two Borel measures. Property~(4) is a consequence of the representation of $\eta^{D\cup\wt D}_N$ as the sum of independent copies of $\eta^D_N$ and $\eta^{\wt D}_N$ --- which itself follows by representing $h^{D_N\cup\wt D_N}$ as the sum of independent fields $h^{D_N}$ and~$h^{\wt D_N}$. The translation invariance in property~(6) is immediate. 

Concerning property~(5), let $D,\wt D\in\mathfrak D_0$ with $\wt D\subset D$ and $\leb(D\smallsetminus\wt D)=0$. The Gibbs-Markov decomposition of the DGFF (see~\eqref{E:A.6}) then yields
\begin{equation}
h^{D_N} \laweq h^{\wt D_N}+\varphi^{D_N,\wt D_N}.
\end{equation}
This means that if~$f\colon\overline D\times\R\to\R$ is continuous with compact support in~$\wt D$, then
\begin{equation}
\label{E:3.35}
\langle\eta^D_N,f\rangle \laweq\langle\eta^{\wt D}_N,f_\varphi\rangle
\end{equation}
where
\begin{equation}
f_\varphi(x,h) := f\bigl(x,h+\varphi^{D_N,\wt D_N}(\lfloor xN\rfloor)\bigr)
\end{equation}
with $\varphi^{D_N,\wt D_N}$ independent of $\eta^{\wt D}_N$ on the right-hand side of \eqref{E:3.35}.
As shown in \cite[Lemma~B.14]{BL3}, for each~$N\ge1$ and each~$\delta>0$, there is a coupling of~$\varphi^{D_N,\wt D_N}(\lfloor \cdot N\rfloor)$ with~$\Phi^{D,\wt D}$ such that the supremum of the difference on~$D^\delta$ tends to zero in probability; see~\eqref{E:A.8}. Thanks to continuity and restriction on the support of~$f$, we thus have
\begin{equation}
\label{E:3.38}
\langle\eta^{\wt D}_N,f_\varphi\rangle \laweq \langle\eta^{\wt D}_N,f_\Phi\rangle+o(1)
\end{equation}
where $o(1)\to0$ in probability (as $N\to\infty$) and where
\begin{equation}
f_\Phi(x,h):=f\bigl(x,h+\Phi^{D,\wt D}(x)\bigr)
\end{equation}
with $\Phi^{D,\wt D}$ independent of~$\eta^{\wt D}_N$ on the right-hand side of \eqref{E:3.38}. Since $x\mapsto\Phi^{D,\wt D}(x)$ is continuous on~$\wt D$ a.s., for any simultaneous subsequential limits $\eta^D$ of $\{\eta^D_N\colon N\ge1\}$ and $\eta^{\wt D}$ of $\{\eta^{\wt D}_N\colon N\ge1\}$, we thus obtain
\begin{equation}
\langle \eta^D,f\rangle \laweq \langle\eta^{\wt D},f_\Phi\rangle\,,
\end{equation}
where $\Phi^{D,\wt D}$ (implicitly contained in $f_\Phi$) is independent of $\eta^{\wt D}$ on the right-hand side. But the representation \eqref{E:3.31} now permits us to write
\begin{equation}
\begin{aligned}
\langle\eta^{\wt D},f_\Phi\rangle 
&= \int_{D\times\R}Z^{\wt D}_\lambda(\textd x)\texte^{-\alpha\lambda h}\textd h \,\,f\bigl(x,h+\Phi^{D,\wt D}(x)\bigr)
\\
&=\int_{D\times\R}Z^{\wt D}_\lambda(\textd x)\texte^{-\alpha\lambda (h-\Phi^{D,\wt D}(x))}\textd h \,\,f(x,h).
\end{aligned}
\end{equation}
As this holds for any continuous $f\colon D\times\R\to\R$ with support in~$\wt D$, and since both $Z_\lambda^{\wt D}$ and $Z^D_\lambda$ assign zero mass to $D\smallsetminus\wt D$ due to the fact that $\leb(D\smallsetminus\wt D)=0$, property~(5) follows.
\end{proofsect}

For property~(7), and also later use, we will need:

\begin{lemma}
\label{lemma-3.9a}
For each~$\lambda\in(0,1/\sqrt2)$ there is~$c\in(0,\infty)$ such that for any open square~$S\subset\C$
\begin{equation}
E\bigl[Z^S_\lambda(S)^2\bigr]\le c \bigl[E Z_\lambda^S(S)\bigr]^2.
\end{equation}
\end{lemma}

\begin{proofsect}{Proof}
Suppose~$S$ is a translation (and rotation) of~$(0,r)^2$. Then \eqref{e:2.3b} in Lemma~\ref{lemma-3.3} (with the help of Fatou's lemma) and a simple scaling argument show that $E[Z^S_\lambda(S)^2]\le c r^{4+4\lambda^2}$ for some constant~$c$ independent of~$r$. On the other hand, Lemma~\ref{lemma-3.2} along with uniform integrability of normalized level-set sizes and the fact that $\psi_\lambda^{rD}(rx) = r^{2\lambda^2}\psi_\lambda^D(x)$ show $E[Z^S_\lambda(S)]\ge \tilde c r^{2+2\lambda^2}$ for some absolute~$\tilde c>0$. The claim follows.
\end{proofsect}

\begin{proofsect}{Proof of property~(7)}
It suffices to prove this for all squares of the form \eqref{E:3.33} as each open set  contains at least one such square. (Note that~$\mathfrak D_0$ contains all these squares.) For $n\in\Z$, let $p_n:=P(Z_\lambda^{S_n}>0)$ where (appealing to translation invariance), $S_n:=(0,2^{-n})$.  The second moment estimate in conjunction with Lemma~\ref{lemma-3.9a}  show that $p:=\inf_{n\in\Z}p_n>0$. Decomposing~$S_n$ into $4^m$ translates $S_{n+m,1},\dots,S_{n+m,4^m}$ of the square~$S_{n+m}$, the Gibbs-Markov property yields
\begin{equation}
Z^{S_n}_\lambda(S_n) \,\laweq\, \sum_{i=1}^{4^m}\int_{S_{n+m,i}}Z^{S_{n+m,i}}_\lambda(\textd x)\,\,\texte^{\alpha\lambda\Phi^{S_n,\wt S_{n,m}}(x)},
\end{equation}
where $\wt S_{n,m}:=\bigcup_{i=1}^{4^m} S_{n+m,i}$. Then $Z^{S_n}_\lambda(S_n)=0$ forces $Z^{S_{n+m,i}}_\lambda(S_{n+m,i})=0$ for all $i=1,\dots,{4^m}$. Since the latter measures are independent of one another as well as of the field~$\Phi^{S_n,\wt S_{n,m}}$, we have
\begin{equation}
1-p_n\le (1-p_{n+m})^{4^m}\le(1-p)^{4^m}.
\end{equation}
Taking $m\to\infty$ we get $p_n=1$ for each $n\in\Z$ and so property~(6) follows.
\end{proofsect}

Using the same notation as in the previous proof, in order to prove uniqueness, we will first characterize $Z^{S_n}_\lambda$ as the limit of the measures
\begin{equation}
Y^{S_n}_m(\textd x):=c\sum_{i=1}^{4^m}\texte^{\alpha\lambda\Phi^{S_n,\wt S_{n,m}}(x)}\psi_\lambda^{S_{n+m,i}}(x)\1_{S_{n+m,i}}(x)\,\textd x,
\end{equation}
where~$c$ is the constant from \eqref{E:2.12}. Indeed, we have:

\begin{lemma}
\label{lemma-conv}
For each $\lambda\ge0$, there exists an a.s.\ finite random measure $Y^{S_n}_\infty$ (possibly degenerate to zero), such that for each bounded, measurable $f\colon \overline D\to\R$,
\begin{equation}
\langle Y^{S_n}_m,f\rangle\,\,\underset{m\to\infty}\longrightarrow\,\, \langle Y^{S_n}_\infty,f\rangle,\qquad \text{\rm a.s.}
\end{equation}
\end{lemma}

\begin{proofsect}{Proof}
Thanks to the structure of the covariances \eqref{E:2.10}, we can write $\Phi^{S_n,\wt S_{n,m}}(x)$ as the sum of independent fields
\begin{equation}
\Phi^{S_n,\wt S_{n,m}}(x)  = \sum_{j=1}^m \Phi^{\wt S_{n,j-1},\wt S_{n,j}}(x),
\end{equation}
where $\wt S_{n,0}:=S_n$. In light of the fact that
\begin{equation}
\sum_{j=1}^{4^m}\psi^{S_{n+m,j}}_\lambda(x)\1_{S_{n+m,j}}(x) = \psi^{\wt S_{n,m}}_\lambda(x),
\qquad x\in \wt S_{n,m}
\end{equation}
and that, for any $\wt D\subset D$,
\begin{equation}
E\texte^{\alpha\lambda\Phi^{D,\wt D}(x)} = \texte^{\frac12\alpha^2\lambda^2 C^{D,\wt D}(x)} = \frac{\psi^D_\lambda(x)}{\psi^{\wt D}_\lambda(x)},\qquad x\in\wt D,
\end{equation}
a straightforward calculation shows that $\{\langle Y^{S_n}_m,f\rangle\colon m\ge1\}$ is a martingale with respect to the filtration
\begin{equation}
\FF_m:=\sigma\Bigl(\Phi^{\wt S_{n,j-1},\wt S_{n,j}}(x)\colon x\in S_n',\,j=1,\dots,m\Bigr)
\quad\text{where}\quad
S_n':=\bigcap_{m\ge1}\wt S_{n,m}.
\end{equation}
(Note that, as $\leb(S_n\smallsetminus S_n')=0$, restricting the measures to~$S_n'$ carries no loss.) Since for $f\ge0$ the martingale is non-negative, the Martingale Convergence Theorem shows that
\begin{equation}
L(f):=\lim_{m\to\infty}\langle Y^{S_n}_m,f\rangle
\end{equation}
exists almost surely for each bounded measurable~$f$ (we treat the positive and the negative part of~$f$ separately). The null set in this statement may depend on~$f$.

In order to show that the limit is an integral of~$f$ with respect to a random measure, we follow a standard argument from multiplicative chaos theory: Fix a countable dense subset $A\subset C(\overline D)$. Fatou's lemma yields
\begin{equation}
\label{E:3.50}
E \bigl|L(f)\bigr|\le c\int_{S_n}|f(x)|\psi_\lambda^D(x)\textd x,\qquad f\in C(\overline D)
\end{equation}
and so by the Markov inequality, on a set of full probability, the linear functional $f\mapsto L(f)$ is well-defined for all~$f\in A$ simultaneously and bounded on~$A$ in the supremum norm. It follows that $f\mapsto L(f)$ extends uniquely to an almost-everywhere defined continuous linear functional $f\mapsto\overline L(f)$ on~$C(\overline D)$ such that $L(f)=\overline L(f)$ holds almost surely for each~$f\in C(\overline D)$. (The null set of this equality may depend on~$f$.) The Riesz Representation Theorem then readily gives existence of a Borel measure~$Y_\infty^{S_n}$ on~$\overline D$ such that $ L(f) = \langle Y_\infty^{S_n},f\rangle$ holds almost surely for each $f\in C(\overline D)$. The measure is finite a.s.\ (albeit possibly trivially zero) thanks to \eqref{E:3.50}.
\end{proofsect}

The desired uniqueness of the law of~$Z^D_\lambda$ will now follow from:

\begin{proposition}
\label{prop-unique}
Suppose $\lambda\in(0,1/\sqrt2)$. Then for each $n\in\Z$,
\begin{equation}
\label{E:3.51}
Z^{S_n}_\lambda(\textd x)\,\laweq Y^{S_n}_\infty(\textd x).
\end{equation}
\end{proposition}

 For the proof of this proposition, let 
$f\colon S_n\to[0,\infty)$ be a bounded, measurable function. Our aim is to prove that
\begin{equation}
\label{E:3.52}
E\texte^{-\langle Z^{S_n}_\lambda,f\rangle} = E\texte^{-\langle Y^{S_n}_\infty,f\rangle}
\end{equation}
We will do this by separately proving~$\ge$ and~$\le$.

\begin{proofsect}{Proof of $\ge$ in \eqref{E:3.52}}
Thanks to the Gibbs-Markov property, we may represent $Z^{S_n}_\lambda$ as
\begin{equation}
\label{E:3.53}
Z^{S_n}_\lambda(\textd x)=
\sum_{i=1}^{4^m}\texte^{\alpha\lambda\Phi^{S_n,\wt S_{n,m}}(x)}\1_{S_{n+m,i}}(x)Z^{S_{n+m,i}}_\lambda(\textd x),
\end{equation}
where $Z^{S_{n+m,i}}_\lambda$, $i=1,\dots,{4^m}$, are independent of one another as well as of~$\Phi^{S_n,\wt S_{n,m}}$ on the right-hand side.  In light of \eqref{E:2.12}
we have
\begin{equation}
\label{E:3.54}
E\Bigl(\langle Z^{S_n}_\lambda,f\rangle\,\Big|\,\sigma\bigr(\Phi^{S_n,\wt S_{n,m}}\bigl)\Bigr)
=\langle Y^{S_n}_m,f\rangle.
\end{equation}
Jensen's inequality then shows
\begin{equation}
\label{E:3.55}
E\texte^{-\langle Z^{S_n}_\lambda,f\rangle}
\ge 
E\texte^{-\langle Y^{S_n}_m,f\rangle}
\,\underset{m\to\infty}\longrightarrow\,E\texte^{-\langle Y^{S_n}_\infty,f\rangle},
\end{equation}
where the limit uses Lemma~\ref{lemma-conv} and the Bounded Convergence Theorem.
\end{proofsect}

For the opposite bound, we will need the following ``reverse Jensen'' inequality:

\begin{lemma}
\label{lemma-3.11a}
If $X_1,\dots,X_n$ are non-negative independent random variables, then for each~$\epsilon>0$,
\begin{equation}
\label{E:3.56}
E\biggl(\exp\Bigl\{-\sum_{i=1}^n X_i\Bigr\}\biggr)\le \exp\Bigl\{-\texte^{-\epsilon}\sum_{i=1}^n E(X_i\,;\, X_i\le \epsilon)\Bigr\}\,.
\end{equation}
\end{lemma}

\begin{proofsect}{Proof}
Thanks to independence, it suffices to prove this for~$n=1$. This is checked by bounding $E(\texte^{-X})\le E(\texte^{-\wt X})$, where $\wt X :=X\1_{\{X\le\epsilon\}}$, writing 
\begin{equation}
-\log E(\texte^{-\wt X})=\int_0^1\textd s\,\, \frac{E(\wt X\texte^{-s\wt X})}{E(\texte^{-s \wt X})}
\end{equation}
and invoking the bounds $E(\wt X\texte^{-s\wt X})\ge \texte^{-\epsilon} E(\wt X)$ and $E(\texte^{-s \wt X})\le1$.
\end{proofsect}

We will also need to invoke an additional truncation: For $\delta\in(0,1/2)$, let $S_k^\delta$ be the translate of $(\delta 2^{-k},(1-\delta)2^{-k})$ centered at the same point as~$S_k$. Analogously, let $S_{n+m,i}^\delta$ be the corresponding truncation of~$S_{n+m,i}$. Then set $\wt S_{n,m}^\delta:=\bigcup_{i=1}^{4^m} S_{n+m,i}^\delta$ and let 
\begin{equation}
f_{\delta}(x):=\1_{\wt S_{n,m}^\delta}(x) f(x).
\end{equation}
Positivity of $f$ implies that $f-f_\delta\ge0$ with $f_\delta\uparrow f$ (as~$\delta\downarrow0$) Lebesgue almost everywhere, the Monotone Convergence Theorem and \eqref{E:3.50}) show $E\langle Y^{S_n}_\infty,f-f_\delta\rangle\to0$ as~$\delta\downarrow0$. In particular,
\begin{equation}
\label{E:3.57}
\langle Y^{S_n}_\infty,f_\delta\rangle\,\underset{\delta\downarrow0}\longrightarrow\,
\langle Y^{S_n}_\infty,f\rangle,\qquad\text{a.s.}
\end{equation}
This will permit us to work with~$f_\delta$ instead of~$f$. Next we will need:

\begin{lemma}
\label{lemma-3.11}
Suppose $\lambda\in(0,1/\sqrt2)$. Given~$\delta\in(0,1/2)$, for each $i=1,\dots,{4^m}$ abbreviate
\begin{equation}
\label{E:3.59}
X_i:=\int_{S_{n+m,i}^\delta}\texte^{\alpha\lambda\Phi^{S_n,\wt S_{n,m}}(x)}f_\delta(x) Z^{S_{n+m,i}}_\lambda(\textd x).
\end{equation}
Then for each~$\epsilon>0$,
\begin{equation}
\label{E:3.59a}
\sum_{i=1}^{4^m} E(X_i\,;\,X_i>\epsilon)\,\underset{m\to\infty}\longrightarrow\,0.
\end{equation}
\end{lemma}

\begin{proofsect}{Proof}
Abbreviate $L:=2^m$ throughout this proof. Since $E(X_i\,;\,X_i>\epsilon)\le \frac1\epsilon E(X_i^2)$, we can bound the sum in \eqref{E:3.59a} by
\begin{equation}
\label{E:3.60}
\frac1\epsilon \sum_{i=1}^{L^2} E(X_i^2)
\le\frac{\Vert f\Vert^2}\epsilon\sum_{i=1}^{L^2} E\int_{S_{n+m,i}^\delta\times S_{n+m,i}^\delta}
\texte^{\alpha\lambda[\Phi^{S_n,\wt S_{n,m}}(x)+\Phi^{S_n,\wt S_{n,m}}(y)]} Z^{S_{n+m,i}}_\lambda(\textd x)
Z^{S_{n+m,i}}_\lambda(\textd y).
\end{equation}
Since $Z^{S_{n+m,i}}_\lambda$, $i=1,\dots,{4^m}$, are independent of the field, we will now take conditional expectation given these measures and invoking that
\begin{equation}
\begin{aligned}
E\texte^{\alpha\lambda[\Phi^{S_n,\wt S_{n,m}}(x)+\Phi^{S_n,\wt S_{n,m}}(y)]}
& = \texte^{\frac12\alpha^2\lambda^2\Var(\Phi^{S_n,\wt S_{n,m}}(x)+\Phi^{S_n,\wt S_{n,m}}(y))}
\\
&\le c \texte^{\frac12\alpha^2\lambda^2 4g\log(L)} = c L^{8\lambda^2}
\end{aligned}
\end{equation}
for some constant~$c$ that arises from the uniform bound $\Var(\Phi^{S_n,\wt S_{n,m}}(x))\le c'+g\log({2^m})$ valid with the same constant~$c'$ for all~$x\in \wt S_{n,m}^\delta$. The right-hand side of \eqref{E:3.60} is thus at most
\begin{equation}
c\,  \frac{\Vert f\Vert^2}{\epsilon}\,L^{8\lambda^2}\sum_{i=1}^{L^2} E\bigl(Z_\lambda^{S_{n+m,i}}(S_{n+m,i})^2\bigr).
\end{equation}
Under the condition $\lambda\in(0,1/\sqrt2)$, Lemma~\ref{lemma-3.9a} implies
\begin{equation}
E\bigl(Z_\lambda^{S_{n+m,i}}(S_{n+m,i})^2\bigr)\le c\, \bigl[E\bigl(Z_\lambda^{S_{n+m,i}}(S_{n+m,i})\bigr)\bigr]^2
\end{equation}
and a trivial scaling argument applied to \eqref{E:2.12} shows
\begin{equation}
\label{E:3.65}
E\bigl(Z_\lambda^{S_{n+m,i}}(S_{n+m,i})\bigr)\le c L^{-(2+2\lambda^2)} \,,
\end{equation}
where $c$ depends on $n$, which is fixed throughout the proof. These observations yield
\begin{equation}
\sum_{i=1}^{L^2} E(X_i\,;\,X_i>\epsilon)
\le\frac c\epsilon \Vert f\Vert^2\,L^{8\lambda^2-2-4\lambda^2} = \frac c\epsilon \Vert f\Vert^2\,L^{-2(1-2\lambda^2)}.
\end{equation}
This tends to zero as~$m\to\infty$ for all $\lambda<(0,1/{\sqrt2})$ as claimed.
\end{proofsect}

\begin{proofsect}{Proof of $\le$ in \eqref{E:3.52}}
Consider the $\sigma$-algebra $\FF:=\sigma(\Phi^{S_n,\wt S_{n,m}}(x)\colon x\in \wt S_{n,m})$ and note that, for~$X_i$ as in \eqref{E:3.59}, we have $\langle Y^{S_n}_m,f_\delta\rangle=\sum_{i=1}^{4^m}E(X_i|\FF)$.  Fix~$\epsilon>0$. The ``reverse Jensen'' inequality \eqref{E:3.56} for the conditional expectation given~$\FF$ then yields
\begin{equation}
\label{E:3.66}
E\texte^{-\langle Z^{S_n}_\lambda,f\rangle}
\le E\texte^{-\langle Z^{S_n}_\lambda,f_\delta\rangle}
\le E\Biggl(\exp\biggl\{-\texte^{-\epsilon}\Bigl[\langle Y^{S_n}_m,f_\delta\rangle-\sum_{i=1}^{4^m} E(X_i\1_{\{X_i>\epsilon\}}|\FF)\Bigr]\biggr\}\Biggr).
\end{equation}
Lemma~\ref{lemma-conv} gives $\langle Y^{S_n}_m,f_\delta\rangle\to\langle Y^{S_n}_\infty,f_\delta\rangle$ almost surely, while Lemma~\ref{lemma-3.11} shows
\begin{equation}
\sum_{i=1}^{4^m} E(X_i\1_{\{X_i>\epsilon\}}|\FF)\,\underset{m\to\infty}\longrightarrow\,0,\quad\text{in probability}.
\end{equation}
Since the square bracket on the right-hand side of \eqref{E:3.66} is non-negative, taking $m\to\infty$ with the help of the Bounded Convergence Theorem followed by $\epsilon\downarrow0$ thus yields
\begin{equation}
E\texte^{-\langle Z^{S_n}_\lambda,f\rangle}
\le E\texte^{-\langle Y^{S_n}_\infty,f_\delta\rangle}.
\end{equation}
From here $\le$ in \eqref{E:3.52} follows in light of the observation \eqref{E:3.57}.
\end{proofsect}

\begin{proofsect}{Proof of Proposition~\ref{prop-unique}}
Thanks to the fact that the Laplace transform determines the law for non-negative random variables, \eqref{E:3.52} implies $\langle Z^{S_n}_\lambda,f\rangle \laweq\langle Y^{S_n}_\infty,f\rangle$ for each bounded, measurable~$f$. This is what is represented by \eqref{E:3.51}.
\end{proofsect}

We are now ready to summarize our conclusions in:

\begin{theorem}
\label{thm-3.13}
For each $\lambda\in(0,1/\sqrt2)$ and each~$D\in\mathfrak D$ there is a random Borel measure~$Z^D_\lambda$ on~$D$  such that the following holds for each~$a_N$ satisfying \eqref{E:1.10} and each sequence~$\{D_N\}$ of lattice domains satisfying \twoeqref{E:2.1}{E:2.2}: The family of measures $\{\eta^D_N\colon N\ge1\}$ from \eqref{E:3.19} obeys
\begin{equation}
\label{E:3.69}
\eta_N^D\,\,\underset{N\to\infty}\Lawarrow\,\,Z_\lambda^D(\textd x)\otimes\texte^{-\alpha\lambda h}\textd h.
\end{equation}
The measures $\{Z^D_\lambda\colon D\in\mathfrak D\}$ obey conditions (1-7) from Theorem~\ref{thm-2.3} and these identify their laws uniquely. In particular, on dyadic squares we have $Z^{S_n}_\lambda(\textd x)\laweq Y_\infty^{S_n}(\textd x)$.
\end{theorem}

\begin{proofsect}{Proof}
That subsequential limits of~$\eta^D_N$ take the form on the right of \eqref{E:3.69} has been shown in Lemma~\ref{lemma-3.7} and that the resulting $Z^D_\lambda$ measures obey properties (1-7) from Theorem~\ref{thm-2.3} is the content of Proposition~\ref{lemma-3.8}. Thanks to the representation in Proposition~\ref{prop-unique}, the law of  $Z^D_\lambda$ is determined for~$D$ being any dyadic square. In order to prove the theorem, it thus suffices to show that the law of $Z^D_\lambda$ is similarly determined for all~$D\in\mathfrak D$.

Let $D\in\mathfrak D$. Since we may assume that~$D$ belongs to the distinguished set~$\mathfrak D_0$ of domains, the properties (1-7) from Theorem~\ref{thm-2.3} apply to~$D$ as well. Thus, in particular, $Z_\lambda^D(\overline D\smallsetminus D^n)\to 0$ for any sequence of measurable sets~$D^n$ such that $D^n\uparrow D$. Thanks to property~(1) in Theorem~\ref{thm-2.3}, the same holds even if we take~$D^n$ to be the union of all open dyadic squares $S_{n,i}$, $i=1,\dots,m(n)$ of side length~$2^{-n}$ whose closure is contained in~$D$. However, in light of properties~(4-5) in Theorem~\ref{thm-2.3}, in this case we may write
\begin{equation}
1_{D^n}(x)Z^D_\lambda(\textd x) \laweq \texte^{\alpha\lambda\Phi^{D,D^n}(x)}\sum_{i=1}^{m(n)} Y^{S_{n,i}}_\infty(\textd x),
\end{equation}
where $\{Y^{S_{n,i}}_\infty\colon i=1,\dots,m(n)\}$ are independent of $\Phi^{D,D^n}$ and of one another and are equidistributed, modulo a shift, to $Y^{S_n}_\infty$. The law of $Z^D_\lambda$ is thus determined solely by those of $\{Y^{S_n}_\infty\colon n\ge1\}$ and the Gaussian fields  $\{\Phi^{D,D^n}\colon n\ge1\}$. We conclude that limit \eqref{E:3.69} exists for all~$D\in\mathfrak D$ and is the same regardless of the approximating sequence of lattice domains~$D_N$ and/or the way~$a_N$ approaches the limit \eqref{E:1.10}.
\end{proofsect}

\section{Beyond untruncated second moments}
\noindent
Our next goal is to eliminate the restriction to $\lambda<1/\sqrt2$ assumed throughout the proofs in the previous section. There were three specific steps where this restriction was crucially used: the non-triviality of the subsequential limits of $\{\eta_N^D\colon N\ge1\}$, the factorization property in Lemma~\ref{lemma-aux} and the estimate of expectations of integrals against~$Z_\lambda$-measures in Lemma~\ref{lemma-3.11} based on Lemma~\ref{lemma-3.9a}. This is because all three rely on the second moment estimate on the size of the level-set in Lemma~\ref{lemma-3.3} which fails when $\lambda\ge1/\sqrt2$.  (Lemmas~\ref{lemma-3.1} and \ref{lemma-3.2} hold for all~$\lambda\in(0,1)$.)

It turns out, and this is no surprise in this subject area, that the lack of the second moment is remedied by introducing a suitable truncation. This will help us fix the above three second-moment calculation while preserving the overall strategy of the proof.

\subsection{Truncated measures}
Let us start with a truncated version of the measures in \eqref{E:3.19}. Pick a sequence of domains $\{D_N\}$ approximating, via \twoeqref{E:2.1}{E:2.2}, a given continuum domain~$D\in\mathfrak D$. Recall our earlier notation $\Lambda_r(x):=\{z\in\Z^2\colon |z-x|\le r\}$ and, for each~$N\ge1$ and each~$x\in D_N$, let
\begin{equation}
n(x):=\max\bigl\{n\ge0\colon\Lambda_{\texte^{n+1}}(x)\subseteq D_N\bigr\}.
\end{equation}
Observe that $\log N - c \leq n(x) \leq \log N + c'$ for all $x \in D_N$ such that $\text{dist}(x, D_N^{\text c}) > \epsilon N$, with the first constant depending only on the choice of $\epsilon > 0$ and the second only on $D$. Define now the sequence of domains
\begin{equation}
\label{e:4.2}
\Delta^k(x) := \begin{cases}
	\emptyset			  & \text{for } k = 0\,,\\
	\Lambda_{\texte^k}(x) & \text{for } k= 1, \dots, n(x)-1 \,, \\
	D_N 				  & \text{for } k = n(x) \,.
\end{cases} 
\end{equation}
In accord with \eqref{E:A.6}, for $V\subseteq U$ let us write $\varphi^{U,V}$ for the conditional field $E(h^U|\sigma(h^U(z)\colon z\in U\smallsetminus V))$. We now set
\begin{equation}
\label{e:2.8b}
S_k(x) := \varphi^{D_N,\Delta^k(x)}(x),\qquad k=0,\dots,n(x).
\end{equation}
Observe that, by definition, $S_0(x) = h^{D_N}(x)$ while $S_{n(x)}(x) = 0$.

Next, for a given sequence $a_N$ such that \eqref{E:1.10} holds for some $\lambda\in(0,1)$ and $M>0$, define the truncation event
\begin{equation}
T_{N,M}(x):=\bigcap_{k= k_N}^{n(x)}\biggl\{\Bigl|S_k(x)- a_N\frac{n(x)-k}{n(x)}\Bigr|\le M(n(x)-k)^{3/4}\biggr\}
\end{equation}
 where
\begin{equation}
\label{E:4.13}
\OL k_N:=\frac18\log(K_N)=\frac14\bigl[(1-\lambda^2)+o(1)\bigr]\log(N)\,. \eOL
\end{equation}
Consider the point measure
\begin{equation}
\label{E:3.19a}
\wh\eta_N^{D,M}:=\frac1{K_N}\sum_{x\in D_N}\1_{T_{N,M}(x)}\,\delta_{x/N}\otimes\delta_{\,h^{D_N}(x)-a_N}
\end{equation}
Obviously, $\langle \widehat\eta_N^{D,M},f\rangle\le\langle\eta_N^D,f\rangle$ for any measurable $f\ge0$.
Let us now re-run the arguments from the previous section replacing the key second-moment lemmas with their truncated versions. Our first point to note is that the difference between the measures $\wh\eta_N^{D,M}$ and $\eta^D_N$ disappears when the truncation is removed  by taking~$M\to\infty$.  For this we introduce the truncated level set
\begin{equation}
\label{e:2.8a}
\wh\Gamma_N^{D,M}(b) :=
\Bigl\{ x \in D_N \colon h^{D_N}(x) \ge a_N+b,\,  T_{N,M}(x) \text{ occurs} \Bigr\} \,.
\end{equation}
Then we have:

\begin{lemma}
\label{lemma-4.1}
For each $\lambda\in(0,1)$ and each~$b_0>0$ there are constants~$c,\tilde c\in(0,\infty)$ such that for all~$D\in\mathfrak D$, all $b\in[-b_0,b_0]$, all $M\ge1$ and all~$N$ sufficiently large,
\begin{equation}
\bigl|\Gamma_N^D(b)\smallsetminus\wh\Gamma_N^{D,M}(b)\bigr|\le c\texte^{-\tilde cM^2}(\diam D)^{2+2\lambda^2}K_N.
\end{equation}
\end{lemma}

\noindent
Using this lemma we immediately get
that for any bounded, measurable $f\colon D\times\R\to\R$, 
\begin{equation}
\label{E:4.10}
\lim_{M\to\infty}\,\limsup_{N\to\infty}\,\bigl|\langle \wh\eta_N^{D,M},f\rangle-\langle\eta_N^D,f\rangle\bigr|=0.
\end{equation}
Since Lemma~\ref{lemma-3.1} and the aforementioned domination of $\wh\eta_N^{D,M}$ by~$\eta^D_N$ show that the family of measures $\{\wh\eta_N^{D,M}\colon N\ge1\}$ is tight in the topology of vague convergence, we can extract a subsequential weak limit $\wh\eta^{D,M}$ and study its properties. 

The first and foremost question is non-triviality of the limit. Here we need an analogue of Lemma~\ref{lemma-3.3}, now without restrictions on~$\lambda$. For $b,b'\in\R$ with $b<b'$, abbreviate
\begin{equation}
\wh\Gamma_N^{D,M}(b,b'):=\wh\Gamma_N^{D,M}(b)\smallsetminus\wh\Gamma_N^{D,M}(b')\,.
\end{equation}
Then we have:

\begin{lemma}
\label{lemma-4.2}
Let $\lambda\in(0,1)$. For all $\epsilon>0$, all $M\ge0$ and all~$b,b'\in\R$ with $b<b'$, there is $c=c(M,b,b',\epsilon)\in(1,\infty)$ such that for all $D\in\mathfrak D$ and all $N$ large enough,
\begin{equation}
E\bigl(|\wh\Gamma_N^{D,M}(b,b')\cap D_N^\epsilon|^2\bigr)\le c (\diam D)^{4+4\lambda^2}K_N^2
\end{equation}
\end{lemma}

\noindent
The second moment calculation spelled out in the proof of Corollary~\ref{cor-3.4} together with Lemma~\ref{lemma-4.1} then show that every subsequential weak limit $\wh\eta^{D,M}$ of measures $\{\wh\eta_N^{D,M}\colon N\ge1\}$ has positive total mass with positive probability, provided $M$ is chosen large enough. In light of the domination of~$\wh\eta^{D,M}_N$ by~$\eta^D_N$, the same applies to any subsequential limit of the measures $\{\eta^D_N\colon N\ge1\}$.

The hardest point to be addressed is factorization. This is the subject of the following lemma which effectively replaces Lemma~\ref{lemma-aux}:

\begin{lemma}
\label{lemma-4.3}
Let $\lambda\in(0,1)$. Using the notation $A_N:=\{x\in\Z^2\colon x/N\in A\}$, for each open~$A\subseteq D$ and each~$b\in\R$, we have
\begin{equation}
\label{E:4.12}
\lim_{N\to\infty}\,
\frac1{K_N}E\Bigl|\bigl|\wh\Gamma_N^{D,M}(0)\cap A_N\bigr|-\texte^{\alpha\lambda b}\bigl|\wh \Gamma_N^{D,M}(b)\cap A_N\bigr|\Bigr|=0.
\end{equation}
\end{lemma}

Lastly, we will need one more lemma dealing with the maximum of the field~$\Phi^{S_n,\wt S_{n,m}}$, where~$S_n$ is a dyadic square of side-length~$2^{-n}$ and $S_{n,m}$ is the disjoint union of $4^m$ dyadic squares of side-length~$2^{-(n+m)}$ that just barely fit into~$S_n$. Recall that $\wt S_{n+m}^\delta$ is the union of ``shrunk'' dyadic squares~$S_{n,m}^\delta$ centered at the same points as $S_{n,m}$, respectively. See the paragraph before \eqref{E:3.57}. 
\begin{lemma}
\label{lemma-4.5}
For each~$\delta>0$ there is a constant $c=c(\delta)$ such that
\begin{equation}
\label{E:4.23}
P\biggl(\,\,\sup_{x\in \wt S_{n,m}^\delta}\Phi^{S_n,\wt S_{n,m}}(x) >  2\sqrt g\log(2^m)+c\sqrt{\log(2^m)}\biggr)\,\underset{m\to\infty}\longrightarrow\,0.
\end{equation}
\end{lemma}

Deferring the proofs of these lemmas to the next subsection, we use them to prove:
\begin{theorem}
\label{thm-4.4}
The statement of Theorem~\ref{thm-3.13} applies to all $\lambda\in(0,1)$.
\end{theorem}
\begin{proofsect}{Proof}
Consider a countable family~$\mathfrak D_0$ of domains in~$\mathfrak D$ which include all dyadic squares. A diagonal argument permits us to to extract a subsequence along which $\wh\eta^{D,M}_N$ tends in law to a limit $\wh\eta^{D,M}$ for every~$D\in\mathfrak D_0$ and every integer~$M\ge1$. Applying monotonicity in~$M$, we can then define $\eta^D:=\lim_{M\to\infty}\wh\eta^{D,M}$. By \eqref{E:4.10}, $\eta^D$ is the limit of $\eta^D_N$ along the chosen subsequence.
Lemma~\ref{lemma-4.3} implies that~$\eta^D$ obeys \eqref{E:3.24} for every~$f$ of the stated form. By Lemma~\ref{lemma-3.7} we then have
\begin{equation}
\label{E:4.15}
\eta^D(\textd x\,\textd h)=Z^D_\lambda(\textd x)\otimes\texte^{-\alpha\lambda h}\textd h
\end{equation}
for some a.s.-finite random Borel measure~$Z^D_\lambda$ which has positive mass with positive probability. By the same reasoning as in the proof of Proposition~\ref{lemma-3.8}, the measures $\{Z^D_\lambda\colon D\in\mathfrak D_0\}$ obey properties~(1-7) in Theorem~\ref{thm-2.3}. In particular,~$Z^D_\lambda$ charges every non-empty open set a.s. 

In order to determine the law of $Z^D_\lambda$ uniquely, and thus prove the existence of the limit, we claim that, on the dyadic square~$S_n$, we have the representation
\begin{equation}
\label{E:4.14}
Z^{S_n}_\lambda(\textd x)\laweq Y_\infty^{S_n}(\textd x),
\end{equation}
where $Y_\infty^{S_n}$ is the measure constructed in Lemma~\ref{lemma-conv}. As the Gibbs-Markov property for~$Z^D_\lambda$ was already proved as part of the properties of Theorem~\ref{thm-2.3} above, the starting equation \eqref{E:3.53} is valid and the argument thereafter applies. We just need to replace Lemma~\ref{lemma-3.11} with a suitable analogue that does not rely on the existence of the second moment of $Z^{S_n}_\lambda(S_n)$.

The idea is to reintroduce the truncation while sticking with the $N\to\infty$ limit measures. Indeed, using the above subsequential limit point $\wh\eta^{D,M}$ of $\{\wh\eta^{D,M}_N\colon N\ge1\}$, we define
\begin{equation}
\wh Z^{D,M}_\lambda(A):=\alpha\lambda\,\wh\eta^{D,M}\bigl(A\times[0,\infty)\bigr)
\end{equation}
for each Borel measurable~$A$. These measures are not expected to obey most of the properties in Theorem~\ref{thm-2.3}. Notwithstanding, we have
\begin{equation}
\label{E:4.16}
\wh Z^{D,M}_\lambda(A)\le Z^D_\lambda(A)\quad\text{and}\quad \wh Z^{D,M}_\lambda(A)\uparrow Z^D_\lambda(A)\text{ as }M\to\infty.
\end{equation}
We will refer to $\wh Z^{D,M}_\lambda$ as a ``truncated measure'' although this not very accurate.

\begin{remark}
Note that Lemma~\ref{lemma-4.3} gives us a factorization property \eqref{E:4.15} for $\wh\eta^{D,M}$ as well; just replace $Z^D_\lambda$ by $\wh Z^{D,M}_\lambda$. One might thus be tempted to think that~$\wh\eta^{D,M}$ also satisfies the Gibbs-Markov property. However, this is false because the addition of the ``binding field''~$\Phi^{D,\wt D}$ changes the truncation events on the subdomain~$\wt D$. In any case, if the properties (1-7) of Theorem~\ref{thm-2.3} were true for $\wh Z^{D,M}_\lambda$, our argument from the previous section would represent this measure using a derivative martingale and, later, by the LQG measure. This would lead to a contradiction because the LQG measure is known to lack the second moment for all $\lambda\in[1/\sqrt2,1)$ yet (by Fatou) $\wh Z^{D,M}_\lambda$ is square integrable for all $\lambda\in(0,1)$.
\end{remark} 

Moving back to the proof of Theorem~\ref{thm-4.4}, we now define a measure~$\wt Z^{S_n,M}_{m}$ by \eqref{E:3.53} with the~$Z^{S_{n+m,i}}_\lambda$ on the right-hand side now replaced by their truncated analogues, 
\begin{equation}
\label{E:3.53a}
\wt Z^{S_n,M}_{\lambda}(\textd x):=
\sum_{i=1}^{4^m}\texte^{\alpha\lambda\Phi^{S_n,\wt S_{n,m}}(x)}\1_{S_{n+m,i}}(x)
\wh Z^{S_{n+m,i},M}_\lambda(\textd x).
\end{equation}
For each bounded, measurable $f\colon \overline D\mapsto[0,\infty)$ and each $\delta>0$ we then have 
\begin{equation}
\label{E:4.18}
E\bigl(\texte^{-\langle Z^{S_n}_\lambda,f\rangle}\bigr)
\le E\bigl(\texte^{-\langle \wt Z^{S_n,M}_{\lambda},f\rangle}\bigr)
\le E\bigl(\texte^{-\langle \wt Z^{S_n,M}_{\lambda},f_\delta\rangle}\bigr),
\end{equation}
where we $f_\delta$ is as defined just before \eqref{E:3.57}. Let $\Phi^{S_n,\wt S_{n,m}}$ be independent of the truncated measures $\{\wh Z^{S_{n+m,i},M}_{\lambda}\colon i=1,\dots,4^m\}$, which are themselves regarded as independent, and set
\begin{equation}
\label{E:3.59b}
\wt X_i:=\int_{S_{n+m,i}}\texte^{\alpha\lambda\Phi^{S_n,\wt S_{n,m}}(x)}f_\delta(x) \wt Z^{S_{n+m,i},M}_\lambda(\textd x).
\end{equation}
Noting that $\langle \wt Z^{S_{n},M}_{\lambda},f_\delta\rangle=\sum_{i=1}^{4^m}\wt X_i$, we then get for each $\epsilon>0$,
\begin{equation}
\label{E:4.20}
E\bigl(\texte^{-\langle \wt Z^{S_{n},M}_\lambda,f_\delta\rangle}\bigr)
\le E\biggl(\exp\Bigl\{-\texte^{-\epsilon}\sum_{i=1}^{4^m}E\bigl(\wt X_i\1_{\{\wt X_i\le\epsilon\}}\big|\Phi^{S_n,\wt S_{n,m}}\bigr)\Bigr\}\biggr)
\end{equation}
from the ``reverse Jensen'' inequality in Lemma~\ref{lemma-3.11a}.

To replace Lemma~\ref{lemma-3.11}, we claim that, for each $\epsilon>0$,
\begin{equation}
\label{E:4.21}
\sum_{i=1}^{4^m}E\bigl(\wt X_i\1_{\{\wt X_i>\epsilon\}}\big|\Phi^{S_n,\wt S_{n,m}}\bigr)
\,\,\underset{m\to\infty}\longrightarrow\,0
\end{equation}
in probability. For this let $A_{n,m}$ denote the event in \eqref{E:4.23} and note that, on~$A_{n,m}$, we can use straightforward calculations to bound
\begin{equation}
\sum_{i=1}^{4^m}E\bigl(\wt X_i\1_{\{\wt X_i>\epsilon\}}\big|\Phi^{S_n,\wt S_{n,m}}\bigr)
\le\frac1\epsilon\texte^{4\alpha\lambda\sqrt g\log(2^m)+c(\delta)\sqrt{\log(2^m)}}\,\sum_{i=1}^{4^m}E\bigl(\wh Z^{S_{n+m,i},M}_\lambda(S_{n+m,i}^\delta)^2\bigr).
\end{equation}
Here we will finally benefit from using the truncated measures $\wh Z^{S_{n+m,i},M}_\lambda$. Indeed,
Lemma~\ref{lemma-4.2}, monotonicity and Lemma~\ref{lemma-3.2} together with the scaling of the integral in~\eqref{E:3.6} ensure that, for some $c,c'\in(0,\infty)$ depending on $M$ and $n$,
\OL
\begin{equation}
\begin{aligned}
E\bigl(\wt Z^{S_{n+m,i},M}_\lambda(S_{n+m,i}^\delta)^2\bigr)
\le c \bigl[E(Z^{S_{n+m,i}}_\lambda(S_{n+m,i}))\bigr]^2
\le c'(2^m)^{-2(2+2\lambda^2)}
\end{aligned}
\end{equation}
\eOL
where the last inequality follows from \eqref{E:3.65}. Since $4\alpha\lambda\sqrt g = 8\lambda$, this yields
\begin{equation}
\sum_{i=1}^{4^m}E\bigl(\wt X_i\1_{\{\wt X_i>\epsilon\}}\big|\Phi^{S_n,\wt S_{n,m}}\bigr)
\le\frac{c'}\epsilon (2^m)^{-4(1-\lambda)^2}\texte^{\tilde c\sqrt m}
\qquad\text{on }A_{n,m}
\end{equation}
and so \eqref{E:4.21} follows from \eqref{E:4.23}.

Plugging \eqref{E:4.21} into \eqref{E:4.20} then gives
\begin{equation}
\label{E:4.27}
\limsup_{m\to\infty}E\bigl(\texte^{-\langle \wt Z^{S_n,M}_\lambda,f_\delta\rangle}\bigr)
\le\limsup_{m\to\infty} E\bigl(\texte^{-\texte^{-\epsilon}E(\langle \wt Z^{S_n,M}_\lambda,f_\delta\rangle
|\Phi^{S_n,\wt S_{n,m}})}\bigr)
\end{equation}
Our next task is to prove
\begin{equation}
\label{E:4.28}
E\bigl(\langle Z^{S_n}_\lambda,f_\delta\rangle
\,\big|\,\Phi^{S_n,\wt S_{n,m}}\bigr)
-E\bigl(\langle \wt Z^{S_n,M}_\lambda,f_\delta\rangle
\,\big|\,\Phi^{S_n,\wt S_{n,m}}\bigr)
\,\longrightarrow\,0
\end{equation}
in probability in the limit as~$m\to\infty$ followed by~$M\to\infty$. Since the left hand side above is non-negative, it thus suffices to show the convergence to $0$ in the mean.

Since $\Var(\Phi^{S_n,\wt S_{n,m}}(x))\le g\log(2^m)+c$ uniformly on $\wt S_{n,m}^\delta$, the expectation on the left-hand side of \eqref{E:4.28} is at most
\begin{equation}
\text{l.h.s.\ of\ \eqref{E:4.28}}
\le
c(2^m)^{2\lambda^2}4^m\Vert f\Vert_\infty
\Bigl[E\bigl(Z^{S_{n+m}}_\lambda(S_{n+m})\bigr)-
E\bigl(\wh Z^{S_{n+m},M}_\lambda(S_{n+m})\bigr)\Bigr].
\end{equation}
Lemma~\ref{lemma-4.1} now implies
\begin{equation}
E\bigl(Z^{S_{n+m}}_\lambda(S_{n+m})\bigr)-
E\bigl(\wh Z^{S_{n+m},M}_\lambda(S_{n+m})\bigr) 
\OL \le c\texte^{-\tilde c M^2}(2^{m})^{-2-2\lambda^2}.  \eOL
\end{equation}
This bounds the expectation of \eqref{E:4.28} by $c\texte^{-\tilde c M^2}\Vert f\Vert_\infty$ which tends to zero as~$M\to\infty$. 

Combining \eqref{E:4.18} with \twoeqref{E:4.27}{E:4.28}, using \eqref{E:3.54}, Lemma~\ref{lemma-conv} and the limits~$\epsilon\downarrow0$ and~$\delta\downarrow0$ (with \eqref{E:3.57}) thus show
\begin{equation}
E\bigl(\texte^{-\langle Z^{S_n}_\lambda,f\rangle}\bigr)
\le E\bigl(\texte^{-\langle Y^{S_n}_\infty,f\rangle}\bigr).
\end{equation}
Jointly with \eqref{E:3.55}, we then get \eqref{E:4.14}. The same argument as in the proof of Theorem~\ref{thm-3.13} now gives uniqueness of the law of~$Z_\lambda^D$ for all~$D\in\mathfrak D$.
\end{proofsect}

\subsection{Truncated moment calculations}
We now move to the technical statements (Lemmas~\ref{lemma-4.1}--\ref{lemma-4.3}) in the proof of Theorem~\ref{thm-4.4} whose proof was deferred from the previous subsection to here. For $\epsilon>0$, we write $D_N^\epsilon:=\{x\in D_N\colon \dist(x,D_N^\cc)>\epsilon N\}$. We will need some observations concerning the law of the random variables~$S_k$ defined in \eqref{e:2.8b}.

\begin{lemma}
\label{lem:var_estimates}
Recall that $g:=2/\pi$. For each $\epsilon > 0$ and each~$r>0$, there is $c=c(\epsilon,r)\in(0,\infty)$ such that for all~$D\in\mathfrak D$ with $\diam D\le r$ and all $N$ large enough, we have:
\settowidth{\leftmargini}{(11)}
\begin{enumerate}
\item For all $x \in D_N$ and all 
$k_N \leq k \leq m < n(x)$, 
\begin{equation}
\Var \bigl(S_k(x) - S_m(x)\bigr) = (m-k)g + o(1)
\end{equation}
where $o(1) \to 0$ when $N \to \infty$ uniformly in $k$. 
\item For all $x \in D_N^\epsilon$ and all~$k$ with $k_N \leq k \leq n(x)$,
\begin{equation}
\Var \bigl(S_k(x)\bigr) - (n(x)-k)g \, \in [0, c]
\end{equation}
\end{enumerate}
Moreover, for all $\ell \geq 1$ there is $c'=c'(\epsilon,\ell) >0$ such that for all $x \in D_N^\epsilon$, all~$k$ with $k_N \leq k \leq n(x)$, all~$m$ satisfying $k-\ell \leq m \leq k$ and all $y \in D_N$ such that $\Delta^{m+1}(y) \subseteq \Delta^k(x) \setminus \{x\}$, we have
\begin{equation}
E \bigl(S_k(x) S_m(y)\bigr) \leq (n(x)-k)g + c' 
\ \  \text{and} \ \ 
\Var \bigl(S_m(y)-S_k(x) \bigr) \in [g/2, c] \,.
\end{equation}
\end{lemma}

\begin{proofsect}{Proof}
Fix~$r>0$ and consider any domain~$D\in\mathfrak D$ with~$\diam D\le r$. By the Gibbs-Markov property, translation invariance and Green function asymptotics~\eqref{G-asymp2},
\begin{equation}
\begin{aligned}
\Var (S_k(x) - S_m(x)) &= \Var \Bigl( \varphi^{\Lambda_{\texte^m}(0), \Lambda_{\texte^k}(0)}(0) \Bigr)
\\
&= G_{\Lambda_{\texte^m}(0)}(0,0) - G_{\Lambda_{\texte^k}(0)}(0,0) = g(m-k) + o(1) \,.
\end{aligned}
\end{equation}
This gives the first statement. For the second, we assume that $k < n(x)$ since otherwise it is  trivially true. Then, since $n(x) \geq \log N - c$ for \OL $c=c(\epsilon) > 0$, we may find $\tilde c = \tilde c(\epsilon,r) > 0$ \eOL such that $D_N \subseteq \Lambda_{\texte^{n(x) + \tilde c}}(x)$. Monotonicity of the Green function with respect to inclusion and similar considerations as above now show
\begin{equation}
\Var \bigl( S_k(x) \bigr) \leq 
G_{\Lambda_{\texte^{n(x)+\tilde c}}(0)}(0,0) - G_{\Lambda_{\texte^k}(0)}(0,0) \leq (n(x)-k)g + g \tilde c + o(1) \,.
\end{equation}
On the other hand, by definition $D_N \supseteq \Lambda_{\texte^{n(x) + 1}}(x)$ 
 and consequently 
\begin{equation}
\Var \bigl( S_k(x) \bigr) \geq 
G_{\Lambda_{\texte^{n(x)+1}}(0)}(0,0) - G_{\Lambda_{\texte^k}(0)}(0,0) \geq (n(x)-k)g + g/2 + o(1) \,.
\end{equation}
This completes the second statement.

Turning to the third statement, observe that the expectation there can be written explicitly as $E \varphi^{D_N, \Delta^k(x)}(x)\varphi^{D_N, \Delta^k(x)}(y)$. By the Gibbs-Markov property, this expectation equals
\begin{equation}
E \bigl(h^{D_N}(x) h^{D_N}(y)\bigr) - E \bigl(h^{\Delta^k(x)}(x) h^{\Delta^k(x)}(y)\bigr) \leq (n(x)-k)g + c' \,,
\end{equation}
where we have used the Green function asymptotics again. The constant $c' > 0$ above depends on the distance of $y$ to the boundary of $\Delta^k(x)$ relative to its diameter. This in turn is governed by the choice of $\ell$. 

Finally, the upper bound on the variance follows from the above bounds together with
\begin{equation}
\Var \bigl(S_m(y)-S_k(x) \bigr) = \Var \bigl(S_m(y)\bigr) + \Var \bigl(S_k(x)\bigr) - 2 E \bigl( S_m(y) S_k(x) \bigr)
\end{equation}
and the relation between $m$ and $k$.
As for the lower bound,
\begin{multline}
\Var\, \bigl(S_m(y)-S_k(x) \bigr) =
\Var\, \bigl(\varphi^{D, \Delta^k(x)}(y) - \varphi^{D, \Delta^k(x)}(x)
+ \varphi^{\Delta^k(x), \Delta^m(y)}(y)\bigr) \\
\geq \Var \bigl(\varphi^{\Delta^k(x), \Delta^m(y)}(y)\bigr)
\geq \Var \bigl(\varphi^{\Delta^{m+1}(x), \Delta^m(y)}(y)\bigr)
\geq g + o(1)
\end{multline}
thus proving the third statement as well.
\end{proofsect}

The next lemma notes that the dependency structure of the process $(S_k(x))_{k,x}$ is tree-like.
\begin{lemma}
\label{lem:independence}
If $x,y \in D_N$ and $0 \leq k_1 < k_2 \leq n(x)$, $0 \leq m_1 < m_2 \leq n(y)$ are such that $\Delta^{k_2}(x) \subseteq \Delta^{m_1}(y)$, then the increments $S_{k_1}(x) - S_{k_2}(x)$ and $S_{m_1}(y) - S_{m_2}(y)$ are independent. In particular, for any $x \in D_N$, the process $\bigl(S_k(x)\bigr)_{k=0}^n(x)$ has independent increments.
\end{lemma}

\begin{proofsect}{Proof}
This is a direct consequence of the definition and the Gibbs-Markov property.
\end{proofsect}

Finally, we will also need the following simple fact:
\begin{lemma}
\label{lemma-cond}
Suppose $X\laweq\NN(0,\sigma_X^2)$ and~$Y\laweq\NN(0,\sigma_Y^2)$ are independent. Then~$(X|X+Y)$, i.e., $X$ conditional on $X+Y$, obeys
\begin{equation}
\label{E:4.35}
(X|X+Y)\laweq\NN\biggl(\frac{\sigma_X^2}{\sigma_X^2+\sigma_Y^2}(X+Y),\,\frac{\sigma_X^2\sigma_Y^2}{\sigma_X^2+\sigma_Y^2}\biggr).
\end{equation}
\end{lemma}

\begin{proofsect}{Proof}
A simple algebra shows
\begin{equation}
X = \frac{\sigma_X^2}{\sigma_X^2+\sigma_Y^2}(X+Y)+\frac{\sigma_Y^2X-\sigma_X^2 Y}{\sigma_X^2+\sigma_Y^2}\,.
\end{equation}
The second expression on the right is a Gaussian random variable that is independent of $X+Y$, has mean zero and variance as the random variable on the right of \eqref{E:4.35}.
\end{proofsect}

We are now ready to control the defect to the level set size caused by the truncation:

\begin{proofsect}{Proof of Lemma~\ref{lemma-4.1}}
Pick $b,b'\in\R$ with $b<b'$. Let~$\epsilon>0$ and note that
\begin{multline}
\label{E:4.34}
\qquad
E\bigl|\Gamma_N^D(b)\smallsetminus\wh\Gamma_N^{D,M}(b)\bigr|
\le E\bigl|\Gamma_N^D(b)\smallsetminus D_N^\epsilon\bigr|+
E\bigl|\Gamma_N^D(b')\bigr|
\\
+\sum_{x\in D_N^\epsilon}\,\sum_{k=k_N}^{n(x)}
P \Bigl( h^{D_N}(x) - a_N \in[b,b') \,,\,
	\bigl| S_k(x) - a_N \tfrac{n(x)-k}{n(x)} \bigr|> M (n(x)-k)^{3/4} \Bigr).
\qquad
\end{multline}
By Lemma~\ref{lemma-3.1}, the first two expectations are bounded desired estimate provided we take~$\epsilon$ small and~$b'$ sufficiently large (depending on~$\diam D$). We thus have to show the bound for the double sum regardless of $\epsilon>0$ and~$b'\in[b,\infty)$.

Fix $x\in D_N^\epsilon$ and let $k \in\{ 0, \dots, n(x)\}$. 
We will estimate the probability on the right-hand side of \eqref{E:4.34} by conditioning on the value of $h^{D_N}(x)$. For this we note that, by Lemma~\ref{lem:independence}, $h^{D_N}(x)$ is the sum of independent random variables $S_k(x)$ and $S_0(x)-S_k(x)$. 
Applying Lemma~\ref{lemma-cond} for $X:=S_k(x)$ and $Y:=S_0(x)-S_k(x)$, we thus get for $s \in [b, b']$
\begin{multline}
\label{E:4.37}
\qquad
\bigl(S_k(x)\,\big|\,h^{D_N}(x)=a_N+s\bigr)
\\
\laweq\NN\biggl(\frac{\Var(S_k(x))}{\Var(S_0(x))}(a_N+s),\,\frac{\Var(S_k(x))\Var(S_0(x)-S_k(x))}{\Var(S_0(x))}\biggr).
\qquad
\end{multline}
Invoking the variance estimates in Lemma~\ref{lem:var_estimates}, we obtain
\begin{equation}
\Bigl|\frac{\Var(S_k(x))}{\Var(S_0(x))}-\frac{n(x) - k}{n(x)}\Bigr|\le \frac{c_1}{n(x)}
\end{equation}
and
\begin{equation}
\frac{\Var(S_k(x))\Var(S_0(x)-S_k(x))}{\Var(S_0(x))}
\le c_2 (n(x)-k)
\end{equation}
where the constants~$c_1$ and~$c_2$ are independent of~$k$,~$x$ and~$n(x)$ as chosen above. Plugging these in \eqref{E:4.37} and using that $a_N$ is proportional to~$n(x)$, a standard Gaussian estimate yields
\begin{equation}
\label{E:4.40}
P \Bigl(\bigl| S_k(x) - a_N \tfrac{n(x)-k}{n(x)} \bigr|> M (n(x)-k)^{3/4}\,\Big|\,h^{D_N}(x) = a_N+s \Bigr)
\le c\texte^{-\tilde c M^2 (n(x)-k)^{1/2}},
\end{equation}
Thanks to the uniformity in $s$ of \eqref{E:4.40}, the last term in \eqref{E:4.34} is bounded by
\begin{equation}
\label{e:4.51}
c\sum_{x\in D_N^\epsilon}\sum_{k=1}^{\infty} \texte^{-\tilde c M^2 k^{1/2}} P\bigl(h^{D_N}(x)-a_N\in[b,b')\bigr)
\le c'\texte^{-\tilde c M^2}\,E\bigl|\Gamma_N^D(b)\bigr|.
\end{equation}
By Lemma~\ref{lemma-3.2}, this obeys the desired bound as soon as~$N$ is sufficiently large. 
\end{proofsect}

Next we move to the proof of the second moment estimate for truncated level sets:

\begin{proofsect}{Proof of Lemma~\ref{lemma-4.2}}
Pick $b,b'\in\R$ with $b<b'$ and fix $\epsilon>0$ and $M > 0$. Given $N \geq 1$ and $x,y\in D_N^\epsilon$, we will first estimate the probability that $x,y\in\wh\Gamma_N^D(b,b')$ \OL for $|x-y| > K_N^{1/4}$. \eOL
Denote 
\begin{equation}
k := \bigl(\lceil \log^+|x-y| \rceil + 1\bigr) \wedge n(x) \,,
\end{equation}
and let $\ell \geq 1$ be the minimal such that
\begin{equation}
\label{E:4.43}
\Delta^{k-\ell}(x) \cap \Delta^{k-\ell}(y) = \emptyset 
\ \ \text{and} \ \ 
\Delta^{k-\ell+1}(x) \cup \Delta^{k-\ell+1}(y) \subseteq \Delta^k(x) \,. 
\end{equation}
Observe that since $n(x) \leq \log N + c$ and $n(y) \geq \log N - c'$ for $c = c(D) > 0$ and $c'=c(\epsilon) > 0$ we must have $\ell \leq \tilde c$ with $\tilde c = \tilde c(\epsilon, D) > 0$.
Also note that
\begin{equation}
\label{E:4.44}
h^{D_N}(x) = S_k(x) + \big (S_{k-\ell}(x) - S_k(x)\bigr) + \bigl(S_0(x) - S_{k-\ell}(x)\bigr) 
\end{equation}
and
\begin{equation}
\label{E:4.45}
h^{D_N}(y) = S_k(x) + \big (S_{k-\ell}(y) - S_k(x)\bigr) + \bigl(S_0(y) - S_{k-\ell}(y)\bigr)  \,.
\end{equation}
By \eqref{E:4.43} and Lemma~\ref{lem:independence}, the three terms on the right of \eqref{E:4.44} are independent of each other, while for the terms on the right of \eqref{E:4.45} we get that the last one is independent of the first two as well as of of the last term on the right of \eqref{E:4.44}. For any $t \in [-M(n(x)-k)^{3/4},M(n(x)-k)^{3/4}]$, $s_1,s_2\in[b,b')$ and any $u_1,u_2\in[-n(x)^{3/4},n(x)^{3/4}]$, we then write
\begin{equation}
\label{e:2.23c}
\begin{aligned}
P \Bigl( h(x) &- a_N \in \textd s_1,\,
		 h(y) - a_N \in \textd s_2 \,\Big|\,   S_k(x) - a_N \tfrac{n(x)-k}{n(x)} = t,
\\ 
& \qquad \qquad \qquad \qquad \qquad \qquad
		 S_{k-\ell}(x) - S_k(x) = u_1,\,
		 S_{k-\ell}(y) - S_k(x) = u_2 \Bigr) \\
& = P \Bigl( S_0(x) - S_{k-\ell}(x) - a_N \tfrac{k}{n(x)} + t + u_1 \in \textd s_1 \Bigr)
\\
&\qquad \qquad \qquad \qquad \times 
P \Bigl(S_0(y) - S_{k-\ell}(y) - a_N \tfrac{k}{n(x)} + t + u_2 \in \textd s_2 \Bigr) \\
& \leq \frac ck
	\exp \biggl\{ -\frac{\bigl( \frac{a_N}n(x) k - t - u_1 + s_1 \bigr)^2 +
			\bigl( \frac{a_N}n(x) k - t - u_2 + s_2 \bigr)^2}
	{2 g k} \biggr\} \textd s_1 \textd s_2 \\
& \leq \frac ck  
	\exp \Bigl\{ -\frac{a_N^2}{gn^2} k + \frac{a_N}{gn} (2t + u_1 + u_2 - s_1 - s_2)\Bigr\}
	\textd s_1 \textd s_2 \,. 
\end{aligned}
\end{equation}
Here in the first inequality we used Lemma~\ref{lem:var_estimates} to replace variances of the random variables $S_0(x) - S_{k-\ell}(x)$ and $S_0(y) - S_{k-\ell}(y)$ by~$gk$. This causes only a change in the multiplicative constant because, by our assumptions on~$t$, $u_1$, $u_2$, $s_1$ and~$s_2$, the quantities in  the squares in the exponent are both at most order~$k$. In the second inequality we opened up the squares and retained, through a bound, only the quantities that depend on~$a_N$.

Our next step is to integrate the above conditional probability with respect to the conditional law of $S_{k-\ell}(y) - S_k(x)$ and~$S_{k-\ell}(x) - S_k(x)$ given~$S_k(x)$. For this we will need to examine the dependency of $S_{k-\ell}(y) - S_k(x)$ on $S_k(x)$.
By Lemma~\ref{lem:var_estimates}, there are $c > 0$ and $c'' > c' > 0$ such that 
$E \bigl((S_{k-\ell}(y) - S_k(x) )S_k(x) \bigr) \leq c$ and 
$\Var \bigl(S_{k-\ell}(y) - S_k(x) \bigr) \in [c', c'']$. Consequently for all $t$ with $|t| \leq M (n(x)-k)^{3/4}$,
\begin{equation}
\Bigl|E \bigl(S_{k-\ell}(y) - S_k(x) \, \big|\,  S_k(x) - a_N \tfrac{n(x)-k}{n(x)} = t \bigr)\Bigr|
\leq c t / (n(x)-k + 1) \, \leq c M
\end{equation}
and
\begin{equation}
\Var \bigl(S_{k-\ell}(y) - S_k(x) \,\big|\,  S_k(x) \bigr)
\leq \Var \bigl(S_{k-\ell}(y) - S_k(x) \bigr) \, \leq c'' \,, 
\end{equation}
with the conditional expectation vanishing and the conditional variance bounded similarly for $S_{k-\ell}(x) - S_k(x)$. Since $a_N/n(x)$ is bounded, the Cauchy-Schwarz inequality shows
\begin{equation}
E \Bigl( \texte^{\frac{a_N}{gn} (S_{k-\ell}(y) - S_k(x)) + \frac{a_N}{gn}(S_{k-\ell}(x) - S_k(x))}
	 \, \Big|\,  S_k(x)  - a_N \tfrac{n(x)-k}{n(x)} = t\Bigr) \leq \tilde c\,,
\end{equation}
uniformly in $t$ as above. In conjunction with \eqref{e:2.23c}, this yields
\begin{equation}
\begin{aligned}
\label{E:4.51}
\qquad
P \Bigl( h^{D_N}(x) - a_N & \in [b,b'),  \, 
		 h^{D_N}(y) - a_N \in [b,b') \, \Big|  \, S_k(x) - a_N \tfrac{n(x)-k}{n(x)} = t \Bigr)
		 \\
 &\leq P \Bigl(|S_{k-l}(x) - S_k(x)| \vee |S_{k-l}(y) - S_k(x)| > n^{3/4}
 \, \Big|  \, S_k(x) - a_N \tfrac{n(x)-k}{n(x)} = t \Bigr) \\
 & \qquad \qquad + 
 \frac ck \int_{[b,b')}\textd s_1\int_{[b,b')}\textd s_2\,
	\exp \Bigl\{ -\frac{a_N^2}{gn^2} k + \frac{a_N}{gn} (2t - s_1 - s_2) \Bigr\}\,
\\
&\leq \frac{c'}k\,\exp \Bigl\{ -\frac{a_N^2}{gn^2} k + 2\frac{a_N}{gn}t \Bigr\}
\end{aligned}
\end{equation}
uniformly in above~$t$ above, where $c'$ depends on~$b$,~$b'$ and~$M$ and we have used the fact that the right hand side is at least $\texte^{-cn}$ for some $c > 0$.

Now if $k=n(x)$ then $S_k(x)=t=0$ and therefore the right hand side above is also a bound on the unconditional probability. Otherwise, we integrate the left-hand side of \eqref{E:4.51} with respect to the distribution of the random variable $S_k(x) - a_N \tfrac{n(x)-k}{n(x)}$ to get (abbreviating $\theta_n(k):=(n(x)-k)^{3/4}$)
\begin{equation}
\label{E:4.52}
\begin{aligned}
P \Bigl( & x,y\in\wh\Gamma_N^{D,M}(b,b') \Bigr) 
\\
&\leq
P \Bigl( h^{D_N}(x) - a_N \in [b,b'),\, h^{D_N}(y) - a_N \in [b,b'), 
		\bigl|S_k(x) -a_N \tfrac{n(x)-k}{n(x)}\bigr| \leq M \theta_n(k) \Bigr) 
\\ 
& \leq \frac{c}{k (n(x)-k)^{1/2}} \,\texte^{-\frac{a_N^2}{gn^2}k}\,\int_{|t| \leq M \theta_n(k)}
	\exp\biggl\{2\frac{a_N}{gn}t-\frac{(a_N\frac{n(x)-k}n(x)+t)^2}{2 g (n(x)-k)}\biggr\}\,\textd t  \\
& \leq \frac{c}{k (n(x)-k)^{1/2}} \,\texte^{-\frac{a_N^2}{gn^2}k-\frac{a_N^2}{2gn^2}(n(x)-k)}\,\int_{|t| \leq M \theta_n(k)}\texte^{+\frac{a_N}{gn}t}\textd t
 \\
& \leq \frac{c'}{k (n(x)-k)^{1/2}}\texte^{-\frac{a_N^2}{2gn^2}k-\frac{a_N^2}{2gn^2}n(x)+\tilde c M\theta_n(k)} \\
&\le c'' \frac{K_N}{N^2} \frac{n(x)^{1/2}}{k(n(x)-k+1)^{1/2}}\,\texte^{-\frac{a_N^2}{2gn^2}k+\tilde c M\theta_n(k)}
\end{aligned}
\end{equation}
for some constants $c,c',c'',\tilde c\in(0,\infty)$. The latter bound applies also to the case $k=n(x)$.
 
The desired expectation is now obtained by summing over~$x,y\in D_N^\epsilon$. This yields
\begin{multline}
\label{E:4.53}
\qquad
E \bigl| \wh{\Gamma}^{D, M}_{N}(b,b') \cap D_n^\epsilon \bigr|^2  
= \sum_{x,y \in D_N^\epsilon} P \Bigl( x,y \in \wh{\Gamma}^{\epsilon, M}_N(b,b') \Bigr)
\\
\le \sum_{\begin{subarray}{c}
x,y \in D_N^\epsilon\\|x-y|\le K_N^{1/4}
\end{subarray}}
P \Bigl( x\in \wh{\Gamma}^{\epsilon, M}_N(b,b') \Bigr)
+\sum_{\begin{subarray}{c}
x,y \in D_N^\epsilon\\|x-y|>K_N^{1/4}
\end{subarray}}
P \Bigl( x,y \in \wh{\Gamma}^{\epsilon, M}_N(b,b') \Bigr).
\qquad
\end{multline}
The first term on the right-hand side is bounded by $K_N^{1/2}E|\wh\Gamma_N^{D,M}(b,b')|=O(K_N^{3/2})$. For the second term we partition the pairs $(x,y)$ further depending on which annulus $\Delta^k(x)\smallsetminus \Delta^{k-1}(x)$ the vertex~$y$ belongs to. As there are order $N^2\texte^{2k}$ such pairs for a given~$k$, the bound \eqref{E:4.52} gives
\begin{equation}
\label{E:4.54}
\sum_{\begin{subarray}{c}
x,y \in D_N^\epsilon\\|x-y|>K_N^{1/4}
\end{subarray}}
P \Bigl( x,y \in \wh{\Gamma}^{\epsilon, M}_N(b,b') \Bigr)
 \leq c' \frac{K_N}{N^2} \sum_{k=k_N}^n
\frac{n^{1/2}}{k(n-k+1)^{1/2}}\,N^2\texte^{2k-\frac{a_N^2}{2gn^2}k+\tilde c M(n(x)-k)^{3/4}},
\end{equation}
where we set~$n$ to be the maximum of $n(x)$ over all $x \in D_N^\epsilon$ (Here we note that a change in~$n$ by a additive constant changes \eqref{E:4.52} only by a multiplicative constant.)

Now $\frac{a_N^2}{2gn^2}$ is asymptotic to $2\lambda^2<2$ in the limit as $N\to\infty$ and so the exponent on the right of \eqref{E:4.54} grows linearly with~$k$. The sum is thus dominated by the $k=n$ term. Since $n = \log N + O(1)$, simple algebra shows that the expression on the right of \eqref{E:4.54} is $O(K_N^2)$. Since all bounds above were uniform in~$D\in\mathfrak D$ with a given diameter, say,~$\diam D\le 1$, we only need to show how to get the diameter dependence explicitly.

A key point is that the bounds were also independent of the approximating sequence of domains~$D_N$, nor of the centering sequence~$a_N$ as long as it obeyed \eqref{E:1.6} and not even much on the sequence~$k_N$ in the cutoff for the event~$T_{N,M}$ as long as $N$ is large enough.
 Fix~$D\in\mathfrak D$ with $r:=\diam D\le1$, let~$D_N$ be a sequence of approximating domains obeying \twoeqref{E:2.1}{E:2.2} and set $D':= r^{-1}D$. Fix~$j\in\{0,1,\dots,\lfloor r^{-1}\rfloor\}$ and set
\begin{equation}
D_N':=D_{\lfloor N/r\rfloor-j},\quad a_N':=a_{\lfloor N/r\rfloor-j}\quad\text{and}\quad k_N':=k_{\lfloor N/r\rfloor-j}.
\end{equation}
Then~$\{D_N'\}$ is a sequence of domains approximating, in the sense of \twoeqref{E:2.1}{E:2.2}, domain~$D'$.
Using~$a_N'$ as the centering sequence and $k_N'$ as the cutoff in the event $T_{N,M}$ in the definition of $\wh\Gamma_N^{D',M}(b)$, we then get
\begin{equation}
\wh\Gamma_{\lfloor N/r\rfloor-j}^{D,M}(b)=\wh\Gamma_N^{D',M}(b)
\end{equation}
A calculation now shows that the normalizing sequence~$K_N'$ defined using~$a_N'$ above obeys
\MB
\begin{equation}
K_N'= \bigl(r^{2+2\lambda^2}+o(1)\bigr)K_{\lfloor N/r\rfloor-j}
\end{equation}
\eMB
Since every integer can be cast to the form $\lfloor N/r\rfloor-j$ for some~$N$ and some~$j$ as above, the claim for~$D$ follows from the claim for~$D'$. 
\end{proofsect}

Our final task is the proof of the factorization property for the truncated level sets:

\begin{proofsect}{Proof of Lemma~\ref{lemma-4.3}}
Using the Cauchy-Schwarz inequality, we bound the expression inside the limit in~\eqref{E:4.12} by the square root of
\begin{equation}
\label{E:4.14A}
\begin{aligned}
\frac1{K_N^2} \sum_{x,y \in A_N} 
E\Bigl( \bigl( 1_{\{h^{D_N}(x) \ge a_N\}} & - \texte^{\alpha \lambda b} 
	1_{\{h^{D_N}(x) \ge a_N + b\}}  \bigr)1_{T_{N,M}(x)}\\
	& \times\bigl( 1_{\{h^{D_N}(y) \ge a_N\}} - \texte^{\alpha \lambda b} 
		1_{\{h^{D_N}(y) \ge a_N + b\}} \bigr)1_{T_{N,M}(y)}\Bigr) \,. 
\end{aligned}
\end{equation}
To bound the sum above, we first consider pairs $x,y$ for which $m = \big \lfloor \log \|x-y\|_\infty \rfloor$ obeys
\begin{equation}
\label{E:4.15A}
m \ge \frac32 k_N
\quad \text{and}\quad 
\|x-y\|_\infty \in \bigl[ \texte^m + 2\texte^{k_N},   \texte^{m+1} - 2\texte^{k_N} \bigr].
\end{equation}
For such $x,y$, we have
\begin{equation}
\Delta^k(x) \cap \Delta^{k_N}(y) = \emptyset,\qquad k=k_N,\dots,n(x),
\end{equation}
and similarly
\begin{equation}
\Delta^k(y) \cap \Delta^{k_N}(x) = \emptyset,\qquad k=k_N,\dots,n(y).
\end{equation}
Consequently, letting $\FF := \sigma(h^{D_N(x)}_z \colon z \in D_N \smallsetminus (
\Delta^{k_N}(x) \cup \Delta^{k_N}(y)))$ the term corresponding to such~$x$ and~$y$ in \eqref{E:4.14A} can be written as
\begin{multline}
\label{E:4.17A}
\qquad E \Big[ \Big(
	P \bigl(h^{D_N}(x) \geq a_N \,\big|\, \FF \bigr) 
- \texte^{\alpha \lambda b} P \bigl(h^{D_N}(x) \geq a_N + b \,\big|\, \FF\bigr) \Big) 1_{T_{N,M}(x)} \\
\Big(
	P \bigl(h^{D_N}(y) \geq a_N \,\big|\, \FF \bigr) 
- \texte^{\alpha \lambda b} P \bigl(h^{D_N}(y) \geq a_N + b \,\big|\, \FF\bigr) \Big) 
 1_{T_{N,M}(y)} \Big] \,.
 \qquad
\end{multline}
We now write $h^D_N(x)$ as $h^D_N(x) = (S_0(x) - S_{k_N}(x)) + S_{k_N}(x)$ and note that the quantity in the parenthesis is independent of~$S_{k_N}(x)$, due to the fact that this term is $\FF$-measurable. Using a similar decomposition for $h^{D_N}(y)$, the expectation in \eqref{E:4.17A} is bounded by
\begin{equation}
\label{E:4.18A}
E \biggl[ 1_{T_{N,M}(x)} F \Bigl( S_{k_N}(x) - a_N \tfrac{\log N - k_N}{\log N} \Bigr)
	1_{T_{N,M}(y)} F \Bigr( S_{k_N}(y) - a_N \tfrac{\log N - k_N}{\log N} \Bigr)\biggr] \,,
\end{equation}
where
\begin{equation}
F(u) = \biggl|P \Bigl( \tilde h(0) \geq a_N \tfrac{k_N}{\log N} - u \Bigr)
- \texte^{\alpha \lambda b} P \Bigl( \tilde h(0) \geq a_N \tfrac{k_N}{\log N} - u + b \Bigr)\biggr| \,.
\end{equation}
with $\tilde h$ denoting the DGFF on $\Lambda_{\texte^{k_N}}$.

Since $a_N k_N/\log N \sim 2 \sqrt{g} \lambda k_N$ and $E \tilde h(0)^2 = g k_N + O(1)$, if we assume that $|u| < (k_N)^{7/8}$ then we can use tail asymptotics for the Gaussian density to estimate the quantities above. In particular, $P \big(\tilde h(0) \geq a_N \tfrac{k_N}{\log N} - u + b\big) / P \big(\tilde h(0) \geq a_N \tfrac{k_N}{\log N} - u \big)$ is asymptotic to
\begin{equation}
\frac{a_N\tfrac{k_N}{\log N} - u}{a_N \tfrac{k_N}{\log N} - u + b}
\exp \biggl\{\frac{-2b \bigl(a_N \tfrac{k_N}{\log N} \OL - \eOL u\bigr) - b^2}{2 E [\tilde h(0)^2]} \biggr\}
\sim \texte^{-\alpha \lambda b},\qquad N\to\infty.
\end{equation}
It follows that
\begin{equation}
\label{E:4.63}
F(u) = o(1) P \bigl(\tilde h(0) \geq a_N \tfrac{k_N}{\log N} - u \bigr)
\end{equation}
with $o(1)\to0$ as $N \to \infty$ uniformly in~$u$ satisfying $|u| < (k_N)^{7/8}$.
Since
\begin{equation}
\Bigl|S_{k_N}(x) - a_N \tfrac{\log N - k_N}{\log N}\Bigr| < k_N^{7/8},\qquad \text{on }T_{N,M}(x)
\end{equation}
whenever $N$ is large enough, with a similar condition holding for $y$ under $T_{N,M}(y)$, we can apply \eqref{E:4.63} in~\eqref{E:4.18A}. Reversing the step \eqref{E:4.18A}, the expectation in~\eqref{E:4.17A} is thus bounded by
\begin{equation}
o(1) P\Bigl( h^{D_N}(x) \ge a_N \,,\,\, h^{D_N}(y) \ge a_N\,,\,\, T_{N,M}(x)\,,\,\, T_{N,M}(y) \Bigr) \,.
\end{equation}
Proceeding as in~\eqref{E:4.54}, the contribution to the sum in~\eqref{E:4.14A} from terms where~\eqref{E:4.15A} holds is therefore at most $o((K_N)^2)$.

Turning to the remaining terms in the sum in~\eqref{E:4.14A}, if $x,y$ satisfy \OL $m < \frac32 k_N$ \eOL, then we bound the corresponding term by $4 \texte^{2\alpha \lambda (b\vee0)} P(h^{D_N}(x) \geq a_N +b\wedge0)$.
As in~\eqref{E:3.17} the contribution to the sum from all such terms is $o(K_N^2)$. For the pairs 
with $m \geq \frac32 k_N$ not satisfying the second restriction in~\eqref{E:4.15A}, we can bound by
$4 \texte^{2\alpha \lambda (b\vee0)} P(h^{D_N}(x) \geq a_N - b\wedge0 \,,\,\, h^{D_N}(y) \geq a_N - b\wedge0)$.  Observing that the number of such pairs for a given $m$ is at most order $N^2 \texte^{m + k_N} = o(1)N^2 \texte^{2m}$, the calculation in~\eqref{E:4.54} again shows that such terms contribute $o((K_N)^2)$ as well. The claim follows.
\end{proofsect}


\section{Proofs of main theorems}
\noindent
The goal of this short section is to give formal proofs of our main theorems. Before we do that, we still have to address one issue that has so been excluded from the discussion so far: the third component of the point process that captures the local behavior of the field near a point of an intermediate level set. 

\subsection{Local structure of intermediate level sets}
Henceforth, let $\eta_N^D$ denote the full three-component process defined in \eqref{E:2.6} and let~$\wh\eta_N^{D,M}$ be its truncation to points~$x$ where $T_{N,M}(x)$ holds.  In addition, define 
\begin{equation}
\wt\eta_N^{D,M}:=\frac1{K_N}\sum_{x\in D_N}\delta_{x/N}\otimes\delta_{\,h^{D_N}(x)-a_N}\otimes\delta_{\{h^{D_N}(x)-h^{D_N}(x+z)\colon z\in\Z^2\}}\1_{\wt T_{N,M}(x)}
\end{equation}
where
\begin{equation}
\wt T_{N,M}(x):=T_{N,M}(x)\cap\Bigl\{\max_{y\in\partial \Delta^{k_N}(x)}\bigl|h^{D_N}(x)-h^{D_N}(y)\bigr|\le k_N^2\Bigr\}.
\end{equation}
and $\Delta^k(x)$, $k_N$ are defined in~\eqref{e:4.2} and~\eqref{E:4.13} respectively.
Obviously, $\eta^D_N$ dominates $\wh\eta^{D,M}_N$ which dominates $\wt\eta_N^{D,M}$. First we note that the truncations do not really matter as soon as proper limits are taken:

\begin{lemma}
\label{lemma-5.1}
For any continuous, compactly-supported function $f\colon\overline D\times\R\times\R^{\Z^2}\to\R$,
\begin{equation}
\label{E:4.10a}
\lim_{M\to\infty}\,\limsup_{N\to\infty}\,\bigl|\langle \wt\eta_N^{D,M},f\rangle-\langle\eta_N^D,f\rangle\bigr|=0.
\end{equation}
\end{lemma}

\begin{proofsect}{Proof}
Thanks to Lemma~\ref{lemma-4.1}, we have \eqref{E:4.10}  for the three-component process as well.  In light of Lemma~\ref{lemma-3.1} it suffices to show that for any $\epsilon > 0$, the probability that $T_{N,M}(x)\smallsetminus\wt T_{N,M}(x)$ occurs at some~$x \in D_N^\epsilon$ goes to $0$ as $N \to \infty$. By the Union Bound, this probability is at most
\begin{equation}
\label{E:5.4}
\sum_{x\in D_N^\epsilon}\sum_{y\in \partial \Delta^{k_N}(x)}P\bigl(h^{D_N}(x)-h^{D_N}(y)>k_N^2\bigr)
\end{equation}
Since $h^{D_N}(x)-h^{D_N}(y)$ has mean zero and variance bounded by a constant times~$k_N$ uniformly for all such pairs whenever $N$ is large enough, the probability on the right is at most $\texte^{-ck_N^3}$. As $k_N\sim c'\log N$ and as the number of terms in the sum is only order $N^2 {\text e}^{k_N}$, the claim follows.
\end{proofsect}

 The principal computation to be done in this section is now the content of:

\begin{proposition}
\label{lemma-5.2}
Let~$\nu$ be the measure in \eqref{E:2.8} and, given any continuous, compactly-supported function $f\colon\overline D\times\R\times\R^{\Z^2}\to\R$, let
\begin{equation}
\label{E:5.5}
f_\nu(x,h):=E_\nu\bigl(f(x,h,\phi)\bigr)
\end{equation}
with the expectation over~$\phi$. Then for any $M>0$,
\begin{equation}
\lim_{N\to\infty}\,E\Bigl|\langle\wt\eta_N^{D,M},f\rangle-\langle\wt\eta_N^{D,M},f_\nu\rangle\Bigr|=0.
\end{equation}
\end{proposition}

For the proof we will need:

\begin{lemma}
\label{lemma-help}
Let~$\epsilon>0$ and, for~$x\in D_N^\epsilon$ and a sample of~$h^{D_N}$, let $\varphi_N$ denote the discrete-harmonic extension of the values of~$h^{D_N}$ on~$\{x\}\cup \Delta^{k_N}(x)^\cc$. Recall the notation $\fraka$ for the potential kernel associated with the simple symmetric random walk on~$\Z^2$. Then
\begin{equation}
\label{E:5.12}
\max_{x\in D_N^\epsilon}\,\max_{z\in\Lambda_r(x)}\,\sup_{\begin{subarray}{c}
h^{D_N}\in \wt T_{N,M}(x)\\|h^{D_N}(x)-a_N|\le\log\log N
\end{subarray}}
\,
\Bigl|
h^{D_N}(x)-\varphi_N(z) - \frac2{\sqrt g}\lambda\fraka(z-x)\Bigr|\,\underset{N\to\infty}\longrightarrow\,0
\end{equation}
\end{lemma}

\begin{proofsect}{Proof of Lemma~\ref{lemma-help}}
To show this, let $H_N(z,y)$ denote the probability that the simple random walk started at~$z$ first returns to $\{x\}\cup \Delta^{k_N}(x)^\cc$ at~$y$. (Note that $H_N(x,x)>0$ in this case.) Then
\begin{equation}
\label{E:5.8a}
h^{D_N}(x)-\varphi_N(z)=\sum_{y\in\partial \Delta^{k_N}(x)}H_N(z,y)\bigl[h^{D_N}(x)-h^{D_N}(y)\bigr].
\end{equation}
Recall the notation for~$S_k(x)$ and note that
\begin{equation}
S_{k_N}(x)\sum_{y\in\partial \Delta^{k_N}(x)}H_N(x,y) = \sum_{y\in\partial \Delta^{k_N}(x)}H_N(x,y)h^{D_N}(y).
\end{equation}
Swapping $h^{D_N}(x)$ for~$S_{k_N}(x)$ on the right hand side of \eqref{E:5.8a} then gives
\begin{multline}
\label{E:5.10a}
\qquad
h^{D_N}(x)-\varphi_N(z) = 
\bigl(1-H_N(z,x)\bigr)\bigl[h^{D_N}(x)-S_{k_N}(x)\bigr]
\\+\sum_{y\in\partial \Delta^{k_N}(x)}\bigl(H_N(z,y)-\OL H_N(x,y)\eOL\bigl)\bigl[S_{k_N}(x)-h^{D_N}(y)\bigr].
\qquad
\end{multline}
We claim that the second term on the right vanishes in the stated limits. Indeed, on $\wt T_{N,M}(x)$ we have $|h^{D_N}(x)-h^{D_N}(y)|\le k_N^2$ for each $y\in\partial \Delta^{k_N}(x)$ and so $|h^{D_N}(y)-S_{k_N}(x)|\le 2k_N^2$. The standard bounds on the regularity of the harmonic measure show $|H_N(z,y)-\OL H_N(x,y) \eOL|\le c r\texte^{-2k_N}$ for all $z\in\Lambda_r(x)$. The second term is thus of order $k_N^2\texte^{-k_N}$. 

Concerning the first term on the right of \eqref{E:5.10a} we note that, on the event $T_{N,M}(x)\cap\{|h^{D_N}(x)-a_N|\le \log\log N\}$ we have
\begin{equation}
h^{D_N}(x)-S_{k_N}(x)=a_N\frac{k_N}n(x)+O(k_N^{3/4})
\end{equation}
where $n=\log N+O(1)$, while
\begin{equation}
gk_N(1-H_N(z,x))=\fraka(z-x)+o(1)
\end{equation}
 uniformly in~$z\in\Lambda_r(x)$. Using the asymptotic \eqref{E:1.10} for~$a_N$, we then get \eqref{E:5.12}.
\end{proofsect}

\begin{proofsect}{Proof of Proposition~\ref{lemma-5.2}}
By way of limit arguments, we may assume that $f$ depends only on a finite number of coordinates of~$\phi$, say, those in~$\Lambda_r(0)$, and that~$f(x,h,\phi)\ne0$ implies $\dist(x,\Delta^\cc)>\epsilon$ and $h\in[b,b')$ for some~$\epsilon>0$ and~$b<b'$. Abbreviate
\begin{equation}
F(x,h,\phi):=f(x,h,\phi)-f_\nu(x,h).
\end{equation}
We then bound the expression by the square root of
\begin{multline}
\label{E:5.8}
\qquad
\frac1{K_N^2}\sum_{x,y\in D_N^\epsilon}E\biggl(F\Bigr(\tfrac xN,h^{D_N}(x)-a_N, 
\big(h^{D_N}(x)-h^{D_N}(x+z)\big)_{z \in \Lambda_r(0)}\Bigr)\1_{\wt T_{N,M}(x)}
\\
\times
F\Bigr(\tfrac yN,h^{D_N}(y)-a_N, \big(h^{D_N}(y)-h^{D_N}(y+z)\big)_{z \in \Lambda_r(0)}\Bigr)\1_{\wt T_{N,M}(y)}\biggr).
\qquad
\end{multline}
The argument at the end of the proof of Lemma~\ref{lemma-4.3} permits us to assume that $x,y$ are such that \eqref{E:4.15A} applies. We also assume that~$N$ is so large that $\texte^{k_N}>r$. Conditioning on the sigma algebra
\begin{equation}
\FF:=\sigma\bigl(h^{D_N}(z)\colon z\not\in \Delta^{k_N}(x)\cup \Delta^{k_N}(y)\bigr)
\end{equation}
then splits the expectation into a product of two parts, one for~$x$ and the other for~$y$. 
 Using the Gibbs-Markov decomposition to write~$h^{D_N}$ on~$\Delta^{k_N}(x)\smallsetminus\{x\}$ as $\varphi_N+\tilde h_N$, where $\tilde h_N$ is the DGFF on $\Delta^{k_N}(x)\smallsetminus\{x\}$ and $\varphi_N$ is as in Lemma~\ref{lemma-help},  we now write the term corresponding to~$x$ as
\begin{equation}
\label{E:5.10}
E \Bigl(\wt F_N\Bigr(\tfrac xN,h^{D_N}(x)-a_N,h^{D_N}(x)-\varphi_N(\cdot-x)\Bigr)\1_{\wt T_{N,M}(x)}\,\Big|\,\FF\Bigr),
\end{equation}
where
\begin{equation}
\wt F_N(x,h,\phi):=E F(x,h,\phi+\tilde h_N)
\end{equation}
 with the expectation with respect to $\tilde h_N$. 
Our aim is to show that the random variable under expectation in \eqref{E:5.10} is small uniformly in~$x\in D_N^\epsilon$ and the part of the configuration measurable with respect to~$\FF$. 

Thanks to uniform continuity of~$f$, the identity \eqref{E:5.12} permits us to replace $h^{D_N}-\varphi_N(\cdot-x)$ in \eqref{E:5.10} by $\frac2{\sqrt g}\lambda\fraka$ with a cost that tends deterministically to zero. The random variable under expectation then depends on the conditional field only through~$h^{D_N}(x)$. We now observe that, by the weak convergence of the DGFF on $\Lambda_{k_N}(0)\smallsetminus\{0\}$ to the DGFF on~$\Z^2\smallsetminus\{0\}$ --- which can be verified, e.g., by comparing covariances --- we get
\begin{equation}
E \wt F_N\bigl(x,h,\tfrac2{\sqrt g}\lambda\fraka\bigr)\,\underset{N\to\infty}\longrightarrow\,0
\end{equation}
uniformly in~$x$ and~$h$. We conclude that the conditional expectation in \eqref{E:5.10} is bounded by
\begin{equation}
o(1)1_{\{h^{D_N}(x)-a_n\in[b,b')\}}\1_{T_{N,M}(x)}
\end{equation}
 with~$o(1)\to0$ uniformly in~$x\in D_N^\epsilon$. In light of Lemma~\ref{lemma-4.3}, the quantity in \eqref{E:5.8} tends to zero as~$N\to\infty$. The claim follows.
\end{proofsect}

\subsection{Proofs of main results}
We are now ready to give the proofs of our main theorems:

\begin{proofsect}{Proof of Theorem~\ref{thm-2.1}}
Let $f\colon\overline D\times\R\times\R^{\Z^2}\to\R$ be continuous with compact support and let~$f_\nu$ be as in \eqref{E:5.5}. Theorems~\ref{thm-3.13} and~\ref{thm-4.4} and ensure that $\langle \eta_N^D,f_\nu\rangle$ tends in distribution to
\begin{equation}
\int Z^D_\lambda(\textd x)\otimes\texte^{-\alpha\lambda h}\textd h \,\,f_\nu(x,h)
=\int Z^D_\lambda(\textd x)\otimes\texte^{-\alpha\lambda h}\textd h\otimes\nu(\textd\phi)\,\,f(x,h,\phi).
\end{equation}
Lemma~\ref{lemma-5.1} and Proposition~\ref{lemma-5.2} then identify this with the distributional limit of $\langle \eta_N^D,f_\nu\rangle$. As this holds for all such~$f$, the claim follows.
\end{proofsect}

\begin{proofsect}{Proof of Corollary~\ref{cor2.2}}
This follows from Theorem~\ref{thm-2.1} with Lemma~\ref{lemma-3.1} used to reduce the problem to level sets between two values of the form~$a_N+b$.
\end{proofsect}

\begin{proofsect}{Proof of Theorem~\ref{thm-2.3}}
This was proved as part of the proofs of Theorems~\ref{thm-3.13} and~\ref{thm-4.4}.
\end{proofsect}

\begin{proofsect}{Proof of Theorem~\ref{thm-2.4}}
This is proved exactly as \cite[Theorem~7.2]{BL2}; one just needs to change the exponent 4 into $2+2\lambda^2$ in suitable places. We only verify the  parts  of this theorem where this change shows up. First off,  the independence of the limit \eqref{E:2.7} on the particular sequence~$a_N$ permits us to assume $a_N:=2\sqrt g\lambda\log N$ for which we then have $K_{rN}/K_N\to r^{2+2\lambda^2}$ as~$N\to\infty$. This  implies
the scaling relation \eqref{E:2.22}. The representation of $Z^S_\lambda\laweq Y_\infty^S$ for any square~$S$ then yields rotation invariance. This is because, under a conformal map~$f\colon D\to f(D)$, we have
\begin{equation}
C^{f(D),f(\wt D)}\bigl(f(x),f(y)\bigr) = C^{D,\wt D}(x,y)
\end{equation}
for any admissible $\wt D\subset D$ and thus
\begin{equation}
\Phi^{f(D),f(\wt D)}\circ f\laweq \Phi^{D,\wt D}.
\end{equation}
The rotation invariance of~$Y_\infty^S$ then follows from the rotation invariance of the function $\psi_\lambda^D$.

With these properties verified, the proof of \cite[Proposition~7.2]{BL2} can then be followed literally to yield,  for any~$u\colon D\to[0,\infty)$ bounded and measurable,
\begin{equation}
\label{E:5.21}
E\bigl(\texte^{-\langle Z^{f(D)}_\lambda,u\circ f\rangle}\bigr)\ge 
E\bigl(\texte^{-\langle Z^{D}_\lambda,|f'\circ f^{-1}|^{2+2\lambda^2}u\rangle}\bigr)
\end{equation}
Iterating this with~$f$ replaced by~$f^{-1}$ then gives equality in \eqref{E:5.21}. The claim follows.
\end{proofsect}

\begin{proofsect}{Proof of Theorem~\ref{thm-LQG}}
Let~$D\in\mathfrak D$ be such that it~$D$ fits an open dyadic square of side~$r$. For any integer~$k\ge0$, let $S_{k,i}$, $i=1,\dots,n(k)$ be open dyadic squares of side~$r2^{-k}$ that lie entirely in~$S$. Clearly, each $S_{k,i}$ has a non-empty intersection and contains exactly 4 squares of the form $S_{k+1,j}$ although there may be squares of the latter form that do not belong to any square of the form~$S_{k,i}$. For each $k\ge1$, let $\cmss H_k$ be the set of functions in $\cmss H_0^1(S)$ that are harmonic on each~$S_{k,i}$, $i=1,\dots,n(k)$, and vanish on~$D\smallsetminus\bigcup_{i=1}^{n(k)}S_{k,i}$. Then,  as is checked by the Gauss-Green formula,  $\{\cmss H_k\colon k\ge0\}$ are orthogonal subspaces of $\cmss H_0^1(S)$ with
\begin{equation}
\cmss H_0^1(S) = \bigoplus_{k\ge0}\cmss H_k.
\end{equation}
A minor complication that arises in this setting is that each~$\cmss H_k$ is infinitely dimensional. Still, by separability of $\cmss H_0^1(D)$, we can find an countable orthonormal basis $\{\wt\varphi_{k,j}\colon j\ge1\}$ in each~$\cmss H_k$. 

Let $\{X_{k,j}\colon k,j\ge1\}$ be i.i.d.\ standard normals and write $S^0:=D$ and $D^k:=\bigcup_{i=1}^{n(k)}S_{k,i}$. A covariance calculation shows
\begin{equation}
\Phi^{D^{k-1},D^k}\laweq \sum_{j\ge1}X_{k,j}\wt\varphi_{k,j}\quad \text{on }D^k,\qquad k\ge1.
\end{equation}
We also have
\begin{equation}
\Phi^{D,D^m}\laweq\sum_{k=1}^m\Phi^{D^{k-1},D^k}
\end{equation}
with the fields on the right-hand side regarded as independent.  Using a suitable coupling to realize these distributional identities as  almost sure equalities,  letting
\begin{equation}
\label{E:5.25}
Y_k^D(\textd x):=c\psi_\lambda^D(x)\sum_{i=1}^{n(k)}\texte^{\alpha\lambda\Phi^{D,D^k}(x)-\frac12\alpha^2\lambda^2E[\Phi^{D,D^k}(x)^2]}\1_{S_{k,i}}(x)\textd x
\end{equation}
and setting $\FF_k:=\sigma\bigl(X_{\ell,j}\colon \ell+j\le k\}$, then for any measurable~$A\subset D$,
\begin{equation}
\label{E:5.27}
E\bigl(Y_m^D(A)\,\big|\,\FF_k\bigr)=\int_A c\psi_\lambda^D(x)\mu_{k(k-1)/2}^{D,\alpha\lambda}(\textd x),\qquad m\ge k,
\end{equation}
where $\mu_n^{D,\beta}$ is the measure defined in \eqref{E:2.17a} for the basis $\{\varphi_n\colon n\ge1\}$ in~$\cmss H_0^1(D)$ which is obtained by reordering $\{\wt\varphi_{k,j}\colon j\ge1\}$ according to the complete order
\begin{equation}
(k,j)\preceq(k',j')\quad\Leftrightarrow\quad k+j< k'+j'\quad\text{or}\quad k+j=k'+j'\,\,\,\&\,\,\, j\le j'. 
\end{equation}
Since $\lambda<1$ and thus $\alpha\lambda<\beta_\cc$, it is known  (cf a remark after Rhodes and Vargas~\cite[Theorem~5.5]{RV-review})  that $\mu_n^{D,\beta}$ converges to a non-trivial measure $\mu_\infty^{D,\beta}$ almost surely and in~$L^1$. It follows that $\{Y_k^D(A),\FF_k\colon k\ge1\}$ is a uniformly integrable martingale. The Martingale Convergence Theorem then gives
\begin{equation}
Y_k^D(A)\underset{k\to\infty}\longrightarrow\, Y_\infty^D(A)\quad\text{ a.s.\ and in }L^1.
\end{equation}
Using this in \eqref{E:5.27} shows
\begin{equation}
E\bigl(Y_\infty^D(A)\,\big|\,\FF_k\bigr)=\int_A c\psi_\lambda^D(x)\mu_{k(k-1)/2}^{D,\alpha\lambda}(\textd x),\qquad k\ge1.
\end{equation}
The Levy Backward Theorem and the convergence $\mu_n^{D,\beta}\to\mu_\infty^{D,\beta}$ along with the fact that $\bigcap_{k\ge1}\FF_k$ is trivial now identify $Y_\infty^D$ with the LQG measure on the right of \eqref{E:2.20}.

To link this to the law of $Z^D_\lambda$ we note that, as part of the proofs of Theorems~\ref{thm-3.13} and~\ref{thm-4.4}, we showed that $Z^D_\lambda\laweq Y_\infty^D$ for~$D$ being a dyadic square. The Gibbs-Markov property and the construction \eqref{E:5.25} then readily extend this to all~$D$.
\end{proofsect}

\newcounter{appendA}
\renewcommand{\theappendA}{A}
\section*{Appendix}
\addcontentsline{toc}{section}
{{\tocsection {}{8}{\!\!\!\!Appendix A: Useful properties and bounds\dotfill}}{}}
\refstepcounter{appendA}
\label{appendB}
\renewcommand{\thesection}{A}
\setcounter{subsection}{0}
\setcounter{theorem}{0}
\setcounter{equation}{0}
\renewcommand{\thesubsection}{A.\arabic{subsection}}
\renewcommand{\theequation}{\thesection.\arabic{equation}}
\noindent
Here we will review some of the needed facts concerning the \DGFF{} as well as the Green function of the simple symmetric random walk on~$\Z^2$. We begin with the latter. 

Given~$D\subsetneq\Z^2$, the Green function $G^D(x,y)$ is the expected number of visits to~$y$ of the simple random walk started at~$x$ before the walk exits~$D$.  Denoting $V_N':=(-N,N)^2\cap\Z^2$,   the potential kernel  can be defined by the limit
\begin{equation}
\fraka(x):=\lim_{N\to\infty}\bigl[G^{V_N'}(0,0)-G^{V_N'}(0,x)\bigr]
\end{equation}
The potential kernel admits the asymptotic form
\begin{equation}
\label{E:A.2}
\fraka(x)=g\log|x|+c_0+O(|x|^{-2}),\qquad |x|\to\infty.
\end{equation}
with~$c_0$ a numerical constant.
For~$D$ finite, the fact that $\fraka$ is discrete harmonic away from~$0$ while $x\mapsto G^D(x,y)$ is harmonic on~$D\smallsetminus\{y\}$, we have 
\begin{equation}
\label{E:A.3}
G^D(x,y) = -\fraka(x-y)+\sum_{z\in\partial D}H^D(x,z)\fraka(y-z),
\end{equation}
where $H^D(x,z)$ is the probability that the simple random walk started at~$x$ exits~$D$ at~$z$. As shown in \cite[Lemma~A.2]{BL2}, the class of domains~$\mathfrak D$ is such that, for any~$D\in\mathfrak D$ and any sequence~$D_N$ approximating~$D$ in the sense of \twoeqref{E:2.1}{E:2.2}, we have
\begin{equation}
\label{E:A.4}
\sum_{z\in\partial D_N}H^{D_N}\bigl(\lfloor xN\rfloor,z\bigr)\delta_{z/N}(\cdot)\,\,\overset{\text{vaguely}}{\underset{N\to\infty}\longrightarrow}\,\,\Pi^D(x,\cdot),\qquad x\in D,
\end{equation}
where $\Pi^D(x,\cdot)$ is the harmonic measure on (i.e., the hitting probability of the Brownian motion started from~$x$ to) the boundary~$\partial D$ of the continuum domain~$D$. Using this in conjunction with \eqref{E:A.2}, for $x,y\in D$ with $x\ne y$ we then get
\begin{equation}
\label{G-asymp2}
G^{D_N}\bigl(\lfloor xN\rfloor,\lfloor yN\rfloor\bigr)
= -g\log|x-y| + g\int_{\partial D}\Pi^D(x,\textd z)\log|x-z|+o(1).
\end{equation}
For $x=y$ we instead get \eqref{G-asymp} effectively replacing $-g\log|x-y|$ by $g\log N+c_0$.

Moving to the properties of the \DGFF{}, one that is most fundamental is the Gibbs-Markov decomposition. If $V\subset U\subset\Z^2$ are finite domains and $h^U$ and $h^V$ the \DGFF{}s on~$U$, resp., $V$, then
\begin{equation}
\label{E:A.6}
h^U\laweq h^V+\varphi^{U,V}
\end{equation}
with $\varphi^{U,V}$ independent of~$h^V$ and having sample paths that are discrete harmonic on~$V$ and equidistributed to $h^U$ on~$U\smallsetminus V$. The law of~$\varphi^{U,V}$ can be alternatively prescribed by its covariance structure, which turns out to be the difference $G^U-G^V$. It is now easy to check from \twoeqref{E:A.3}{E:A.4}, this difference admits a scaling limit in the sense that, for any $D,\wt D\in\mathfrak D$ with~$\wt D\subset D$ and locally uniformly in $x,y\in\wt D$
\begin{equation}
\label{E:A.7}
G^{D_N}\bigl(\lfloor xN\rfloor,\lfloor yN\rfloor\bigr)-G^{\wt D_N}\bigl(\lfloor xN\rfloor,\lfloor yN\rfloor\bigr)
\,\underset{N\to\infty}\longrightarrow\, C^{D,\wt D}(x,y),
\end{equation}
with $C^{D,\wt D}$ as in \eqref{E:2.10}. Letting $\Phi^{D,\wt D}$ be the Gaussian process with covariance $C^{D,\wt D}$, for each $N\ge1$ and each~$\delta>0$, there is a coupling of~$\varphi^{D_N,\wt D_N}$ with~$\Phi^{D,\wt D}$ such that
\begin{equation}
\label{E:A.8}
\lim_{N\to\infty}\,P\biggl(\sup_{\begin{subarray}{c}
x\in\wt D\\\dist(x,\wt D^\cc)>\delta
\end{subarray}}
\bigl|\varphi^{D_N,\wt D_N}(\lfloor xN\rfloor)-\Phi^{D,\wt D}(x)\bigr|>\delta\biggr)=0,
\end{equation}
see~\cite[Lemma~B.14]{BL3}. 

As our final item of business, we will prove a lemma that was used in the proof of Theorem~\ref{thm-4.4}. The proof is standard; we include it merely for completeness of exposition.

\begin{proofsect}{Proof of Lemma~\ref{lemma-4.5}}
To lighten the notation suppose~$S$ is an open square of side~$r$ and let $S_i$, $i=1,\dots,L^2$ be disjoint open squares of side~$r/L$ that just barely fit into~$S$. Denote $\wt S:=\bigcup_{i=1}^{L^2}S_i$ and let~$x_i$ be the center point of~$S_i$, for each~$i$. Then $\Var(\Phi^{S,\wt S}(x_i))\le g\log L+c$ for some constant~$c$ and so, by a straightforward union bound,
\begin{equation}
\label{E:A.9}
P\biggl(\,\,\max_{i=1,\dots,L^2}\Phi^{S,\wt S}(x_i)>2\sqrt g\log L\biggr)\le\frac{c'}{\sqrt{\log L}}.
\end{equation}
Next let $S_i^\delta:=\{z\in S_i\colon\dist(z,S_i^\cc)>\delta\}$ and note that 
\begin{equation}
\max_{i=1,\dots,L^2}\,\sup_{z\in S_i^\delta}\,\Var\bigl(\Phi^{S,\wt S}(z)-\Phi^{S,\wt S}(x_i)\bigr)\le c
\end{equation}
with~$c$ independent of~$L$. Letting 
\begin{equation}
M_L^\star:=\max_{i=1,\dots,L^2}\,\sup_{z\in S_i^\delta}\bigl(\Phi^{S,\wt S}(z)-\Phi^{S,\wt S}(x_i)\bigr)
\end{equation}
the Borell-Tsirelson inequality (see  Adler~\cite[Theorem~2.1]{Adler}) shows that~$M_L^\star$ has a uniform Gaussian tail and so
\begin{equation}
\label{E:A.12}
P\bigl(M_L^\star-EM_L^\star>\sqrt{\log L}\bigr)\le L^{-c}
\end{equation}
for some~$c>0$ independent of~$L$. It thus remains to control the growth rate of $E M_L^\star$. For this we consider the pseudometric space $(\scrX,\rho)$, where $\scrX:=\{(i,z)\colon i=1,\dots,L^2,\,z\in S_i^\delta\}$ and $\rho((i,z),(i',z')):=E[\Phi^{S,\wt S}(z)-\Phi^{S,\wt S}(z')\bigr]$ when $z\in S_i^\delta$ and $z'\in S_{i'}^\delta$. Writing $B_\rho((i,z),r)$ for the closed ball in~$\scrX$ of radius~$r$ centered at $(i,z)$ and using~$m$ for the normalized Lebesgue measure on~$\bigcup_{i=1}^{L^2}S_i^\delta$, the Fernique criterion (cf Adler~\cite[Theorem~4.1]{Adler}) then gives
\begin{equation}
\label{E:A.13}
EM_L^\star\le c\int_0^\infty\sqrt{\log\frac1{m(B_\rho(x,r))}}\,\textd r
\end{equation}
for some universal constant~$c$. The fact that $(x,y)\mapsto C^{S,\wt S}(x,y)$ is uniformly Lipschitz on each~$S_i^\delta$ gives $m(B_\rho((i,z),r))\ge c'(r\wedge L^{-1})^4$ with~$c'>0$ independent of~$L$ as soon as~$\delta$ is sufficiently small. Plugging this into \eqref{E:A.13}, we get $EM_L^\star\le c''\sqrt{\log L}$. Combining this with \eqref{E:A.12} and \eqref{E:A.9}, the claim follows.
\end{proofsect}

\section*{Acknowledgments}
\noindent
This research has been partially supported by European Union's Seventh Framework Programme (FP7/2007-2013]) under the grant agreement number 276923-MC--MOTIPROX, and also by NSF grant DMS-1407558 and GA\v CR project P201/16-15238S.  M.B.\ wishes to thank Jean-Dominique Deuschel, and other participants of a summer school in Hejnice, for useful remarks during presentation of an early version of this material.

\end{document}